\documentstyle[fullpage,12pt,amssymb]{amsart}
\catcode`\@=12

\def\Empty{}
\newcommand\oplabel[1]{
  \def\OpArg{#1} \ifx \OpArg\Empty {} \else
  	\label{#1}
  \fi}
		
\long\def\realfig#1#2#3#4{
\begin{figure}[htbp]
\centerline{\psfig{figure=#2,width=#4}}
\caption[#1]{#3}
\oplabel{#1}
\end{figure}}

\newcommand{\comm}[1]{}

\input{psfig}
\newtheorem{thm}{Theorem}[section]
\newtheorem{cor}[thm]{Corollary}
\newtheorem{lem}[thm]{Lemma}
\newtheorem{prop}[thm]{Proposition}

\theoremstyle{definition}

\theoremstyle{remark}

\newtheorem{ex}{Example}[section]

\newcommand{\diam}{\operatorname{diam}}
\newcommand{\area}{\operatorname{area}}
\newcommand{\dist}{\operatorname{dist}}
\newcommand{\cl}{\operatorname{cl}}
\newcommand{\orb}{\operatorname{orb}}
\newcommand{\myint}{\operatorname{int}}
\newcommand{\mymod}{\operatorname{mod}}
\newcommand{\geo}{\operatorname{geo}}
\newcommand{\ltime}{\operatorname{ltime}}

\newcommand{\sgn}{\operatorname{sgn}}
\newcommand{\depth}{\operatorname{depth}}
\newcommand{\wtl}{\widetilde}

\newcommand{\upto}[1]{\stackrel{{#1}}{\to}}
\newcommand{\eps}{\epsilon}
\newcommand{\es}{\emptyset}
\newcommand{\setm}{\setminus}
\newcommand{\bd}{\partial}
\newcommand{\union}{\cup}
\newcommand{\seq}[1]{{\langle}{#1}{\rangle}}

\newcommand{\cF}{{\cal F}}
\newcommand{\cB}{{\cal B}}
\newcommand{\cT}{{\cal T}}
\newcommand{\cN}{{\cal N}}

\newcommand{\cC}{{\cal C}}
\newcommand{\cH}{{\cal H}}
\newcommand{\cR}{{\cal R}}

\newcommand{\cD}{{\cal D}}
\newcommand{\cE}{{\cal E}}

\newcommand{\cK}{{\cal K}}
\newcommand{\Hol}{{\cal P}}
\newcommand{\Tow}{Tow}
\newcommand{\inter}{\cap}
\newcommand{\C}{{\Bbb C}}
\newcommand{\R}{{\Bbb R}}
\newcommand{\N}{{\Bbb N}}
\newcommand{\Hp}{{\Bbb H}}
\newcommand{\D}{{\Bbb D}}

\newcommand{\Z}{{\Bbb Z}}
\newcommand{\Chat}{\widehat{{\Bbb C}}}
\newcommand{\qtr}{1/4}

\newcommand{\bff}{{\bold f}}
\newcommand{\bfg}{{\bold g}}
\newcommand{\bfh}{{\bold h}}
\newcommand{\bR}{{\bold R}}
\newcommand{\bT}{{\bold T}}

\newcommand{\bX}{{\Bbb X}}
\newcommand{\bV}{{\Bbb V}}
\newcommand{\bU}{{\Bbb U}}

\numberwithin{equation}{section}
\newcommand{\thmref}[1]{Theorem~\ref{#1}}
\newcommand{\propref}[1]{Proposition~\ref{#1}}
\newcommand{\secref}[1]{\S\ref{#1}}
\newcommand{\lemref}[1]{Lemma~\ref{#1}}
\newcommand{\corref}[1]{Corollary~\ref{#1}} 
\newcommand{\figref}[1]{Fig.~\ref{#1}}
\newcommand{\exref}[1]{Example~\ref{#1}}

\def\SBIMSMark#1#2#3{
 \font\SBF=cmss10 at 10 true pt
 \font\SBI=cmssi10 at 10 true pt
 \setbox0=\hbox{\SBF Stony Brook IMS Preprint \##1}
 \setbox2=\hbox to \wd0{\hfil \SBI #2}
 \setbox4=\hbox to \wd0{\hfil \SBI #3}
 \setbox6=\hbox to \wd0{\hss
             \vbox{\hsize=\wd0 \parskip=0pt \baselineskip=10 true pt
                   \copy0 \break%
                   \copy2 \break% 
                   \copy4 \break}}
 \dimen0=\ht6   \advance\dimen0 by \vsize \advance\dimen0 by 8 true pt
                \advance\dimen0 by -\pagetotal
 \dimen2=\hsize \advance\dimen2 by .25 true in
  \openin2=publishd.tex
  \ifeof2\setbox0=\hbox to 0pt{}
  \else 
     \setbox0=\hbox to 3.1 true in{
                \vbox to \ht6{\hsize=3 true in \parskip=0pt  \noindent  
                \input publishd.tex 
                \vfill}}
  \fi
  \closein2
  \ht0=0pt \dp0=0pt
 \ht6=0pt \dp6=0pt
 \setbox8=\vbox to \dimen0{\vfill \hbox to \dimen2{\copy0 \hss \copy6}}
 \ht8=0pt \dp8=0pt \wd8=0pt
 \copy8
 \message{*** Stony Brook IMS Preprint #1, #2 ***}
}

\begin{document}

\title[Renormalization ]{Parabolic Limits of Renormalization}
\author { Benjamin Hinkle }

\date{July 10, 1997}

\maketitle
\SBIMSMark{1997/7}{July 1997}{}

\begin{abstract}
In this paper we give a combinatorial description 
of the renormlization limits of infinitely renormalizable unimodal maps
with {\it essentially bounded} combinatorics admitting quadratic-like
complex extensions. 
As an application we construct a natural analogue of the period-doubling
fixed point.
Dynamical hairiness is also proven for maps in this class.
These results are proven by analyzing {\it parabolic towers}: sequences
of maps related either by renormalization or 
by {\it parabolic renormalization}.
\end{abstract}

\tableofcontents

\section{Introduction}
\label{sec:introduction}

In this paper we extend the well-known combinatorial description 
of renormalization limits of unimodal maps with bounded
combinatorics to renormlization limits of maps
with {\it essentially bounded} combinatorics. 
This class of maps was introduced
by Lyubich \cite{L2} and their complex geometry was studied in \cite{LY}.
Roughly speaking the high renormalization periods of such maps
are due to their renormalizations being small perturbations of parabolic
maps. Although this leads to the creation of unbounded combinatorics,
the {\it essential geometry} of these maps remains bounded away from zero.

Let us state our results (see \secref{sec:back1} for background). 
In \secref{sec:special} we construct a countable collection of 
maximal tuned Mandelbrot copies $\{M^{(3)}_n\}_{n=1}^\infty$ 
that accumulate at $c = -1.75$, the root point of $M^{(3)}$.
These copies have ``essentially period tripling'' combinatorics.
Our first result produces the
analog of the renormalization fixed point in the 
essentially period tripling situation:

\begin{thm}
\label{per3}
There is a unique quadratic-like germ $F$ such that
$$
\cR^n(f) \to F
$$
for any quadratic-like map $f$ in the hybrid class of an
infinitely renormalizable real quadratic with a tuning invariant
$$
\tau(f) = (M^{(3)}_{n_1},M^{(3)}_{n_2},\dots,M^{(3)}_{n_k},\dots)
$$
satisfying $n_k \to \infty$ as $k \to \infty$.
Any quadratic-like representative of $F$ is hybrid equivalent to
$z^2 - 1.75$ and hence has a period three parabolic orbit.
\end{thm}

In order to state our second theorem we need to fix some notation.
Let $\Omega$ denote the space of 
{\it unimodal non-renormalizable permutations}, or {\it shuffles},
and let $p_e(\sigma)$ be the essential period of $\sigma \in \Omega$.
Let 
$$
\Omega_p = \{\sigma \in \Omega: p_e(\sigma) \le p\}
$$ 
and let $\Omega^{cpt}_p$ be the compactification of $\Omega_p$
defined in \secref{sec:shuffles}.
Let 
$$
\Sigma_p = \Pi_{-\infty}^{\infty} \Omega^{cpt}_p
$$
with coordinate projections $\pi_n:\Sigma_p \to \Omega^{cpt}_p$
and let $\omega:\Sigma_p \to \Sigma_p$ be the left shift operator.
Let $\Omega^{cpt}_{p,\ast}$ denote the space $\Omega^{cpt}_p$
with the symbol $\ast$ adjoined and let
$\Sigma^{\ast}_p = \Pi_{-\infty}^{\infty} \Omega^{cpt}_{p,\ast}$.
We will denote the left shift on $\Sigma^{\ast}_p$ by $\omega$ as well.
For any quadratic-like map $f$ hybrid
equivalent to an infinitely renormalizable real polynomial
let 
$$
\bar{\sigma}(f) = (\dots,\ast,\ast,\sigma_0,\sigma_1,\sigma_2,\dots)
$$ 
where $\sigma_n$ is the shuffle corresponding to the 
$n$-th Mandelbrot copy in $\tau(f)$.
Let $\bar{p}_e(f) = \sup_{n \ge 0} p_e(\pi_n(\bar{\sigma}(f)))$.
Let $GQuad(m)$ be the space of quadratic-like
germs with modulus at least $m$.
We can now state the combinatorial classification of all limits
of renormalization of an infinitely renormalizable real quadratic with
essentially bounded combinatorics.

\begin{thm}
\label{generalThm}
There is an $m > 0$ so that for any $p > 1$ there exists a
continuous map 
$$
h: \Sigma_p \to GQuad(m)
$$ 
with the following property.
Let $f$ be a quadratic-like map
in the hybrid class of an $\infty$-renormalizable real quadratic with 
$\bar{p}_e(f) \le p$ and
let ${\bar \sigma} \in \Sigma_p$ be a limit point of
$\bar{\sigma}_n = \omega^n(\bar{\sigma}(f))$.
If $\bar{\sigma}_{n_i} \to \bar{\sigma}$ then
\begin{equation}
\label{eqn:convergeR}
\cR^{n_i}(f) \to h({\bar \sigma}).
\end{equation}
Furthermore, if $\sigma_0 = \pi_0(\bar{\sigma}) \in \Omega_p$ then 
$h(\bar{\sigma})$ is renormalizable by the shuffle type 
$\sigma_0$ and $h$ is a conjugacy between $\omega$ and
$\cR_{\sigma_0}$.
If $\sigma_0 \not\in \Omega_p$ then 
the inner class of $h(\bar{\sigma})$ is the root of a 
maximal tuned Mandelbrot copy $M(\sigma_0)$.
\end{thm}

Let us comment on the ideas involved in the paper. Recall that 
the central objects of McMullen's argument \cite{McM2} are 
{\it towers:} sequences of quadratic-like maps related by renormalization.
A {\it forward tower} is a one-sided sequence and a
{\it bi-infinite tower} is a two-sided infinite sequence.
The question of convergence of renormalization
is equivalent to the question of
{\it combinatorial rigidity} of the corresponding 
limiting bi-infinite towers.
However, for maps with essentially bounded combinatorics the limiting
towers may contain parabolic maps and we lose the
renormalization relation between levels.
In this case a new relation appears: {\it parabolic renormalization}.
That is, the maps in the limiting towers are related by
either classical or parabolic renormalization. A tower which contains
a parabolic renormalization is called a {\it parabolic tower}.
Our proof of the rigidity of bi-infinite parabolic towers with
definite modulus and essentially bounded combinatorics
consists of first analyzing forward towers
and then analyzing bi-infinite towers. 

Our analysis of forward parabolic towers was motivated by the work of
A. Epstein \cite{E}, which considered general
holomorphic dynamical systems (with maximal 
domains of definition) and their geometric limits.
The phenomenon studied there was the 
renormalization  (different from the sense used in this paper)
of a parabolic orbit at the ends of its
\`Ecalle-Voronin cylinders. The phenomenon we study occurs away from the ends
and as a result the forward infinite
towers in this paper look in many ways
like infinitely renormalizable real quadratic maps.

The combinatorial rigidity of forward parabolic towers 
with polynomial base map follows from
the theory of {\it quadratic-like families} and from
the combinatorial rigidity of quadratic polynomials with
complex bounds and real combinatorics (see \propref{genCombRigid}).
After analyzing the Julia set of a foward tower we prove any
quasiconformal conjugacy of a forward infinite parabolic tower 
with essentially bounded combinatorics and complex bounds
is a hybrid conjugacy (see \secref{sec:forrigid}).

Then following the arguments
of McMullen we prove in \secref{sec:biTowers}
the rigidity of bi-infinite towers. That is, we first prove

\begin{thm}[Dynamical Hairiness]
\label{dense-lemma}
The union of the Julia sets of the forward infinite sub-towers
of a bi-infinite tower with essentially bounded combinatorics and
complex bounds is dense in the plane.
\end{thm}

Then we prove

\begin{thm}
\label{affine}
Any quasiconformal equivalence
of a bi-infinite tower with essentially bounded combinatorics
and complex bounds is affine.
\end{thm}

Let us mention a parallel with critical circle maps.
The theory of renormalization of unimodal maps is closely related to
renormalization theory of critical circle maps. The rotation number 
$\rho$, more specifically its continued fraction expansion, determines
the combinatorics of a circle map.
If the factors in its expansion
are bounded then the map has bounded combinatorics 
and has unbounded combinatorics otherwise. If a circle
map has unbounded combinatorics then the rotation numbers of the 
renormalizations contain rational limit points and the 
corresponding limit of renormalization contain parabolic periodic points.
That is, the only kind of unbounded combinatorics in the theory
of critical circle maps is the essentially bounded combinatorics.
DeFaria \cite{deF} analyzed the renormalization limits of critical circle 
maps with bounded combinatorics and Yampolsky \cite{Y} proves
complex bounds for arbitrary combinatorics.
We expect the techniques in this paper can be adapted to analyze 
renormalization limits of critical circle maps with
arbitrary combinatorics.

The author specially thanks 
Misha Lyubich for his suggestions and guidance,
Adam Epstein and Misha Yampolsky for the many useful conversations, 
and UNAM at Cuernavaca for their gracious hospitality.

\subsection{Notation}

\begin{itemize}
\item $\Hp \subset \C$ denotes the complex upper half-plane,
	$\Chat$ the Riemann sphere, $\N = \N_0$ the non-negative integers and
	$\N_+$ the positive integers.
\item $[a,b]$ will also denote the interval $[b,a]$ if $b < a$.
\item $\diam(U)$ denotes the euclidean diameter of $U \subset \C$
	and $|I|$ the diameter of $I \subset \R$.
\item $\cl(X)$, $\myint(X)$ and $\bd X$
	denote the closure, interior and boundary of $X$ in $\R$
	if $X \subset \R$ and in $\C$ otherwise.
\item $U \Subset V$ means $U$ is compactly contained in $V$. Namely 
	$\cl(U)$ is compact and $\cl(U) \subset V$.
\item in a dynamical context $f^n$ denotes $f$ composed with
	itself $n$ times.
\item if $V$ is a simply connected domain and $U \subset V$ then
$\mymod(U,V) = \sup_A \mymod(A)$ where $A$ is an annulus separating
$U$ from $\bd V$.
\item $Dom(f)$ and $Range(f)$ denote the domain and range of $f$.
\item $Comp(X)$ denotes the collection of connected components of $X$
and $Comp(X,Y)$ denotes the components of $X$ intersecting $Y$.
\item $P_c(z) = z^2 + c$.
\end{itemize}

\section{Background}
\label{sec:back1}
\subsection{Quadratic-like maps}

We will assume the reader is familiar with the theory of quasiconformal maps
and the Measurable Riemann Mapping Theorem (see \cite{LV}). 

A holomorphic map $f:U \to V$ is {\it quadratic-like} if 
$U$ and $V$ are topological disks in $\C$ with $U \Subset V$
and $f$ is a branched double cover of $U$ onto $V$.
By topological disk we mean a simply connected domain in $\C$.
A topological disk whose boundary is a Jordan curve will be called a
{\it Jordan disk}. Unless otherwise indicated we will assume the
critical point of a quadratic-like map is at the origin.
A point $z \in U$ is {\it non-escaping} if
$f^n(z)$ is defined for all $n \ge 0$.
For a quadratic-like $f:U \to V$ define
\begin{itemize} 
	\item The {\it filled Julia set}
		$K(f) = \cl\{z \in U : z \mbox{ is non-escaping}\}$
	\item The {\it Julia set}
		$J(f) = \bd K(f)$
	\item The {\it post-critical set}
		$P(f) = \cl\{\bigcup_{n \ge 1} f^n(0)\}$
\end{itemize}

An actual quadratic polynomial can be considered quadratic-like by taking
$V=\{z : |z| < R\}$ for some large $R$. 
Following \cite{McM2}, define
$Quad$ to be the union of all quadratic-like maps $f:U \to V$
and all quadratic polynomials $f:\C \to \C$
with a non-escaping critical point at the origin.

Impose on {\it Quad} the {\it Carath\'eodory topology}.
That is, a sequence $f_n:U_n
\to V_n$ converges to $f:U \to V$ iff 
$(U_n,0)$ and $(V_n,f_n(0))$ converge to $(U,0)$ and $(V,f(0))$, respectively,
in the Carath\'eodory topology on pointed domains
in the Riemann sphere $\Chat$
and $f_n$ converges uniformly to $f$ on compact subsets of $U$.
We note the facts:
\begin{enumerate}
\item For any compact connected $U \subset V$
if $\mymod(U,V) \ge m$ then $V$ contains an $\eps(m)$-scaled
neighborhood of $U$, where
an $\eps$-scaled neighborhood of a domain $U$ is an $\eps \cdot \diam(U)$
neighborhood of $U$
\item If the domains are $K$-quasidisks then the Carath\'eodory convergence of pointed
domains is equivalent to the Hausdorff convergence of their closures.
\item The set of $K$-quasidisks in $\C$ containing a definite neighborhood 
of the origin
and with bounded diameter is compact in the Hausdorff topology.
\item If a sequence of pointed domains $(U_n,u_n)$ in $\Chat$
converges in the Carath\'eodory topology to the domain 
$(U,u)$, and if $U_n$ and $U$ are all hyperbolic Riemann surfaces,
then the hyperbolic metrics on $U_n$ converge in the $C^\infty$ norm
uniformly on compact set of $U$ to the hyperbolic metric on $U$.
\end{enumerate}
Define the subspaces
$$
Quad(m) = \{f \in Quad : f \mbox{ is a polynomial or }
\mymod(U,V) \ge m \}
$$
and
$$
RQuad = \{f \in Quad : f(\bar z)= {\overline {f(z)}}\}.
$$
The following compactness lemma is a basic tool in renormalization theory:
\begin{lem}[{\cite[Theorem 5.8]{McM1}}]
\label{Compact}
For any $C_0 > 0$, $C_1 < \infty$, $m > 0$, the set
$$
\{f \in Quad(m) : C_0 \le \diam K(f) \le C_1\}
$$
is compact.
\end{lem}

Define $GQuad(m)$ to be the quotient space of $Quad(m)$ by the relation
$f \sim g$ iff $f = g$ on a neighborhood of zero.
Define the set of {\it quadratic-like germs} to be
$GQuad = \union_m GQuad(m)$. Convergence of germs will always take place
in some $GQuad(m)$. The germ of $f$ will be denoted by $[f]$.
From \cite[Lemma 7.1]{McM2} the (filled) Julia set of a quadratic-like germ
is well defined and consequently if $f$ and $g$ 
are two quadratic-like representatives of a germ in $GQuad(m)$ then
$f = g$ on an $\eps(m/2)$-scaled neighborhood of $K(f) = K(g)$.
Since $K(f)$ is an upper semi-continuous function on $Quad$, 
if $f_k \in GQuad(m)$ converges to $f$ 
then for any sequence of representatives $g_k \in Quad(m)$ it follows
$g_k$ converges to $g$ on a definite neighborhood of $K(f)$.
Given $f \in GQuad$ let 
$$
\mymod(f) = \sup \mymod(U,V)
$$ 
where the supremum
is taken over all quadratic-like representations of $f$.

Let $f \in Quad$.
For a given $x \neq 0$ let $x' = f^{-1}(f(x)) \setm \{x\}$.
If $x = 0$ let $x' = 0$. There are two fixed points $\alpha$
and $\beta$ of $f$ counted with multiplicity and labeled
so that $J(f) \setm \{\beta\}$ is connected.
The only case when $\alpha = \beta$ is when $I(f) = \qtr$.
We say $f \in Quad$ is {\it normalized} if $\beta(f) = 1$.
We normalize a germ by normalizing any quadratic-like representative.

A {\it quasi-conformal equivalence} $\phi$ between quadratic-like maps $f$ and
$g$ is a quasiconformal
map from a neighborhood of $K(f)$ to a neighborhood of $K(g)$
such that $\phi \circ f = g \circ \phi$.
A quasi-conformal equivalence is a {\it hybrid} equivalence if 
$\bar{\partial} \phi|_{K(f)} = 0$ as a distribution.

\begin{prop}[Straightening,\cite{DH2}]
\label{Straighten}
Any quadratic-like map $f$ is hybrid equivalent to a quadratic polynomial.
If $K(f)$ is connected the polynomial is unique up to affine conjugacy.
Moreover, if $f \in Quad(m)$ then the equivalence can be chosen
to be a conjugacy on an $\epsilon(m)$-scaled neighborhood of $K(f)$ and 
with dilatation bounded above by $K(m) < \infty$.
\end{prop}

The {\it inner class} of a map $f \in Quad$, denoted $I(f)$,
is the unique $c$ value such that $f$ is hybrid equivalent to $P_c$. 
The inner class of a germ $I([f])$ is the inner class of any
quadratic-like representative.
The {\it Mandelbrot set}, $M$, is the set of $c \in \C$ such that the Julia
set of $z^2+c$ is connected.
Let $\cH(c) = I^{-1}(c)$ for $c \in M$. Note 
$\cH(c) \subset Quad$.

\begin{prop}[\cite{DH2},{\cite[Proposition 4.7]{McM2}}]
\label{StraightenCont}
$I:Quad \to M$ is continuous.
\end{prop}

\subsection{Renormalization}

A parameter value $c \in \C$ is called 
{\it super-stable} if $0$ is periodic under $P_c$.
To each super-stable $c \neq 0$ there is associated a 
homeomorphic copy of $M$
containing $c$ called the {\it Mandelbrot set 
tuned by $c$}, or, briefly, an {\it $M$-copy}, and denoted by $c \star M$.
The {\it root} of $c \star M$ is the point corresponding to $\qtr$ 
and the {\it center} is the point $c$.
For every copy $c \star M$ there is a $p > 1$
such that for any $c' \in c \star M$,
except possibly the root, and any $f \in \cH(c')$
there is a domain $U \ni 0$ such that $f^p|_U \in Quad$.
The map $f^p|_U$ is called a 
(complex) {\it pre-renormalization} of $f$ and
$f$ is said to be {\it renormalizable} of period $p$.
This pre-renormalization is always {\it simple}, meaning the iterates of
$J(f^p|_U)$ under $f$ are either are disjoint or intersect only along
the orbit of $\beta(f^p|_U)$.
The {\it period} of the copy, $p(c \star M)$, is the maximal such $p$
and we say $c \star M$ is {\it maximal} if there is only one such $p$.
We say $c \star M$ is {\it real} if $c$ is real.
The only real maximal $M$-copy for which the root point is not
renormalizable is the period two copy $M^{(2)}$.
We will denote the real period three copy by $M^{(3)}$.
Define $\cH(c \star M)$ to be the set of 
renormalizable $f \in I^{-1}(c \star M)$. 
In \figref{Mandel} we have drawn the Mandelbrot set 
highlighting $M^{(2)}$ and $M^{(3)}$. The root points are 
$c = -0.75$ and $c = -1.75$, respectively.

\realfig{Mandel}{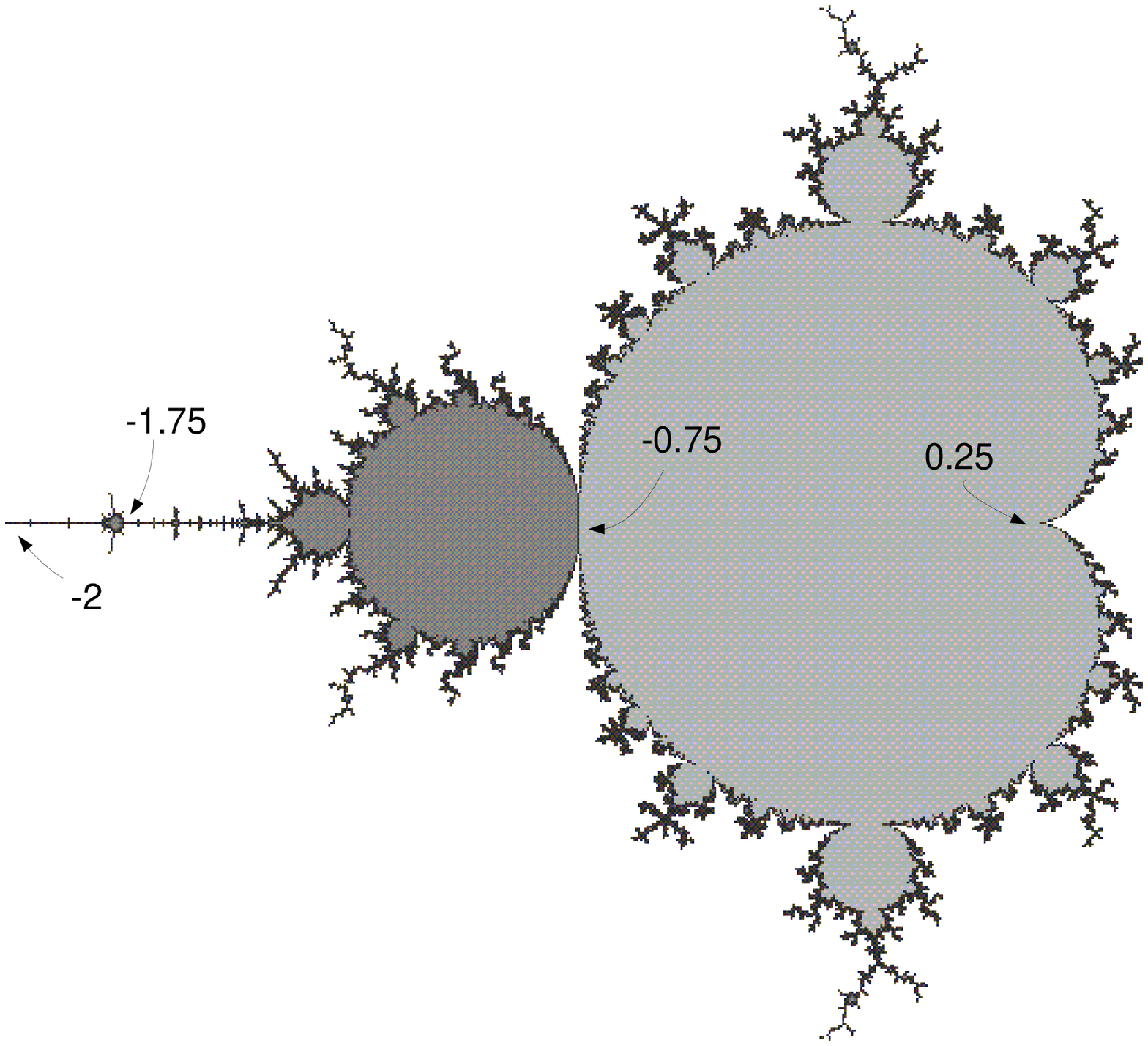}{The Mandelbrot set.}{0.6\hsize}

Let $c \star M$ be a maximal $M$-copy with period $p$ and
suppose $f \in \cH(c \star M)$.
If $f^p|_U$ and $f^p|_{U'}$ are two pre-renormalizations then
$[f^p|_U] = [f^p|_{U'}]$. Hence we can
define the {\it renormalization} $\cR(f)$ to be the
normalized quadratic-like germ of any pre-renormalization of period $p$.
We define the renormalization of a germ $\cR([f])$
to be the renormalization of a quadratic-like representative.
A map $f \in Quad$ is {\it infinitely renormalizable} if $\cR^n(f)$ is 
defined for all $n \ge 0$, or, equivalently, if 
$I(f)$ is contained in infinitely many $M$-copies.
The tuning invariant of an infinitely renormalizable map $f \in Quad$ is 
$$
\tau(f) = (M_0,M_1,M_2,\dots)
$$ 
where $M_n$ is the maximal $M$-copy containing $I(\cR^n(f))$. 
We say $f$, even if it is only finitely renormalizable,
has {\it real combinatorics} if all $M$-copies in $\tau(f)$
are real. 
See \cite{D2} for a more complete description of tuning.

Let us turn to renormalization in real dynamics.
Let $I \subset \R$ be a closed interval.
A continuous map $f:I \to I$ is {\it unimodal} if 
$f(\bd I) \subset \bd I$ and there is a unique extremum $c$ of $f|_I$.
For $f \in RQuad$ let $B(f)=[\beta,\beta']$, and
$A(f)=[\alpha, \alpha'] \subset B(f)$.
Note that $K(f) \inter \R = B(f)$. 
The next lemma follows from \lemref{Compact}
and the continuity of $\beta$ and $\beta'$.
We say $A$ is {\it $C$-commensurable} to $B$
if $C^{-1} \le A/B \le C$.

\begin{lem}
\label{DiamCompare}
For $m > 0$, $|\beta(f)|$ and $|B(f)|$ are $C(m)$-commensurable to
$\diam K(f)$  for any $f \in RQuad(m)$.
\end{lem}

Any $f \in RQuad$ is unimodal on $B(f)$ (and this interval is maximal).
We say a unimodal map $f|_{[a,b]}$ with $a < b$
is {\it positively oriented} if $f(b) = b$. 
The quadratic family $P_c$ is positively oriented.
A unimodal map $f:I \to I$ is {\it real-renormalizable} if there is an 
interval $I' \ni c$ 
and an $n > 1$ such that $f^n|_{I'}$ is unimodal.
%Any complex-renormalizable $f \in RQuad$ is real-renormalizable 
%with the same period and
%any real-renormalizable $f \in RQuad \setm \cH(-0.75)$ is 
%complex-renormalizable with the same period if it is minimal.
%However, any $f \in RQuad \inter \cH(-0.75)$ is real-renormalizable
%but not complex-renormalizable, since $f'(\alpha) = -1$.
Unlike complex renormalization, we can canonically define
real-renormalization as acting on unimodal maps as follows.
Define the real pre-renormalization $f_1$ of a unimodal map $f$ as
$f^n|_{I'}$ where $n$ is minimal and $I$ is maximal and define
the real-renormalization $\cR(f)$ as $f_1$ conjugated by
$x \mapsto x/\beta(f_1)$ where $\beta(f)$ is the boundary fixed point
of $f$.

Suppose $f \in RQuad$ is real-renormalizable and positively oriented.
Let $f_1$ be a pre-renormalization and
let $\sigma(f)$ be the permutation
induced on the orbit of $B(f_1)$ labeled from left to right. 
Any permutation that can be so realized is called a 
{\it unimodal non-renormalizable permutation}, or a {\it shuffle}.
The permutation on two symbols we will denote by $\sigma^{(2)}$.
If $\sigma(f) = \sigma^{(2)}$ we say $f$ is {\it immediately renormalizable}.
The map $c \star M \mapsto \sigma(P_c)$ from the set of real maximal
$M$-copies to the set of shuffles is a bijection. We will denote
the shuffle corresponding to $c \star M$ by $\sigma(c \star M)$
and the real maximal $M$-copy 
corresponding to $\sigma$ by $M(\sigma)$.
We will occasionally use the notation
$\cR_{\sigma}$ to denote the complex renormalization operator acting
on $\cH(M(\sigma))$ and on its germs.
If $g \in \cH(M(\sigma))$ then define $\sigma(g) = \sigma$.
For an infinitely renormalizable $f \in Quad$ with real
combinatorics define
$$
\bar{\sigma}(f) = (\sigma(M_0),\sigma(M_1),\sigma(M_2),\dots)
$$ 
where $\tau(f) = (M_0,M_1,M_2,\dots)$.

An $\infty$-renormalizable map $f \in Quad$ has
{\it complex bounds} if there is some $m > 0$ such that 
the domain $U_k$ and range $V_k$ of the $k$-th 
complex pre-renormalization $f_k$
can be chosen to satisfy $\mymod(U_k,V_k) \ge m$ for all $k \ge 1$.
The following theorem establishes combinatorial rigidity 
of infinitely renormalizable maps with real combinatorics and
complex bounds.

\begin{thm}[\cite{L3}]
\label{hyper-dense}
If $P_c$ and $P_{c'}$ are two $\infty$-renormalizable
quadratics with complex bounds and the same real combinatorics
%$\bar{\sigma}(P_c) = \bar{\sigma}(P_{c'})$
then $c = c'$.
\end{thm}

Complex bounds are proven to exist for real quadratics:

\begin{thm}[\cite{LY,L2,S,LS}]
\label{complex-bounds}
Real infinitely renormalizable quadratics have complex bounds.
That is, if $f_k:U_k \to V_k$ is a complex renormalization of an
$\infty$-renormalizable real quadratic $f$ then $[f_k] \in GQuad(m)$
for some $m > 0$ independent of $k$.
Moreover, $U_k$ and $V_k$ can be chosen to be $K$-quasidisks,
$$
\diam(V_k) \le C \cdot |B(f_k)|,
$$
and, if $\sigma(\cR^{k-1}(f)) \neq \sigma^{(2)}$
then the unbranched condition holds: 
$$
P(f) \inter V_k = P(f_k).
$$
The values $m$, $C$ and $K$ are independent of $f$.
\end{thm}

When we make an additional assumption on the combinatorics we obtain
the unbranched condition on all levels.

\begin{lem}
\label{unbranched}
Let $\eps >0$.
Suppose $f$ is an infinitely renormalizable real quadratic
with $I(\cR^k(f)) \ge -2+\eps$ for all $k \ge 0$.
Then there is an $m > 0$ such that 
the domain $U_k$ and range $V_k$ 
of the $k$-th pre-renormalization can be chosen to 
satisfy 
\begin{itemize}
\item $\mymod(U_k,V_k) \ge m$
\item $U_k$ and $V_k$ are $K$-quasidisks
\item $\diam(V_k) \le C \cdot |B(f_k)|$
\item $P(f) \inter V_k = P(f_k)$
\end{itemize}
for all $k \ge 1$.
The constants $m$ and $K$ then depend on $\eps$.
\end{lem}
\begin{pf}
If $\sigma(\cR^{k-1}(f)) \neq \sigma^{(2)}$
then let $U_k$ and $V_k$ be from \thmref{complex-bounds}.
So assume $\sigma(\cR^{k-1}(f)) = \sigma^{(2)}$.
Let $h:U_{k-1}' \to V_{k-1}'$ 
and $h_1:U_k' \to V_k'$ be the $(k-1)$-st 
and $k$-th pre-normalization from \thmref{complex-bounds}
rescaled so that $\diam(K(h)) = 1$.
Let $E = P(h) \setm P(h_1)$.
From the following lemma, \propref{StraightenCont} and
the assumption $I(\cR^k(f)) \ge -2+\eps$ 
we obtain
$$
\dist(E,B(h_1)) = |h^3(0)-\alpha(h)|
	\ge C(\eps,m) > 0.
$$
From a construction of Sands, $V_k'$ can be chosen to be the union
of a euclidean disk centered at $0$ of radius
$|\beta(h_1)|$ and two small euclidean
disks centered at $\pm \beta(h_1)$ of radius $\eps' > 0$. The modulus
$\mymod(U_k',V_k')$ is bounded below by a function $m'(\eps') > 0$.
Choose $\eps' < C(\eps,m)$.
\end{pf}

\subsection{Generalized quadratic-like maps}
A holomorphic map $f$ is {\it generalized \linebreak quadratic-like} if 
$Range(f) = V$ is a topological disk,
each $U \in Comp(Dom(f))$ is 
a topological disk compactly contained in $V$
and $f|_U$ is a conformal isomorphism except for 
a distinguished component $U_0$, the {\it central component}, where
$f|_{U_0}$ is a branched double cover onto $V$.
We will consider only generalized quadratic-like maps whose domain
has of finitely many components.
Define the {\it filled Julia set}, $K(f)$,
the {\it Julia set}, $J(f)$, and the {\it post-critical set},
$P(f)$, as for quadratic-like maps.

Let $Gen$ be the union of $Quad$
and the space of generalized quadratic-like maps
with a non-escaping critical point at the origin.
Let $RGen$ be the space of real-symmetric maps in $Gen$
with real-symmetric domains.
Define 
$$
Gen(m) = \{f \in Gen : \mymod(Dom(f),Range(f)) \ge m \}.
$$ 

Impose on $Gen$ the Carath\'eodory topology as follows.
For a given $f \in Gen$ let $f(0)$ be the basepoint of $Range(f)$ and let
$u_f = f^{-1}(f(0))$ be the basepoints of $Comp(Dom(f))$.
A sequence $f_n \in Gen$ converges to $f$ iff 
\begin{itemize}
\item $u_{f_n}$ converges in the Hausdorff topology to $u_f$ 
\item if $X$ is any Hausdorff limit point of $\Chat \setm Dom(f_n)$ then 
	$$Dom(f) = Comp(\Chat \setm X, u_f)$$
\item $f_n \to f$ on compact subsets of $Dom(f)$.
\end{itemize}
The space of generalized quadratic-like germs 
is the quotient space of $Gen$
by the relation $f \sim g$ iff 
$f = g$ on a neighborhood of $u_f = u_g$.

Define the {\it geometry} of $f \in Gen$ as
$$
\geo(f) = \inf_{U \in Comp(Dom(f))} \frac{\diam(K(f) \inter U)}{\diam K(f)}.
$$
The following lemma is a
direct generalization of \lemref{Compact}.
\begin{lem}
\label{genCompact}
For a given $m > 0$, $C_0 > 0$ and $C_1$ the set
$$
\{f \in Gen(m): \geo(f) \ge C_0 \mbox{ and }
C_0 \le \diam K(f) \le C_1 \}
$$
is compact.
\end{lem}

Suppose $f:\union U_j \to V$ is a generalized quadratic-like map
with critical point at the origin and suppose $U \subset \C$ is open.
Define the open sets $D_0$ and $D_+$ by
$$
D_{0/+} = \{z : f^n(z) \in U \mbox{ for some } n \in \N_{0/+}\}
$$
and the maps $L_{0/+}(f,U) = L_{0/+}:D_{0/+} \to U$ by
$L_{0/+}(z) = f^n(z)$ for the minimal {\it landing time} 
$\ltime(z) = n \in \N_{0/+}$
such that $f^n(z) \in U$.
We call $L_0$ the {\it first landing map} and $L_+$
the {\it strict first landing map}.

Define the {\it first return map to $U$}, $R(f,U) = R$ by 
$$
R = L_+|_{D_+ \inter U}.
$$
Let $R = R(f,U_0)$ and suppose
$0 \in Dom(R)$. Define the {\it generalized renormalization} of $f$ as
$R$ restricted to $Comp(Dom(R),P(f) \union \{0\})$.
\begin{lem}
\label{returngeo}
Let $\lambda > 0$, $m > 0$ and $r \in \N$.
Suppose $f \in Gen(m)$ satisfies $\geo(f) \ge \lambda$
and suppose $g \in Gen$ is a restriction of $R(f,U_0)$ such that 
$$
Dom(g) = Comp(Dom(R),u_g)
\mbox{ and } \sup_{z \in Dom(g)} \ltime(z) \le r.
$$
Then there exists $C(\lambda,m,r) > 0$ such that $\geo(g) \ge C$.
\end{lem}
\begin{pf}
Assume $\diam K(f) = 1$.
Let $\union_j U_j = Dom(f)$ and let $K_j = U_j \inter K(f)$.
Since $\diam K_j \ge \lambda$ and $\mymod(K_j,U_j) \ge m$
it follows that $U_j$ contains an $\eps(\lambda,m)$ neighborhood
of $K_j$. For each $z \in u_g$ let 
$U_{z,0},U_{z,1},\dots,U_{z,k} = U_0$ be the pull back of $U_0$ along
the orbit $z, f(z),\dots,f^k(z) = g(z)$. 
Assume $f(z) \in \union_{j \neq 0} U_j$. From the Koebe Distortion
Theorem and the fact that $k \le r$
it follows $U_{z,1}$ contains a definite neighborhood of $f(z)$
and that $f^{k-1}:U_{z,1} \to U_0$ has bounded distortion.
Hence each $U_{z,0}$ contains a definite neighborhood of $z$
by \lemref{genCompact}.
The lemma follows by
pulling $\union_{z \in u_g} U_{z,0}$ back to each $U_{z,1}$ by
a map with bounded distortion and bounded derivative
and then to $U_{z,0}$.
\end{pf}

Let $L_0 = L_0(f,U)$ and suppose $0 \in Dom(L_0)$.
Define the {\it first through map}, $T = T(f,U)$, of $f$ by
$$
T = f \circ L_0.
$$
We shall analyze the geometry of certain first through maps
in \secref{sec:ParabolicBackground}.

\subsection{Families of generalized quadratic-like maps}
\label{sec:family}

In this section we summarize the
theory of holomorphic families of generalized quadratic-like maps.
For further details see \cite{L4}.
Let $D \subset \C$ be a Jordan disk and fix $\ast \in D$.
Let $\pi_1$ and $\pi_2$ be the coordinate projections of $\C^2$
to the first and second coordinates.
Given a set $\bX \subset \C^2$ let
$X_\lambda = \pi_2(\bX \inter \pi_1^{-1}(\lambda))$.
An open set $\bX \subset \C^2$ is a {\it Jordan bidisk} over $D$ if 
$\pi_1(\bX) = D$ and
$X_\lambda$ is a Jordan disk for all $\lambda \in D$.
We say $\bX$ admits an extension to the
boundary if $\cl(\bX)$ is homeomorphic over $\cl(D)$ to 
$\cl(D) \times \cl(\D)$. 
A section $\Psi:\cl(D) \to \cl(\bX)$ is a {\it trivial section} if there
is a fiber-preserving 
homeomorphism $h:\cl(\bX) \to \cl(D) \times \cl(\D)$ such
that $(h \circ \Psi)(\lambda) = (\lambda,0)$.
Given a Jordan bidisk $\bX$ which admits an extension to the boundary
we define the {\it frame} $\delta \bX$ as the torus
$\union_{\lambda \in \bd D} \union_{z \in \bd X_\lambda} (\lambda,z)$. 
A section $\Phi:D \to \bX$
is {\it proper} if it admits a continuous extension to
$\bd D$ and $\Phi(\bd D) \subset \delta \bX$.
Let $\Phi$ be a proper section and let $\Psi$ be a 
trivial section.
Let $\phi = \pi_2 \circ \Phi$ and $\psi = \pi_2 \circ \Psi$.
Define the {\it winding number} of $\Phi$ to be the winding
number of the curve $(\phi - \psi)|_{\bd D}$ 
around the origin.

\begin{lem}[Argument Principle]
\label{ArgPrinc}
Let $\bX$ be a Jordan bidisk over $D$ that admits an extension to the
boundary. Let $\Phi:D \to \bX$ be a proper section and let 
$\Psi:\cl(D) \to \cl(\bX)$ be a continuous section, holomorphic on $D$.
Let $\phi = \pi_2 \circ \Phi$ and $\psi = \pi_2 \circ \Psi$.
Suppose there are no solutions to $\phi = \psi$ on $\bd D$.
Then the number of solutions to $\phi = \psi$ counted with multiplicity is 
equal to the winding number of $\Phi$.
\end{lem}

Let $\union_j \bU_j$ be a pairwise disjoint collection of 
Jordan bidisks over $D$ with 
$0 \in U_\lambda = U_{0,\lambda}$. Let $\bV$ be a Jordan
bidisk over $D$ such that each $U_{j,\lambda}$
is compactly contained in $V_\lambda$.
Let 
$$
\bff:\union_j \bU_j \to \bV
$$ 
be a fiber-preserving holomorphic map such that each fiber map
$f_\lambda:\union_j U_{j,\lambda} \to V_\lambda$ is a generalized
quadratic-like map with critical point at the origin
and which on each branch $f_\lambda|_{U_{j,\lambda}}$
admits a holomorphic extension to a neighborhood of $U_{j,\lambda}$.
Let $\bfh$ be a holomorphic motion 
$$
h_\lambda : (\bd V_\ast,\union_j \bd U_{j,\ast}) \to
(\bd V_\lambda,\union_j \bd U_{j,\lambda})
$$
over $D$ with basepoint $\ast \in D$ which respects the dynamics.
We say $(\bff,\bfh)$ is a {\it holomorphic family of generalized 
quadratic-like maps over $D$}. When $\union \bU_j$ consists of only one
bidisk then the family is a {\it DH quadratic-like family}.
A family is {\it proper} if 
\begin{enumerate}
\item $\bV$ admits an extension to the
boundary
\item for each $z \in \union_j \bd U_{j,\ast}$ the
section $\lambda \mapsto (\lambda,h_\lambda(z))$ extends 
continuously to $\bd D$ and is a trivial section
\item the critical-value section
$\Phi(\lambda) = (\lambda,f_\lambda(0))$ is proper.
\end{enumerate}
The {\it winding number} of a proper family is the 
winding number of the critical value section.

\begin{thm}[\cite{DH2}]
\label{DHFamily}
If $(\bff,\bfh)$ is a proper
DH quadratic-like family over $D$ with winding number 1
then
$$
M(\bff,\bfh) = \{\lambda \in D :
	J(f_{\lambda}) \mbox { is connected}\}
$$
is homeomorphic to the standard Mandelbrot set $M$.
The homeomorphism is given by
the inner class map $\lambda \mapsto I(f_\lambda)$.
\end{thm}

We finish this section with the renormalization of a family.
Let $(\bff:\union_j \bU_j \to \bV,\bfh)$ 
be a proper holomorphic family of generalized
quadratic-like maps over $D$ with winding number 1.
If $0 \in R(f_\lambda, U_{0,\lambda})$
let $\bar{i}_\lambda$ be the {\it return itinerary} of 
$f_\lambda$: the 
(possibly empty) sequence of indices of off-critical pieces 
$\{U_{j,\lambda}\}$ through which
the critical point passes before returning to $U_{0,\lambda}$.
For such an $f_\lambda$ we can define a holomorphic motion $\bfh'$
of the boundaries of the domain and range of the return map to 
$U_{0,\lambda}$
by pulling back the holomorphic motion $\bfh$ by $f_\lambda$.
The motion $\bfh'$ has basepoint $\lambda$ and is defined over 
the neighborhood of $\lambda$ having the itinerary $\bar{i}_\lambda$.

\begin{lem}[{\cite[Lemma 3.6]{L4}}]
Let $(\bff:\union_j \bU_j \to \bV_j,\bfh)$ 
be a proper generalized quadratic-like family over $D$ with winding number 1.
Let $\ast \in D$ be the basepoint and let $g_\ast = R(f_\ast,U_{0,\ast})$.
Suppose $0 \in Dom(g_\ast)$.
Then the set 
$$
D' = \{\lambda \in D : \bar{i}_\lambda = \bar{i}_\ast\}
$$
is a Jordan disk and the family of first return maps
$(\bfg,\bfh')$ over $D'$ is proper and has winding number 1.
\end{lem}

\subsection{Parabolic periodic points}
\label{sec:ParabolicBackground}

The limits of maps with unbounded but
essentially bounded combinatorics are maps with
parabolic periodic points. This section reviews the local theory
near parabolic orbits and their perturbations. 
The main results are the existence and continuity of Fatou coordinates.
These results were proven in \cite{DH1} and \cite{La} for
perturbations lying in an analytic family and
later generalized in \cite{Sh}. Our presentation is based on \cite{Sh}.

Throughout this section we give
the space of holomorphic maps the ``compact-open topology with domains''.
A basis for this topology is given by the sets
$$
\cN(f,K,\eps) = \{g : |g(z)-f(z)| < \eps \mbox{ for }z \in K\}
$$
where $K \subset Dom(f)$ is compact and $\eps > 0$.
If a sequence of quadratic-like maps
converges to $f:U \to \C$ in the Carath\'eodory topology then it also
converges to $f:U \to \C$ in this topology.

Let $\Hol_0$ be the space of
holomorphic maps $f_0$ with a fixed point $\xi_0$
that is {\it parabolic} and {\it non-degenerate}: 
$f_0'(\xi_0) = 1$ and $f_0''(\xi_0) \neq 0$.
For example, choose any quadratic-like map
$f_0$ hybrid equivalent to $z^2+\qtr$.
Choose a neighborhood $N \ni \xi_0$ so that 
$f_0|_N$ is a diffeomorphism
and maps $N$ onto a neighborhood $N' \ni \xi_0$.

\begin{prop}[Fatou coordinates]
\label{FatouCoords}
Let $f_0 \in \Hol_0$ and choose $N$ and $N'$ as above.
Then there exist topological disks $D_\pm \Subset N \inter N'$, 
whose union forms
a punctured neighborhood of $\xi_0$ and which satisfy
$$
f_0^{\pm 1}(\cl(D_\pm)) \subset D_\pm \union \{\xi_0\}
\mbox{ and }
\inter_{n \ge 0} f_0^{\pm n}(\cl(D_\pm)) = \{\xi_0\}.
$$
Moreover, there exist univalent maps $\Phi_\pm : D_\pm \to \C$ 
such that
\begin{enumerate}
\item $\Phi_\pm$ are unique up to post-composing with a translation
\item $Range(\Phi_+)$ and $Range(\Phi_-)$ contain a right and left
half-plane, respectively
\item $\Phi_\pm(f_0(z)) = \Phi_\pm(z) + 1$
\end{enumerate}
%If we fix points $z_\pm \in D_\pm$ and normalize $\Phi_\pm$ by
%$\Phi_\pm(z_\pm) = 0$ then $\Phi_\pm$
%depend continuously on $f_0$ within some neighborhood of $f_0$
%in $\Hol_0$.
\end{prop}

%\realfig{ecalle1.fig}{ecalle01.eps}{Fatou coordinates.}{0.4\hsize}

The disks $D_\pm$ are called {\it incoming} and 
{\it outgoing petals} and the maps $\Phi_\pm$ 
are called the {\it Fatou coordinates}.
The Fatou coordinates induce conformal isomorphisms between the 
{\it \`Ecalle-Voronin cylinders} $\cC_\pm = D_\pm/{f_0}$
and $\C / \Z$. 
Let $\pi_\pm$ denote the projection of $D_\pm$ to $\cC_\pm$ and
extend $\pi_+$ to the attracting basin of $\xi_0$ by 
$\pi_+(z) = \pi_+(f_0^n(z))$ for a large enough $n$.

A {\it transit map} $g:\cC_+ \to \cC_-$ 
is a conformal isomorphism which respects the ends $\pm \infty$. 
A holomorphic map $\wtl{g}:U \to \C$ is a {\it local lift} of 
a transit map $g$ if
$\cl(U) \subset D_+$, $Range(\wtl{g}) \subset D_-$, and
$$
g \circ \pi_+ = \pi_- \circ \wtl{g}.
$$
When written in Fatou coordinates,
$\wtl{g}$ is a translation $T_a$ by a complex number $a$.
The quantity
$\bar{a} = a \mod \Z $, called the
{\it phase}, depends only on $g$ (and the normalization of Fatou
coordinates) and uniquely specifies $g$. 
We will use the notation $g_{\bar{a}}$ to denote the transit map with 
phase $\bar{a}$.

To simplify future notation, let 
$\Phi = \Phi_+$ and $\phi = \Phi_-^{-1}$.
Also, we shall freely use the notation
$\Phi_n$, $\Phi_f$, ${\cal C}_{n,\pm}$, etc to indicate a 
dependence on an index $n$ or map $f$.

We now consider perturbations of $f_0 \in \Hol_0$. 
Since $\xi_0$ is a non-degenerate parabolic fixed point the 
generic perturbation will cause it
to bifurcate into two nearby fixed points $\xi_f$ and $\xi_f'$
with multipliers $\lambda_f$ and $\lambda'_f$, respectively.
Let $N$ be the neighborhood of $\xi_0$ chosen
for \propref{FatouCoords} and let $\Hol$ be 
the space of holomorphic maps which are
diffeomorphisms of $N$.
Let $\Hol_1$ be the set of $f \in \Hol$ with exactly
two fixed points $\xi_f$ and $\xi'_f$ in $N$ satisfying
\begin{equation}
\label{eq:goodmultiplier}
\arg(1-\lambda_f), \arg(1-\lambda'_f) \in 
[\pi/4,3\pi/4] \union [-3\pi/4,-\pi/4].
\end{equation}

\begin{thm}[Douady coordinates]
\label{pertFatou}
Let $f_0 \in \Hol_0$.
There is a neighborhood $\cN$ of $f_0$ such
that if $f \in (\cN \inter \Hol_1)$
then there exist univalent maps $\Phi_f = \Phi_{f,+}$
and $\phi_f = (\Phi_{f,-})^{-1}$, 
unique up to translation,
and a constant $a_f \in \C$ satisfying
\begin{enumerate}
	\item $\Phi_f(f(z)) = \Phi_f(z)+1$ and
		$\phi_f(w+1) = f(\phi_f(w))$ where defined
	\item ${\cal C}_{f,+} = Dom(\Phi_f) / f$ and 
		${\cal C}_{f,-} = Range(\phi_f) / f$ are
		conformally cylinders and one can choose
		fundamental domains $S_{f,\pm}$
		to depend on $f \in \Hol_1$ 
		continuously in the Hausdorff topology.
	\item (see \figref{ecalle2.fig}) 
		for $z \in S_{f,+}$ there is an $n > 0$ such that
		$f^n(z) \in S_{f,-}$ and for $n$ minimal
\begin{equation}
\label{eq:transit}
		f^n(z) = (\phi_f \circ T_{a_f+n} \circ \Phi_f)(z).
\end{equation}
\end{enumerate}
If we fix points $z_\pm \in D_\pm$ and normalize $\Phi_{f,\pm}$
by $\Phi_{f,\pm}(z_\pm) = 0$ then $\Phi_{f,\pm}$
depend continuously on $f \in \cN \inter (\Hol_0 \union \Hol_1)$.
\end{thm}

\realfig{ecalle2.fig}{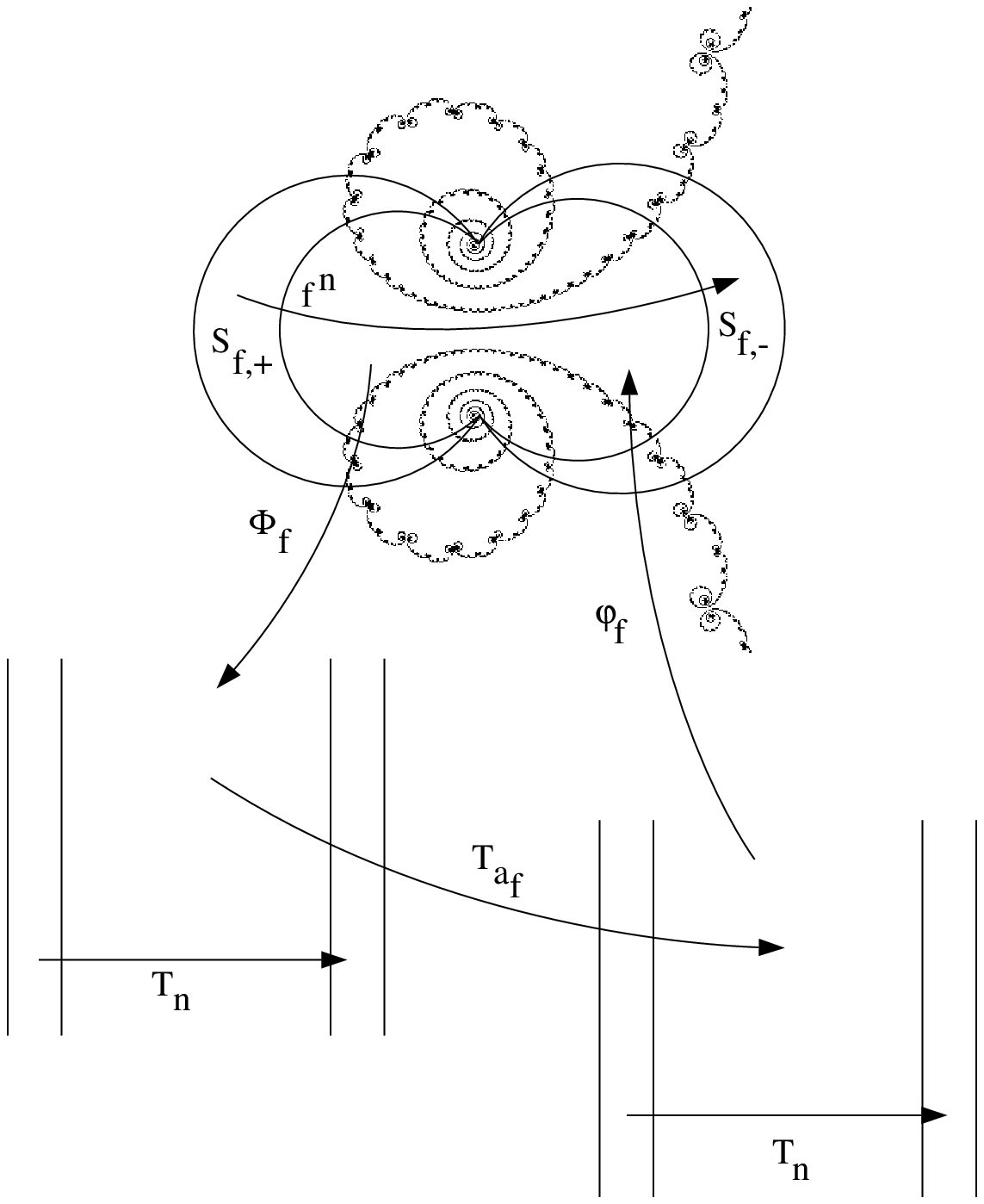}{Perturbed Fatou coordinates.}{0.8\hsize}

Suppose $f_0 \in \Hol_0$ and
$f \in \Hol_1 \inter \cN$ where $\cN$ is from 
\thmref{pertFatou}.
The discontinuous map from $S_{f,+}$ to $S_{f,-}$ defined by
equation \ref{eq:transit} projects to 
a transit map $g_f:\cC_{f,+} \to \cC_{f,-}$
with phase $\bar{a}_f = a_f \mod \Z$. This map describes how a long orbit
of $f$ ``passes though the gate'' between $\xi_f$ and $\xi'_f$.
The following lemma relates the convergence
of $\bar{a}_f$ to the convergence of local lifts.

\begin{lem}
\label{converge1}
Let $f_k \in \Hol_1$ converge to $f_0 \in \Hol_0$
and suppose $\bar{a}_{f_k} \to \bar{a}$. Then
for any local lift $\wtl{g}$ of $g_{\bar{a}}$ there exists a sequence $n_k$
such that
$$
f^{n_k}_k \to \wtl{g}
$$
uniformly as $k \to \infty$.
\end{lem}
\begin{pf}
Let $K = \cl(Dom(\wtl{g}))$ and define $a \in \C$ by
$\wtl{g} = \phi \circ T_a \circ \Phi$.
Let $K_1$ be a compact set in $\C$ containing $\Phi(K)$ in its interior
and let $K_2$ be a compact set in $Dom(\phi)$
containing $T_a(K_1)$ in its interior.
Let $a_k$ be the constant $a_{f_k}$ in \propref{pertFatou}. 
Since $\bar{a}_{f_k} \to \bar{a}$
there exists a sequence $n_k$ so that $a_k+n_k \to a$.

For $k$ large enough $K \subset Dom(\Phi_{f_k})$,
$\Phi_{f_k}(K) \subset K_1$,
$T_{a_k+n_k}(K_1) \subset K_2$ and
$K_2 \subset Dom(\phi_{f_k})$. 
The lemma follows from equation (\ref{eq:transit}) and since
$\Phi_{f_k},T_{a_k+n_k},\phi_{f_k}$
converge to $\Phi,T_a,\phi$ uniformly on $K,K_1,K_2$, respectively.
\end{pf}

The following lemma gives a simple condition under which perturbed 
Fatou coordinates exist.
\begin{lem}
\label{Tangent}
Suppose $f_n$ is a sequence of quadratic-like maps converging 
in the Carath\'eodory topology to a quadratic-like map
$f \in \Hol_0$. Suppose the fixed points of $f_n$ are repelling. Then
$f_n \in \Hol_1$ for $n$ large enough.
\end{lem}
\begin{pf}
Using the holomorphic index (see \cite{M1}) one can prove that 
$$
\frac{1}{1-\lambda_{f_n}} + \frac{1}{1-\lambda'_{f_n}}
$$
converges as $f_n \to f$. Since $\lambda$, $\lambda' \in \C \setm \D$
it follows
$$
|\arg(1-\lambda_{f_n})| \to \pi/2 \mbox{ and }
|\arg(1-\lambda'_{f_n})| \to \pi/2
$$
as $n \to \infty$. In particular, $f_n \in \Hol_1$ for $n$ large.
\end{pf}

For any $z \in D_+ \inter D_-$
define the \`Ecalle-Voronin transformation $\cE$ by
$$
\cE(\pi_-(z)) = \pi_+(z).
$$
One can show that $\cE$ extends holomorphically to the two ends of $\cC_-$
by using the Fatou coordinates and the standard isomorphism
$\pi(z)=\exp(2\pi iz)$ of $\C / \Z$ to ${\C} \setm 0$.
The following lemma is useful for controlling
the dynamics near the ends of the \`Ecalle-Voronin cylinders.

\begin{lem}
\label{EcalleExpands}
Suppose $f_0 \in \cH(\qtr)$ 
and $g:\cC_+ \to \cC_-$ is a transit map such that
the critical point of $f_0$ escapes $K(f_0)$ under iterates
of $f_0$ and local lifts of $g$. Then
$$
|(g \circ \cE)'(\pm\infty)| > 1.
$$
\end{lem}
\begin{pf}
We will prove the lemma with the critical point escaping after just
one iterate of a local lift of $g$.
Assume the critical point of $f$ is at the origin.
Let $R = g \circ \cE$ and $J_- = \pi_-(J(f_0))$.
Let $V_{\pm \infty}$ denote the connected components of
$(\cC_- \setm J_-) \union \{\pm\infty\}$
containing $\pm\infty$ and
let $U_{\pm \infty} = g^{-1}(V_{\pm \infty})$ (see
\figref{ecalleExpand.fig}).

\realfig{ecalleExpand.fig}{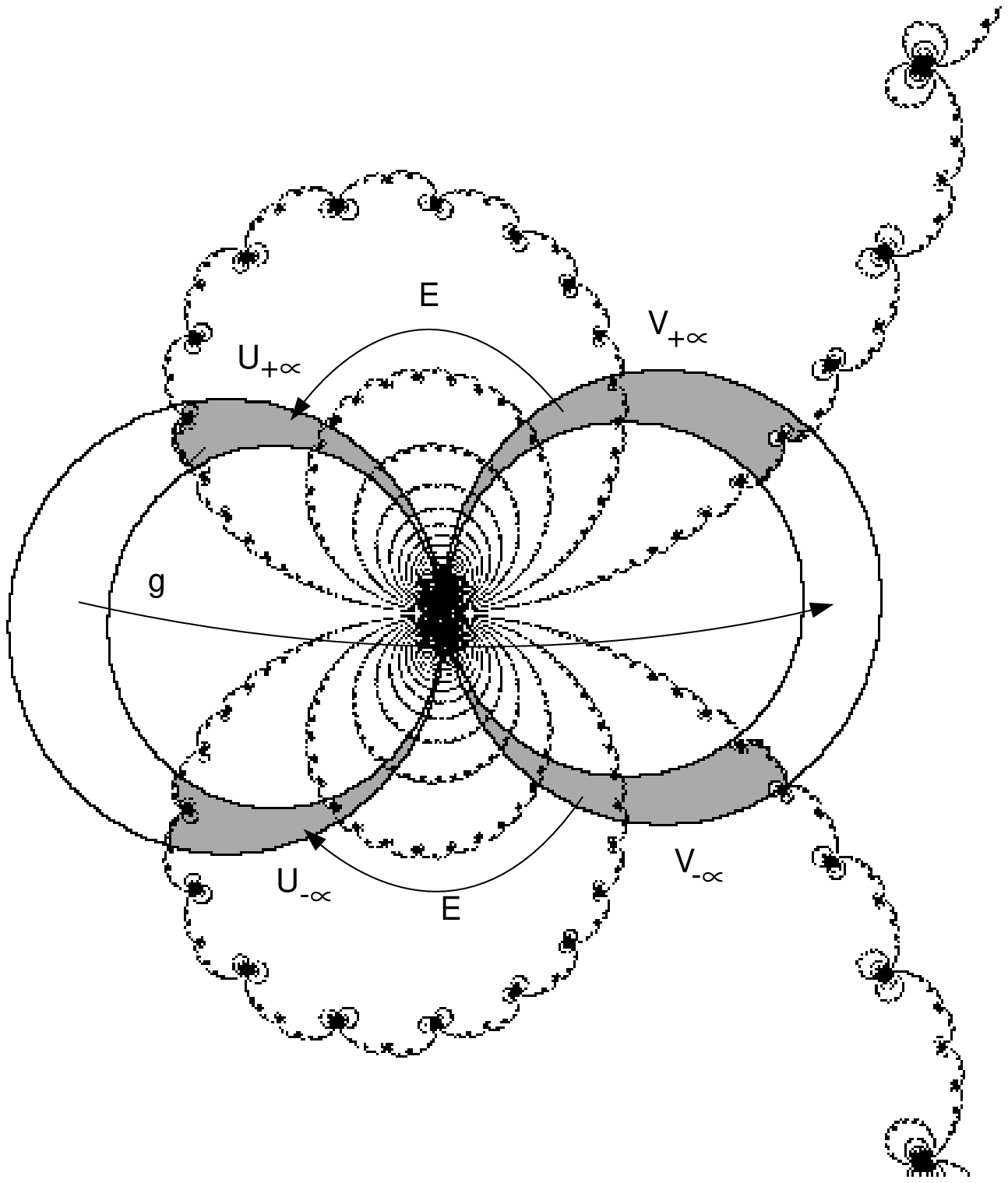}{A blow-up of the Julia set of 
$f_0 = z^2+\qtr$ with pre-images by $f_0$, $g$ and $\cE$ 
highlighting the sets $U_{\pm\infty}$ and $V_{\pm\infty}$.}{0.6\hsize}

Note that $\cE$ can be extended to $V_{\pm\infty}$ as a branched cover.
The set of critical points is the backward orbit of $0$
and the only critical value is $\pi_+(0)$.

Since $\pi_+(0) \not\in U_{\pm \infty}$ and each $U_{\pm\infty}$
is simply connected there is a branch of
$\cE^{-1}$ defined on $U_{\pm \infty}$ preserving $\pm\infty$.
Composing $\cE^{-1} \circ g^{-1}$ we have constructed 
a branch of 
$R^{-1}$ which maps each $V_{\pm \infty}$ strictly inside itself 
and fixes $\pm \infty$. The lemma follows from the Schwarz lemma.
\end{pf}

We close this section with a lemma on the geometry of 
some particular first through maps. Let $m > 0$ and let 
$$
X = \{f \in Quad(m) : f \in \cH(\qtr) \mbox{ and } \diam K(f) = 1 \}.
$$
From \lemref{Compact}, $X$ is compact (in the Carath\'eodory topology).
For each $f \in X$ choose a neighborhood $N \ni \beta(f)$ on which
$f$ is a diffeomorphism and let $\cN_1, \dots, \cN_k$ 
be a finite cover of $X$ by the neighborhoods from \thmref{pertFatou}. 
In order to preserve certain compactness properties, we will
need $\cN_i$ to be closed neighborhoods.
%Define the space $Quad' \supset Quad$ as the set of quadratic-like maps
%with (possibly escaping) critical point at the origin and let
%$$
%X_i' = \{f \in Quad'(m) : f \in \cN_i \inter (\Hol_0 \union \Hol_1)
%\mbox{ and } \diam K(f) = 1 \}.
%$$
By rescaling we can extend the 
neighborhoods $\cN_i$ to be a finite cover of
$\{f \in Quad(m) : f \in \cH(\qtr) \}$.
Note the coordinates do not necessarily agree on the overlaps 
$\cN_i \inter \cN_j$.

Now suppose $f \in Gen(m)$ is a generalized
quadratic-like map such that the critical point escapes the central 
component $U_0$. Suppose $f \in (\cN_i \inter \Hol_1)$ 
for some $1 \le i \le k$ and suppose the critical point of 
$f$ passes once through the gate before landing
in the off-critical pieces.
In this case we say $f$ has a {\it saddle-node cascade}. 
Let $T = T(f,\union_{j \neq 0} U_j)$ be the first through map of $f$ and
define the {\it modified landing time} $l(z)$ of $z \in Dom(T)$ as follows. 
For each $\cN_i \ni f$, there is a 
choice of fundamental domains $S_{f,\pm}$.
Write $T(z)$ as a composition of $f$, $(f|_{N_i})^{-1}$ 
and of the discontinuous
map $\wtl{g}_f:S_{f,+} \to S_{f,-}$ defined in equation (\ref{eq:transit}).
Define $l_i(z)$ as the minimal number of maps in this composition
and define $l(z) = \max_i l_i(z)$.

\begin{lem}
\label{rootlimit}
Let $f_k \in Gen(m)$ be a sequence of maps with saddle-node
cascades such that the modified landing times $l(0)$ are bounded.
Suppose $f_k \to f$ and $K(f|_{U_0})$ is connected.
Then $f|_{U_0} \in \cH(\qtr)$.
\end{lem}
\begin{pf}
Fix some neighborhood $\cN_i$ containing $f_k$ and $f$.
Let $n_k$ be the transit time defined by equation \ref{eq:transit}
for the orbit of the origin.
If $f|_{U_0} \not\in \cH(\qtr)$ then $n_k$ is bounded from above.
But then the origin escapes $U_{k,0}$ under a bounded number of
iterates of $f_k$, which is a contradiction.
\end{pf}

\begin{lem}
\label{throughgeo}
Let $\lambda > 0$, $m > 0$ and $r \in \N$.
Suppose $f \in Gen(m)$ has a saddle-node cascade and $\geo(f) \ge \lambda$.
Then the phase $\bar{a}_f$ of the induced
transit map lies in a pre-compact subset of $\C / \Z$.
Suppose $g \in Gen$ is a restriction of the first through map $T$ such that
$Dom(g) = Comp(Dom(T),u_g)$ and for $z \in u_g$,
$$
l(z) \le r.
$$
Then there exists $C(\lambda,m,r) > 0$ such that $\geo(g) \ge C$.
\end{lem}
\begin{pf}
Let us prove the first statement. 
Suppose $f \in \cN_i$.
Let $c_1 = f^{r_1}(0)$ be the first moment when the orbit of $0$
lands in $S_{f,+}$. We can assume $r_1$ is uniform over the
neighborhood $\cN_i$. Then $c_1$ lies in a pre-compact subset
of $\cC_{f,+}$. Let $c_2 = \wtl{g}_f(c_1)$.
Since $f^n(c_2) \in \union_{j \neq 0} U_j$ for some
$n \le r$, it follows $c_2$ lies in a pre-compact subset of
$\cC_{f,-}$. Hence the phase $\bar{a}_f$, measured in the coordinates
from $\cN_i$, lies in a pre-compact subset of $\C / \Z$. 

The bound on the geometry is clear since the perturbed Fatou coordinates
converge and the transit maps $g_f$ lie in a pre-compact subset and the
number of iterates of $g_f$, $f$ and $(f|_{N_i})^{-1}$ is bounded.
\end{pf}

\section{Combinatorics}
\subsection{Essentially period tripling}
\label{sec:special}

An $\infty$-renormalizable map $f$ 
has {\it bounded
combinatorics} if $\tau(f)$ contains a finite number of distinct maximal tuned
Mandelbrot sets, or,
equivalently, if $\bar{\sigma}(f)$ contains a finite number of distinct shuffles.
In this section we construct an infinite set of maximal tuned Mandelbrot
sets with bounded essential period. 
Hence any map whose tuning invariant is chosen
from these Mandelbrot sets will have essentially bounded combinatorics. 
On the other hand if the tuning invariant contains an infinite number of
distinct Mandelbrot sets then the map will not have bounded combinatorics.
The simplest way to construct such a collection of maximal tuned Mandelbrot
sets is by perturbing in a particular way
the map $z^2-1.75$, the root point of the period three
tuned copy. 
For this reason we say these copies are essentially period tripling.

Let $f(x) = x^2-1.75$ and let $\xi$ be the parabolic periodic orbit of
period three.
Recall $A = A(f) = [\alpha,\alpha']$.
Let $g$ be the first return map of $f$ on $A$ (see \figref{picshuffle2}). 
Let $I^1$ and $I^1_1$ be the two indicated intervals
satisfying $g|_{I^1} = f^3$ and $g|_{I^1_1} = f^2$.

\realfig{picshuffle2}{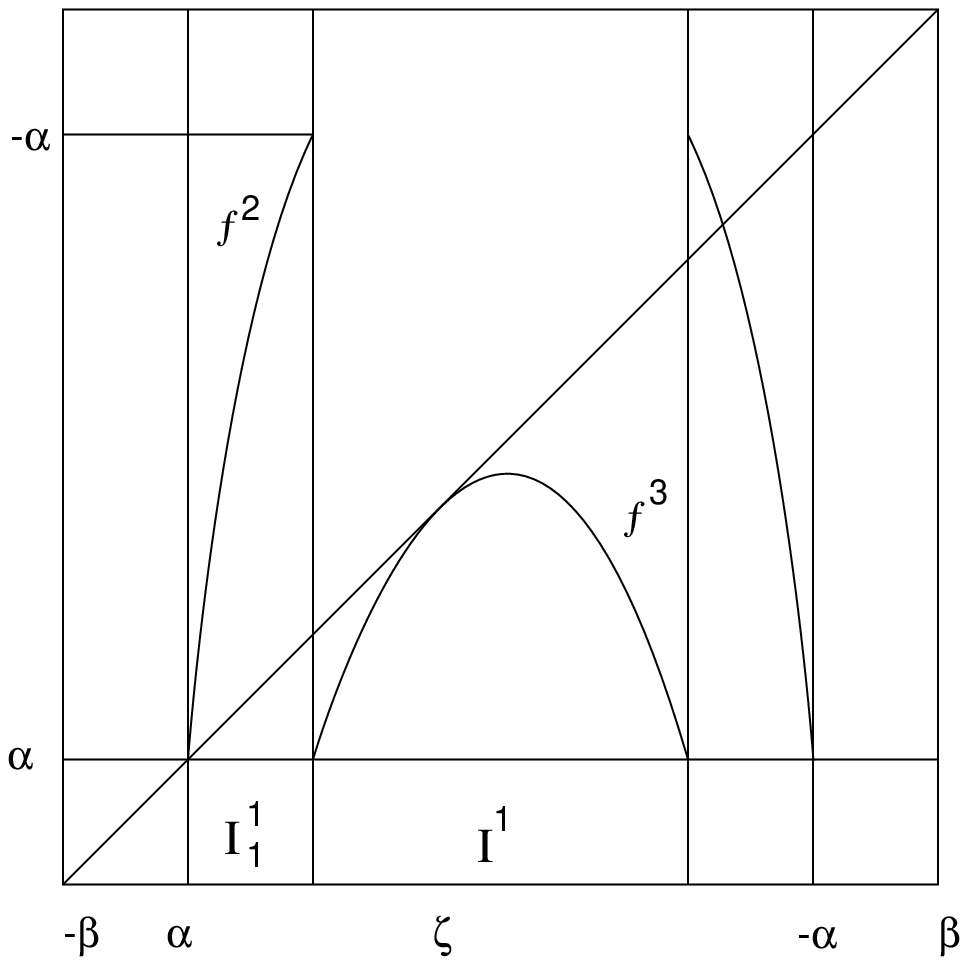}{The first return map for $x^2-1.75$.}
{0.4\hsize}

Fix a small $\eps > 0$ and consider $c \in (-1.75,-1.75+\eps)$. The 
periodic point $\xi$ bifurcates and the orbit of 
the critical point under $f^3_c$ now
escapes the interval $I^1$.
Let $c_n$ be the parameter value (see \figref{picshuffle}) so that for $f =
f_{c_n}$,
\begin{itemize}
\item $f^{3i}(0) \in I^1$ for $i = 1,\dots,n-1$,
\item $f^{3n}(0) \in I^1_1$,
\item $f^{3n+2}(0) = 0$.
\end{itemize}
\realfig{picshuffle}{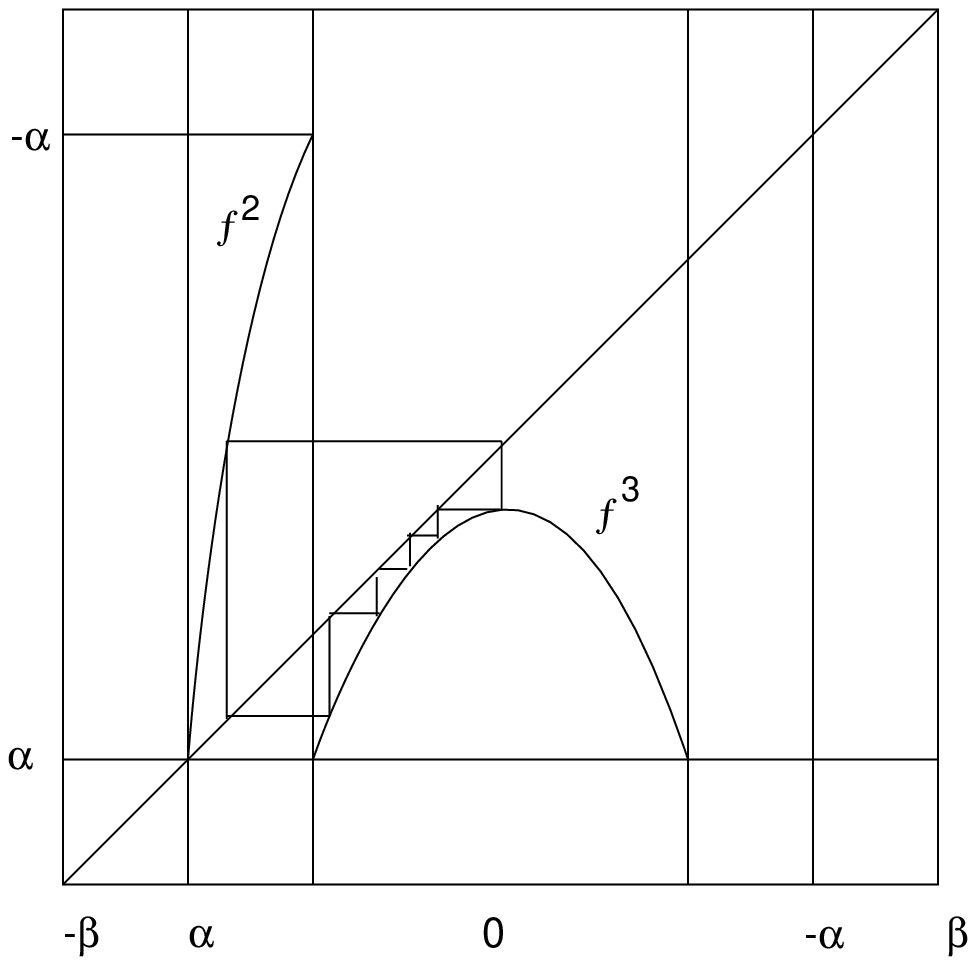}{The first return map for $x^2+c_5$
and the orbit of the origin.}
{0.4\hsize}

In the next section we will justify the claim
that $c_n$ exists and
is the center of a maximal tuned Mandelbrot set, denoted
$M^{(3)}_n$. Equivalently, if we
let $\sigma^{(3)}_n$ be the permutation induced on the orbit 
by $f_{c_n}$ of the origin
labeled from left to right, then $\sigma^{(3)}_n$ is a shuffle. 
In \figref{fig:per3M} we have drawn
the period three tuned Mandelbrot set and a few of the 
$M^{(3)}_n$ accumulating
at its root point. Any map $f_c$ 
for $c \in M^{(3)}_n$ will be renormalizable with 
essential period $p_e(f_c) = 5$.

\realfig{fig:per3M}{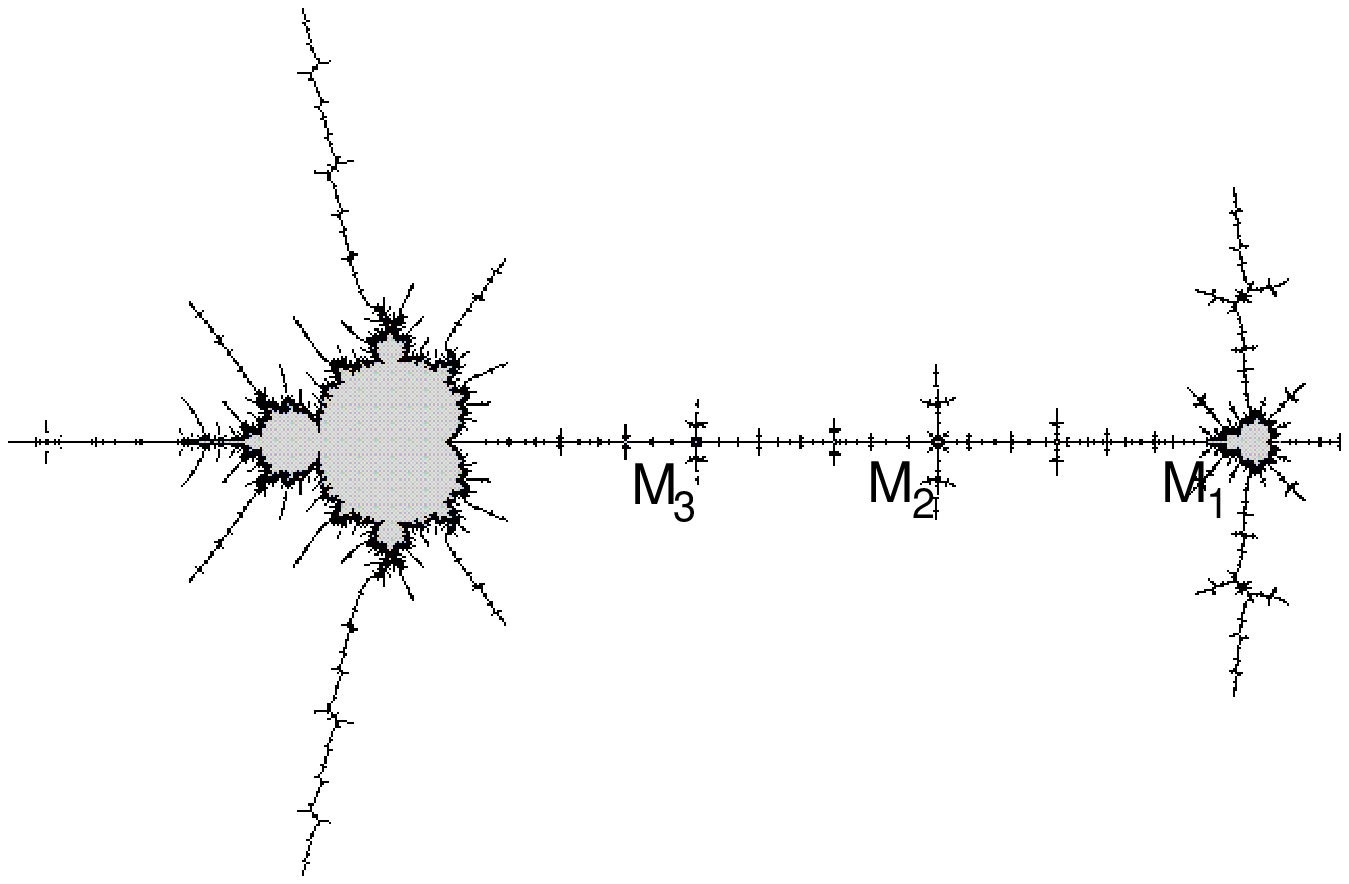}{The Mandelbrot set near the real period three tuned copy.}
{0.4\hsize}

\begin{figure}[hbtp] 
         \centerline{\psfig{figure=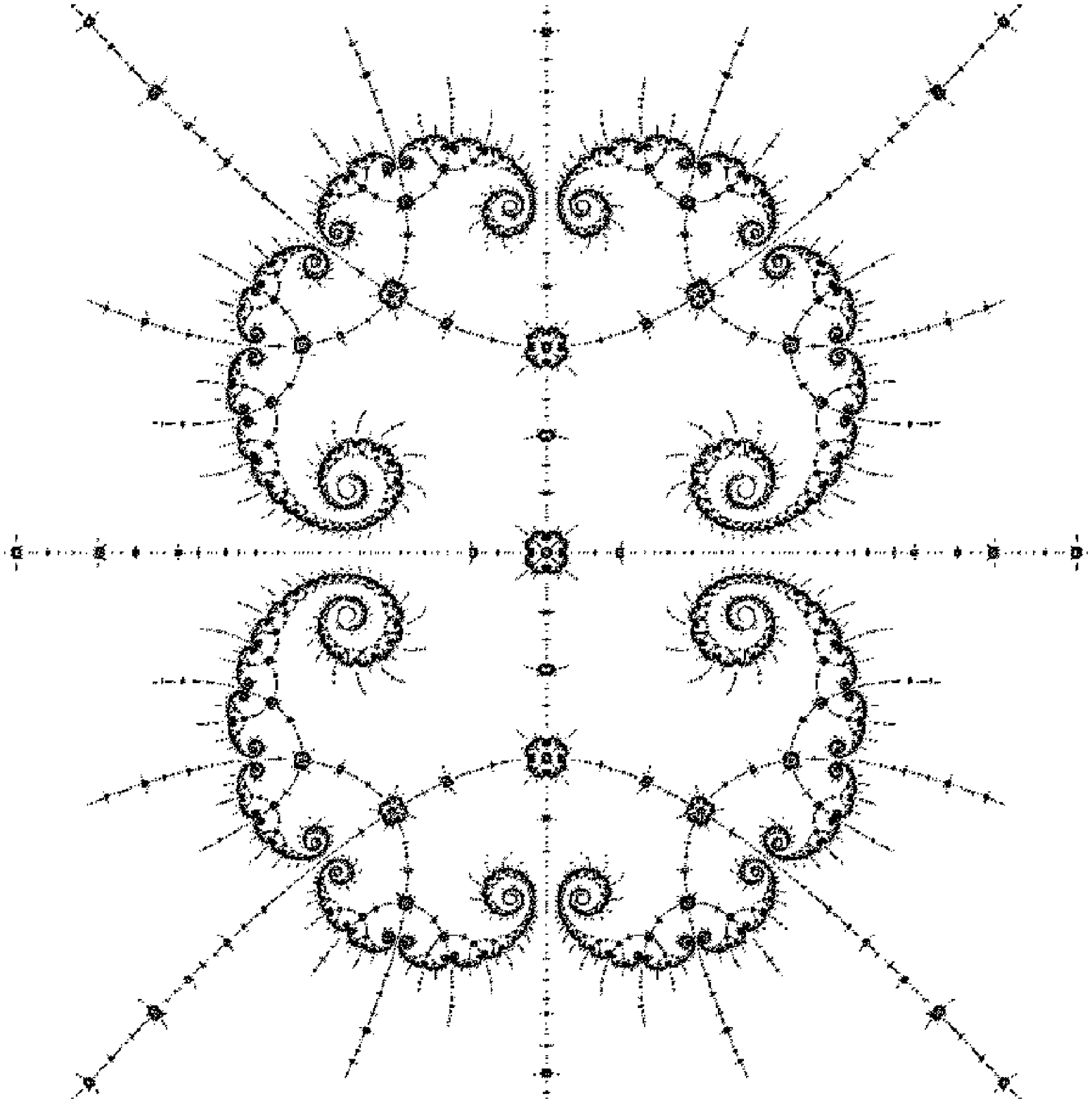,width=6cm}
\hskip 1cm  
         \psfig{figure=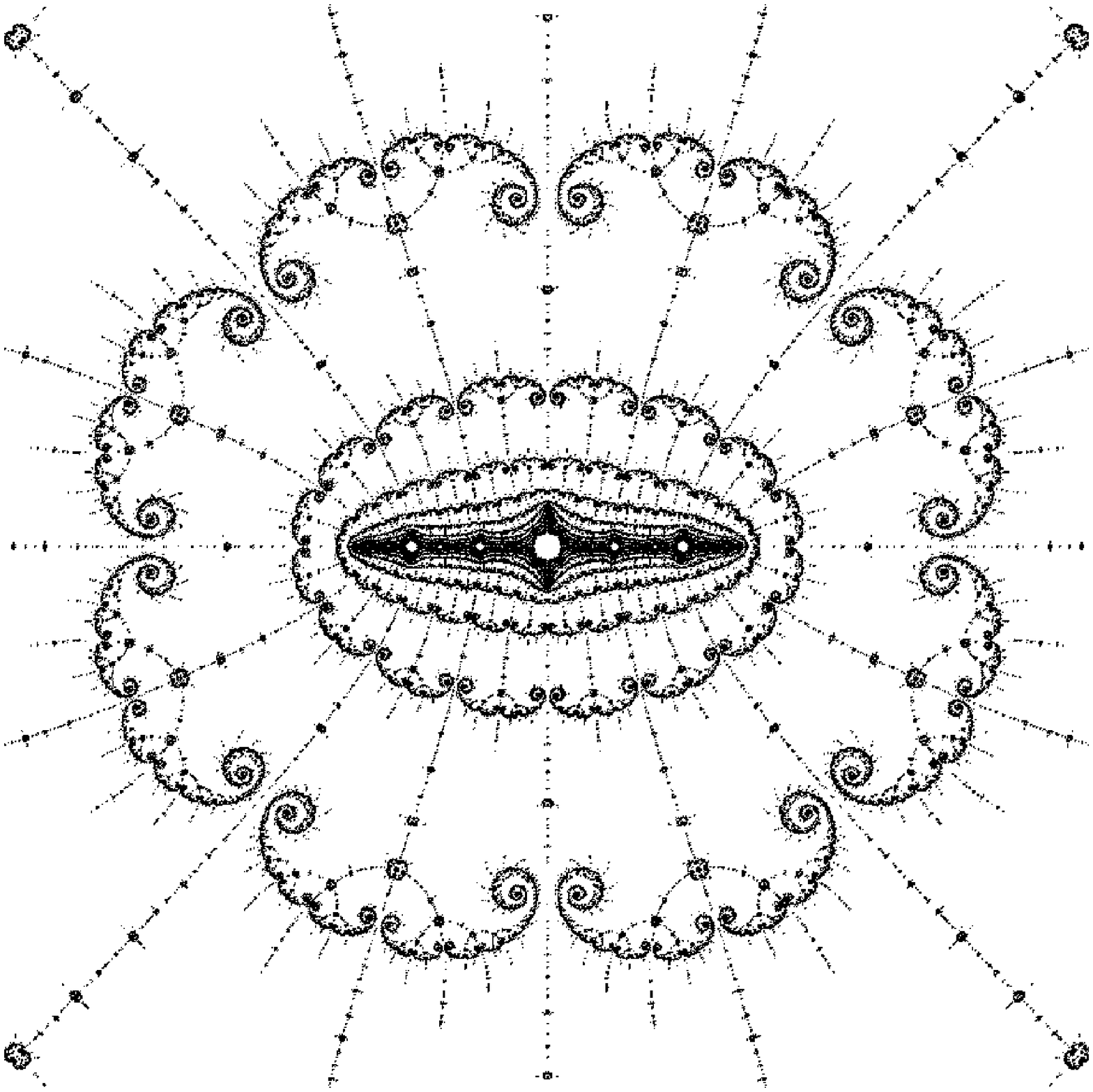,width=6cm}}
         \caption{Blow-ups of $J(f)$ near the origin. }
\label{blow-up}
	\end{figure}

In \figref{two-jul} we have drawn the filled Julia sets for $z^2-1.75$
and for $z^2-c_n$ for some $c_n$ with $n$ large.
\figref{blow-up} shows two blow-ups of the Julia set of 
$f = z^2-c_n$. The ``ghost'' boundary of the basin of
$\xi$ is visible in the left picture and the pre-images of this ghost
boundary nest down to $J(\cR(f))$ in the right picture.

\realfig{two-jul}{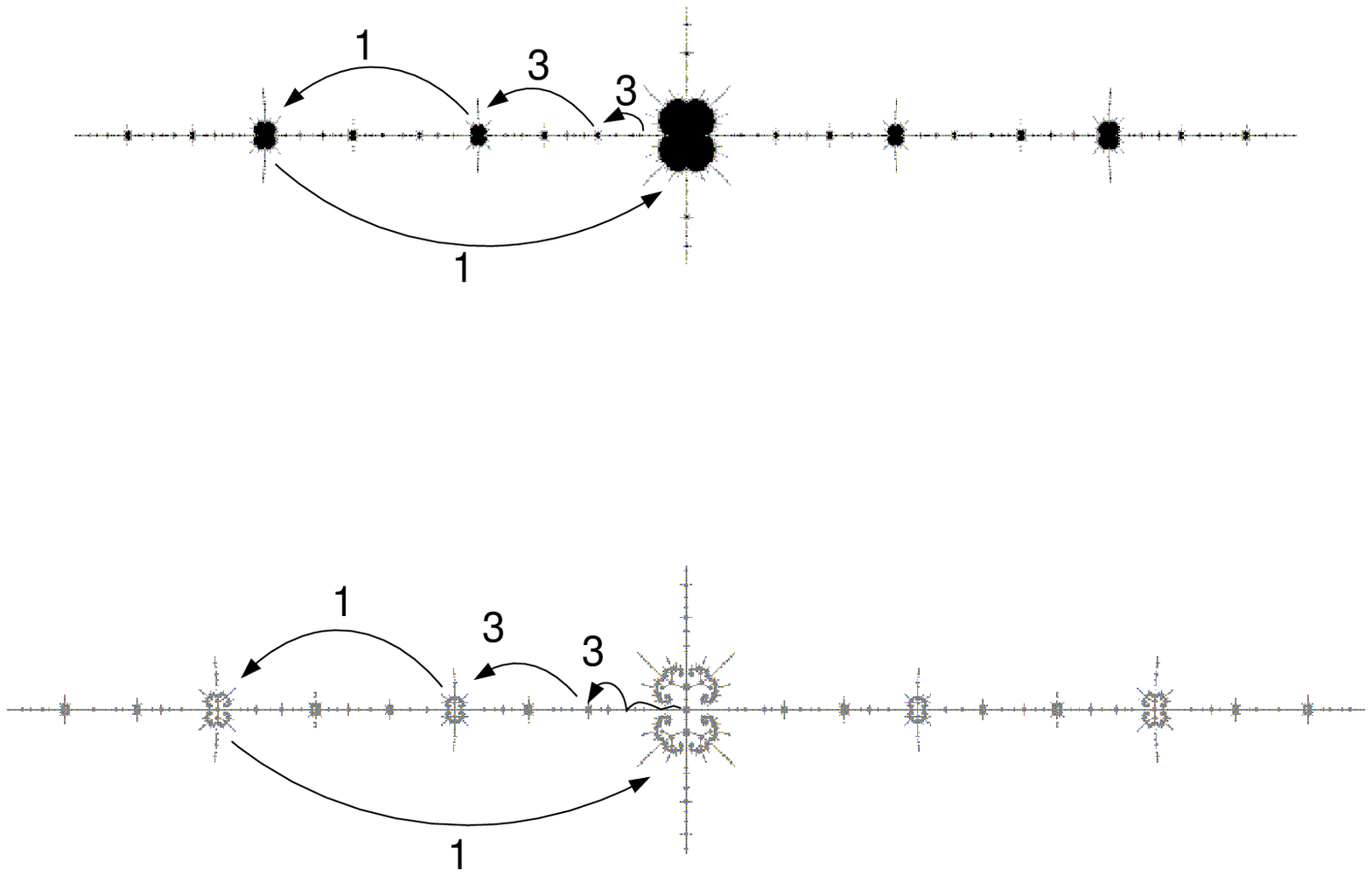}{The filled Julia set of the map $z^2-1.75$ (above)
and of $f = z^2-c_n$ for some large $n$.}{0.7\hsize}

\subsection{Essentially bounded combinatorics}
\label{essentially-bounded}

In this section we describe the {\it return type sequence}
of a given shuffle $\sigma$ 
and define the {\it essential period} $p_e(\sigma)$.

Suppose $f \in Quad$ is renormalizable, has 
real combinatorics, and $\sigma(f) \neq \sigma^{(2)}$.
In real contexts we will assume $f \in RQuad$.
Define the complex principal nest 
$V^0 \supset V^1 \supset V^2 \supset \dots$ 
of $f$ as follows. Choose a straightening of
$f$ to a polynomial $f_c$ and pull the equipotential and external ray
foliations of $f_c$ back to $f$. Cut the domain $D$ 
bounded by a fixed equipotential level by the closure of the rays 
that land at $\alpha$ and at $\alpha'$. 
The resulting set of connected components is called the 
{\it initial Yoccoz puzzle}.
Let $V^0$ be the component
containing $0$ and let 
$$
V^m = Comp(Dom(R(f,V^{m-1})),0).
$$
For $m \ge 0$ let $I^m = V^m \inter \R$.

For $m \ge 1$, let $g_m: \union_{i} V^{m}_i \to V^{m-1}$
be the {\it generalized renormalization} of $f$ on $V^{m-1}$.
We will also denote the restriction to the real line
$g_m:\union_i I^m_i \to I^{m-1}$ by $g_m$. 
Number the intervals $\union_i I^m_i$ (and domains $V^m_i$) from left to
right and so that $0 \in I^m_0 = I^{m+1}$. 
See \figref{picshuffle} for
an example of the first two levels of the real principal nest and the graph
of $g_1$.

\begin{lem}[\cite{L3}]
\label{initialgeo}
Let $m > 0$, $n > 0$ and let $f \in Quad(m)$. Suppose $f$ has real combinatorics,
is not immediately renormalizable, and the return time of any $z \in Dom(g_1)$
through the initial Yoccoz puzzle until the first return to $V^0$
is bounded above by $n$. Then $\mymod(Dom(g_1),Range(g_1)) \ge m_0(m,n) > 0$,
$\geo(g_1) \ge C(m,n) > 0$
and $\diam K(g_1) / \diam K(f) \ge C'(m,n) > 0$.
\end{lem}

The {\it return type} of $g_m$ is defined as follows
(see \cite{L6} for details).
Let $g \in RGen$ have finite type and
let $\union_i I_i = Dom(g) \inter \R$ numbered from left to right with
$0 \in I_0$. Let $(\Gamma,\eps)$ be the free ordered signed semigroup generated
by $\{I_i\}$ where $\eps:\{I_i\} \to \{\pm1\}$ 
is the sign function defined for $i \neq 0$ by
$\eps(I_i) = +1$ iff
$g|_{I_i}$ is orientation preserving and
for $i=0$ by $\eps(I_0) = +1$ iff $0$ is a local minimum of $g$.
Let $h \in RGen$ be
a restriction of $R(g,I_0)$ to finitely many components of its domain
and let $\union_j J_j = Dom(h) \inter \R$.
Let $(\Gamma',\eps')$ be the corresponding signed semigroup for $h$.
Let $\chi:(\Gamma',\eps') \to (\Gamma,\eps)$ 
be the homomorphism generated by assigning to each $J_j$ the word
$I_{i_1} I_{i_2}\cdots I_{i_n}$
where $I_{i_k}$ is the interval containing $g^k(J_j)$ and $n$ is
the return time of $J_j$ to $I_0$.
The homomorphism $\chi:\Gamma' \to \Gamma$ is the
{\it return type} of $h$.

A homomorphism $\chi:(\Gamma',\eps') \to (\Gamma,\eps)$
between free ordered signed semigroups
is called {\it unimodal} if the image of every generator
is a word ending with the central interval and if the map is strictly
monotone on the intervals to the right and left of center and has
an extremum at the center. 
We say a unimodal $\chi$ is {\it admissible} if 
$$
\eps'(I_j') = \sgn(j) \eps(\chi(I'_j)) \mbox{ for } j \neq 0 \mbox{ and }
\eps'(I_j') = \eps(\chi(I'_j)) \mbox{ for } j = 0.
$$

Let us describe the initial combinatorics of $f$.
Let $(\Gamma_0,\eps_0)$ be the signed semigroup generated by
$+I^0$ and two intervals $-I^0_{-1}$
and $+I^0_1$.
We say a unimodal homomorphism $\chi:(\Gamma,\eps) \to (\Gamma_0,\eps_0)$ 
is {\it zero-admissible} if it is admissible and additionally
for each $I_i$ there is a $p_i \ge 0$ with $p_0 \ge 1$ and
such that 
$$
\chi(I_i) = I_{-1} I_1^{p_i} I_0.
$$
The initial combinatorics of $f$ is described by the homomorphism
assigning to each $I^1_i$ its itinerary by $f$ through the
intervals $I^0$ and the connected components of $B(f) \setm I^0$.
In general if $h_1$ is any restriction of the first return map to $V^0$
then the return type of $h_1$ is the homomorphism mapping to any interval
in its domain its itinerary through the above intervals.
Note that if $f$ has negative orientation then repeat the
construction with all signs reversed.

The combinatorics of $f$ up to level $m$ is described
by the sequence $S_m$ of admissible unimodal homomorphisms
$$
\Gamma_m \upto{\chi_m} \Gamma_{m-1} 
 \upto{\chi_{m-1}} \cdots \upto{\chi_2} \Gamma_1 \upto{\chi_1} \Gamma_0
$$
where $\chi_m$ is the return type of $g_m$ and
$\chi_1$ is zero-admissible.
Each $S_m$ is {\it irreducible}, meaning the orbit of the critical point
enters every interval $I^m_i$.
Since $f$ is renormalizable there exists an
$m'$ such that $\Gamma_m$ is the semigroup with one generator
for all $m \ge m'$.
Let $S(\sigma) = S_{m'}$ for the smallest such value of $m'$.
Then the shuffle $\sigma(f)$ is uniquely specified by $S_{m'}$.
Moreover, we have the following

\begin{thm}[\cite{L6}]
\label{shuffleexist}
Let $S$ be an irreducible
finite sequence of admissible unimodal homomorphisms:
$$
\Gamma_m \upto{\chi_m} \Gamma_{m-1} 
 \upto{\chi_{m-1}} \cdots \upto{\chi_2} \Gamma_1 \upto{\chi_1} \Gamma_0.
$$
Suppose $\Gamma_m$ is the only semigroup with one generator,
$\Gamma_0$ is as above and $\chi_1$ is zero-admissible.
Then there is a  unique shuffle $\sigma$ such that
$S(\sigma) = S$.
\end{thm}

We can now justify our construction of the essentially period
tripling shuffles $\sigma^{(3)}_n$ from \secref{sec:special}.
Consider the following signed semigroups generated by the specified
intervals
\begin{equation}
\label{eq:groups}
\begin{array}{l}
\Gamma = \seq{+I_{-1}, -I_0} \\
\Gamma' = \seq{-I_0}
\end{array}
\end{equation}
and consider the following homomorphisms 
\begin{equation}
\label{eq:homo}
\begin{array}{lll}
\chi_0:\Gamma \to \Gamma_0 & \mbox{ generated by }
& I_{-1} \mapsto I_{-1} I_0 \mbox{ and } I_0 \mapsto I_{-1} I_1 I_0 \\
\chi:\Gamma \to \Gamma & \mbox{ generated by }
& I_{-1} \mapsto I_{-1} I_0 \mbox{ and } I_0 \mapsto I_0 \\
\chi':\Gamma' \to \Gamma & \mbox{ generated by } 
& I_0 \mapsto I_{-1} I_0.
\end{array}
\end{equation}
Then the sequence corresponding to
the essentially period tripling combintorics $\sigma^{(3)}_n$ is
$$
\Gamma' \upto{\chi'} \Gamma \upto{\chi} \Gamma  \upto{\chi} \cdots
 \upto{\chi} \Gamma \upto{\chi_0} \Gamma_0
$$
where $\chi$ is repeated $n-1$ times.

A level $m > 0$ is called {\it non-central} iff 
$$
g_{m}(0) \in V^{m-1} \setm V^{m}.
$$
Let $m(0) = 0$ and let $0 < m(1) < m(2) < \cdots < m(\kappa)$ 
enumerate the non-central levels, if any exist, and let
$h_k \equiv g_{m(k)+1}$, $k=0,\dots,\kappa$.
%The number $\kappa = \kappa(f)$ 
%is called the {\it height } of $f$. In the immediately renormalizable
%case set $\kappa = -1$.

The nest of intervals (or the corresponding nest of pieces $V^m$)
\begin{equation}\label{cascade}
I^{m(k)+1}\supset I^{m(k)+2}\supset\ldots\supset I^{m(k+1)}
\end{equation}
is called a {\it central cascade}. The {\it length} $l_k$ of the
cascade is defined as $m(k+1)-m(k)$. Note that a cascade of length 1
corresponds to a non-central return to level $m(k)$.

 A cascade \ref{cascade} is called
{\it saddle-node} if $0 \not\in h_k I^{m(k)+1}$.
 Otherwise it is called
{\it Ulam-Neumann}. For a long saddle-node cascade the map $h_k$
 is combinatorially close
to $z\mapsto z^2+1/4$.  For  a long Ulam-Neumann cascade it is close to 
$z\mapsto z^2-2$. 

The next lemma shows that for a long saddle-node cascade, the map 
$h_{k}:I^{m(k)+1}\rightarrow I^{m(k)}$ is  a 
small perturbation of a map with a parabolic fixed point.

\begin{lem}[\cite{L2}]
\label{pert}
Let $h_k:U_k \to V_k$ 
be a sequence of real-symmetric
quadratic-like maps with $\mymod(h_k) \ge \eps > 0$
having saddle-node cascades of length $l_k\to\infty$. Then
any limit point 
of this sequence in the Carath\'eodory topology $f: U \to V$
is hybrid eqivalent to $z \mapsto z^2+1/4$, and thus
has a parabolic fixed point. 
\end{lem}
\begin{pf} 
It takes $l_k$ iterates for the critical point to escape  $U_k$
under iterates of $h_k$. Hence the critical point does not 
escape $U$ under iterates of $f$. By the kneeding theory \cite{MT}
$f$ has on the real line topological type of $z^2+c$ with $-2\leq c\leq 1/4$.
Since small perturbations of $f$ have escaping critical point,
the choice for $c$  boils down to only two boundary 
parameter values, $1/4$ and $-2$.
Since the cascades of $h_k$ are of saddle-node type, $c=1/4$.  
\end{pf}

Since both fixed points of such a sequence $h_k$ are repelling, it follows
from \lemref{Tangent} that for $k$ large enough $h_k$ has perturbed
Fatou coordinates and so $h_k$ has a saddle-node cascade in the
sense described in \secref{sec:ParabolicBackground}. 

Let $x\in  P(f)\cap (I^{m(k)}\setminus I^{m(k)+1})$ and let
$h_k x\in I^j\setminus I^{j+1}$. Set 
$$d(x)=\min\{j-m(k), m(k+1)-j\}.$$
This parameter shows how deep the orbit of $x$ lands inside the cascade. 
Let us  now define 
 $d_k$ as the maximum of  $d(x)$ over all
 $x\in P(f)\cap (I^{m(k)}\setminus I^{m(k)+1})$. 
Given a saddle-node cascade (\ref{cascade}), let us call all levels
$m(k)+d_k< l < m(k+1)-d_k$ {\it neglectable}. 

Let $f$ be renormalizable and $f_1$ a pre-renormalization of $f$.
Define the {\it essential period} $p_e=p_e(f)$ as follows.
Let $p$ be the period of the periodic interval $J=B(f_1)$,
and set $J_k=f^k J$, for $0\leq k\leq p-1$.
Let us remove from the orbit $\{J_k\}_{k=0}^{p-1}$ all intervals 
whose first landing to some $I^{m(k)}$ belongs to a neglectable level,
to obtain a sequence of intervals $\{J_{n_i}\}_{i=1}^m$.
The essential period is the number of intervals which are left,
$p_e(f)=m$.
Note the essential period of a shuffle is well-defined and in this way
we can define $p_e(f)$ for any $f \in Quad$ with real combinatorics.

Let us give some examples of combinatorial types involving long
saddle-node cascades with neglectable levels.
Let $\Gamma$, $\Gamma'$, $\chi$, $\chi'$ and $\chi_0$ be 
from \ref{eq:groups} and \ref{eq:homo}.
\begin{ex}[Goes Through Twice]
\label{ex:casc1}
Let $\chi_2:\Gamma \to \Gamma$ be the homomorphism generated by
$I_0 \mapsto I_0$ and $I_{-1} \mapsto I_{-1}^2 I_0$.
Then any sequence of the form
$$
\Gamma' \upto{\chi'} \Gamma \upto{\chi} \cdots
 \upto{\chi} \Gamma \upto{\chi_2} \Gamma \upto{\chi} \cdots
 \upto{\chi} \Gamma \upto{\chi_0} \Gamma_0
$$
will correspond to a shuffle where the critical orbit moves up through
the cascade until the top,
returns to the level of $\chi_2$, moves up through the cascade again and
then returns to the renormalization interval. If the total number
of levels in the sequence is $m$ then the number of neglectable
levels will be roughly $m - 2 \min(d,m-d)$ where $d$ is the level
of $\chi_2$. 
%Let $\sigma_{m,d}$ denote the corresponding shuffles.
\end{ex}

\begin{ex}[Two Cascades]
\label{ex:twocasc}
As a second example imagine perturbing the right-hand
picture in \figref{blow-up}
so that the renormalization becomes hybrid equivalent to
$z^2 + \frac{1}{4}$. Now any further perturbation will cause
the parabolic orbit to bifurcate and we can create another long cascade.
More specifically, let
$\chi_3: \Gamma \to \Gamma$ be the homomorphism generated by
$I_0 \mapsto I_{-1} I_0$ and 
$I_{-1} \mapsto I_{-1}^2 I_0$ and consider a sequence of the form
$$
\Gamma' \upto{\chi'} \Gamma \upto{\chi} \cdots
 \upto{\chi} \Gamma \upto{\chi_3} \Gamma \upto{\chi} \cdots
 \upto{\chi} \Gamma \upto{\chi_0} \Gamma_0.
$$
Since $\chi_3$ has a non-central return the two long
sequences of $\chi$ form two seperate saddle-node cascades, each
with a long sequence of neglectable levels. 
%Let $\sigma'_{l_1,l_2}$ denote the shuffle with cascades of length
%$l_1$ and $l_2$.
\end{ex}

\subsection{Parabolic shuffles}
\label{sec:shuffles}

Let $\Omega_p$ be the space of shuffles $\sigma$ satisfying
$p_e(\sigma) \le p$.
In this section we construct a compactification 
$\Omega^{cpt}_p$ of $\Omega_p$ which
will form the elements of our combinatorial description of renormalization limits.

Suppose $f \in RQuad$ is renormalizable
and let 
$$
\Gamma_m \upto{\chi_m} \Gamma_{m-1} 
 \upto{\chi_{m-1}} \cdots \upto{\chi_2} \Gamma_1 \upto{\chi_1} \Gamma_0
$$
be its sequence of return types.
Let $l$ be a neglectable level and let 
$\chi_l:(\Gamma_l,\eps_l) \to (\Gamma_{l-1},\eps_{l-1})$ be the return type
of $g_l$. It is clear that if both level $l-1$ and $l+1$ are
neglectable then $(\Gamma_l,\eps_l)$ and $(\Gamma_{l-1},\eps_{l-1})$
are generated by configurations of the form
$$
\pm I_{-p}, \pm I_{-p+1}, \dots, \pm I_{-1}, -I_0
$$ 
or by
$$
+I_0, \pm I_1, \dots, \pm I_{p-1}, \pm I_p
$$ 
for some $p \ge 1$.
We claim that $(\Gamma_1,\eps_l) \cong (\Gamma_{l-1},\eps_{l-1})$
and that $\chi_l$ is defined by
$I_i \mapsto I_iI_0$ for $i \neq 0$ and $I_0 \mapsto I_0$.
First it is clear $I_0 \mapsto I_0$.
Now if $\chi_l(I_i)$ contained more than one off-critical interval
then $l$ would not be a neglectable level. Since $\chi_l$ is unimodal
it follows $\Gamma_{l-1}$ contains at least as many intervals as $\Gamma_l$.
Since the return type sequence is irreducible $\Gamma_{l-1}$
contains exactly the same number of intervals as $\Gamma_l$. Hence 
$I_i \mapsto I_iI_0$. The claim that the signs agree follows from the
condition that $\chi_l$ be admissible.

%In particular, $\chi_l$ is an admissible unimodal homomorphism.
Hence we can ``insert'' another neglectable level into $S$ before $l$
to obtain another irreducible sequence $S'$ of return types:
$$
\Gamma_m \upto{\chi_m} \cdots
 \Gamma_l \upto{\chi_l}  \Gamma_{l-1} \cong \Gamma_l \upto{\chi_l} \Gamma_{l-1}
 \cdots \upto{\chi_1} \Gamma_0.
$$
From \thmref{shuffleexist} there is a unique shuffle $\sigma'$
such that $S(\sigma') = S'$. 

We say two shuffles $\sigma$ and $\sigma'$ in $\Omega_p$ are
{\it essentially equivalent} if one can insert a finite 
number of neglectable levels into $\sigma$ and $\sigma'$ and obtain
equal shuffles. 
Let $\Xi$ be the partition of $\Omega_p$ into
essentially-equivalent equivalence classes. 
Let $U \in \Xi$ be a non-trivial equivalence class. 
Then there is an $n = n_U > 0$ such that for any $\sigma \in U$ 
the return type sequence $S(\sigma)$ has exactly
$n$ different cascades $S_1,S_2,\dots,S_n$, canonically ordered,
containing neglectable levels. Let $l_k$, $k =1,\dots,n$, 
denote the number of neglectable levels in the cascade $S_k$.
The map $\theta_U:U \to \N_+^n$ given by
$\sigma \mapsto (l_1,l_2,\dots,l_n)$
is a homeomorphism. 
Let 
$$
\overline{\N}_+ = \N_+ \union \{+\infty\}
$$
be the one-point compactification of $\N$.
Define $U^{cpt} \supset U$ as the unique space such that
$\theta_U$ extends to a homeomorphism 
$\theta_U:U^{cpt} \to \overline{\N}_+^n$.
Define $\Omega^{cpt}_p \supset \Omega_p$ as the 
union of the trivial classes of $\Xi$ and of the
spaces $U^{cpt}$ for non-trivial $U \in \Xi$.
An element of $\Omega^{cpt}_p \setm \Omega_p$ is called an {\it end}
and can be represented by a ``sequence'' of return types
where infinitely long sequences of neglectable levels are allowed:
$$
\Gamma_m \upto{\chi_m} \cdots \upto{\chi_{l+2}} \Gamma_{l+1} \upto{\chi_{l+1}}
 (\Gamma_l \upto{\chi_l}  \Gamma_{l-1})^\infty \upto{\chi_{l-1}}
\Gamma_{l-2} \upto{\chi_{l-3}} \cdots \upto{\chi_1} \Gamma_0.
$$

The following lemma is evident from the definition of essential period and
$\Omega^{cpt}_p$. 
\begin{lem}
\label{shufflecompact}
For any $p >1$ the space
$\Omega^{cpt}_p$ is metrizable and compact.
\end{lem}

Let ${\cal M}_p = \{M(\sigma)\}_{\sigma \in \Omega_p}$
be the collection of $M$-copies corresponding to $\Omega_p$
and let $\cC_p = \{c(\sigma)\}_{\sigma \in \Omega_p}$
be the corresponding collection of centers.
We now describe the topology of $\cC_p$ and how $\cl(\cC_p)$
compares to $\Omega^{cpt}_p$.
For any $U \in \Xi$ with $n = n_U \ge 1$
let $\cC_U \subset \cC_p$ denote the collection 
of centers of $\{M(\sigma)\}_{\sigma \in U}$. Since $\Xi$ is a finite
partition it suffices to describe the topology of the sets $\cC_U$.
We claim for each non-trivial $U \in \Xi$ there is a homeomorphism
of $\R$ which maps $\cC_U$ to the image of the function
$F:\N_+^n \to \R$ given by
$$
F(x_1,x_2,\dots,x_n) = 
 2^{-x_1} + 2^{-x_1 x_2 - 1} + \cdots + 
2^{-x_1 x_2 \cdots x_n - n + 1}
$$
where $n = n_U$ (see \figref{pic:centers}).
\realfig{pic:centers}{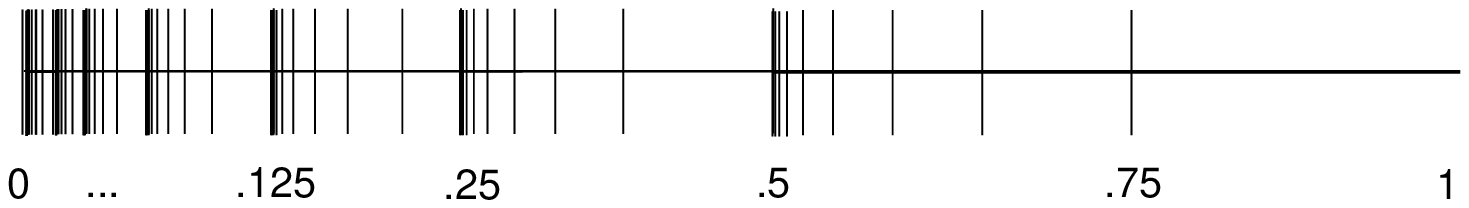}{The image of $F$ for $n=2$.}
{0.8\hsize}

To be more specific the limit points of $\cC_U$ 
are root points of the $M$-copies obtained by ``truncating''
the return type sequences of $\sigma \in U$ at the neglectable
levels.
Let us describe how to truncate a return type sequence
$$
\Gamma_m \upto{\chi_m} \Gamma_{m-1} \upto{\chi_{m-1}} \cdots
\upto{\chi_1} \Gamma_0
$$
at a level $l$.
Let $(\Gamma_T,\eps_T)$ be the semigroup generated by $I_0$
with $\eps_T(I_0) = \eps_l(I^l_0)$
and let $\chi_T$ be the homomorphism defined by
$I_0 \mapsto \chi_l(I^l_0)$. 
Let $S'$ be the sequence
$$
\Gamma_T \upto{\chi_T} \Gamma_{l-1} \upto{\chi_{l-1}} \cdots
\upto{\chi_2} \Gamma_1 \upto{\chi_1} \Gamma_0.
$$
One can check that $S'$ is a sequence
of admissible unimodal return types. If $S'$ is not irreducible then
simply remove all intervals $I^m_i$ not in the combinatorial orbit
of the critical point and shorten the sequence if necessary.
We obtain a unique shuffle 
$\sigma' = \lfloor \sigma \rfloor_l$, the shuffle $\sigma$ truncated at 
level $l$.

Let $U \in \Xi$ satisfy $n = n_U \ge 1$.
Any shuffle $\sigma \in U$ has $n$ cascades with neglectable levels of
lengths $x_1,\dots,x_n$ respectively.
As $x_1 \to \infty$, the corresponding centers accumulate at the root of
the tuned $M$-copy corresponding to any $\sigma \in U$ truncated
at the first neglectable level. If we fix $x_1$ and let
$x_2 \to \infty$ the corresponding centers accumulate at the root of
the $M$-copy corresponding to truncating at the second cascade of
neglectable levels. In general if we
fix the lengths of the first $k$ sequences
of neglectable levels and let the length of the $k+1$-st sequence grow
the centers converge to the root of the $M$-copy corresponding
to truncating at the ($k+1$)-st neglectable sequence.

Given an end $\tau \in \Omega_p^{cpt}$ let
$$
c(\tau) = root(\lfloor \sigma \rfloor_l)
$$ 
where $\sigma \in \Omega_p$ 
is in a sufficiently small neighborhood of $\tau$,
$l$ is a neglectable level of $\sigma$ which belongs to the
first infinitely long cascade of $\tau$, and $root(\sigma)$ is the root
of the $M$-copy $M(\sigma)$.
The map $c:\Omega^{cpt}_p \to \R$ is continuous and its image
is $\cl(\cC_p)$.

We return to our examples. Choose a large $p$ so that
the shuffles from \exref{ex:casc1} and \exref{ex:twocasc}
are contained in $\Omega_p$.

First consider the essentially period tripling shuffles $\sigma^{(3)}_n$.
Then $c(\sigma^{(3)}_n) \to root(\sigma^{(3)})$
where $\sigma^{(3)}$ is the period tripling shuffle.
Moreover, $\sigma^{(3)}_n$ converges to an end $\tau_1 \in \Omega^{cpt}_p$.

Now consider the shuffles $\sigma_{m,d}$ from \exref{ex:casc1}
(Goes Through Twice).
First fix $d > 1$ and let $m \to \infty$. Then 
$c(\sigma_{m,d}) \to root(\lfloor \sigma_{m,d} \rfloor_l) \neq 
root(\sigma^{(3)})$ 
where $l$ is any
negectable level and, in much the same spirit as essential period
tripling, $\sigma_{m,d}$ converges in $\Omega^{cpt}_p$ to an end.
Now fix $m-d > 1$ and let $m \to \infty$. Then 
$c(\sigma_{m,d}) \to root(\sigma^{(3)})$ and
$\sigma_{m,d}$ converges to an end $\tau_2 \in \Omega^{cpt}_p$.
%Since $\tau_1 \neq \tau_2$ it follows $\Omega^{cpt}_p$ is not
%homeomorphic to the closure in $\R$ of 
%the centers $\cC_p$.

Finally consider the shuffles $\sigma_{l_1,l_2}$ from 
\exref{ex:twocasc} (Two Cascades).
Fix $l_1 > 1$ and let $l_2 \to \infty$. Then $c(\sigma_{l_1,l_2})
 \to r_{l_1} = root(\lfloor \sigma_{l_1,l_2} \rfloor_l)$ where $l$
is any negectable level in the second cascade.
The sequence $r_{l_1} \to root(\sigma^{(3)})$ as $l_1 \to \infty$.
Moreover, for any sequence of $l_2$ if we let $l_1 \to \infty$
then $c(\sigma_{l_1,l_2}) \to root(\sigma^{(3)})$.
Now consider the limits of $\sigma_{l_1,l_2}$ in $\Omega^{cpt}_p$.
If we fix $l_2$ and let $l_1 \to \infty$ the shuffles will
converge to an end $\tau_{\infty,l_2}$.
%Since all of the ends $\tau_{\infty,l_2}$ are distinct
%it follows the compactification of 
%the equivalence class $U$ containing the
%collection $\sigma'_{l_1,l_2}$ is not
%homeomorphic to the closure in $\R$ of 
%the centers $\cC_U$.

This completes our description of the topology of ${\cal M}_p$
and how $\cl(\cC_p)$ compares to $\Omega^{cpt}_p$.
\section{Parabolic Renormalization}
\label{sec:parabRenorm}

Let $c \star M$ be a maximal tuned Mandelbrot set 
with root $c'$ and suppose
$f \in \cH(c')$ is renormalizable
and let $f_0$ be a pre-renormalization.
Let $\xi = \beta(f_0)$. 
Choose incoming and outgoing petals
$D_{\pm}$ around the parabolic point $\xi$ and
let $\cC_{\pm}$ denote the respective \`Ecalle-Voronin cylinders
and $\pi_{\pm}$ the projections with $\pi_{+}$ extended to 
$B = \myint(K(f_0))$.
Fix a transit map $g:\cC_+ \to \cC_-$
satisfying 
$$
g(\pi_{+}(0)) \not\in \pi_{-}(K(f_0)).
$$

Given a collection 
$\{f_\alpha\}$ of homolorphic maps let
$\seq{f_\alpha}$ denote the set of 
restrictions of all finite compositions
of $\{f_\alpha\}$.
Let 
$$
\cF(f,g) = \seq{f \union \{\mbox {all local lifts of }g 
	\mbox{ to } D_\pm\}}.
$$
Note that $\cF(f,g)$ is independent of the choice of petals $D_\pm$.
%and we will often abuse notation by writing $g(z)$ when we mean 
%$\wtl{g}(z)$ for some local lift $\wtl{g}$ of $g$ to some petals $D_\pm$.
A collection $\cF$ of holomorphic maps closed
under composition and restriction is called a 
{\it conformal dynamical system}.
Define the orbit of a point $z \in \C$ as
$$
\orb(z) = \orb(\cF,z) = \bigcup_{h \in \cF} h(z).
$$
We say $\cF$ is contained in any {\it geometric limit} of
a sequence $\cF_n$ if for any $f \in \cF$ there are $f_n \in \cF_n$
such that $f_n \to f$ on compact sets.

We say the pair $(f,g)$ is {\it parabolic renormalizable} if
there is a neighborhood $U \ni 0$ and a
map $h \in (\cF(f,g) \setm \seq{f})$ such that 
$$
h|_U \in Quad.
$$
%It is clear there is a minimal such $h$ in the sense that
%if $h_1|_{U_1}$ is another such map then there is an $n \ge 1$
%such that $h^n = h_1$ on a neighborhood of the origin.
We call such an $h|_U$ a {\it parabolic pre-renormalization}
of $(f,g)$
and we call the germ of a normalized pre-renormalization
a {\it parabolic renormalization} of $(f,g)$.
%denoted $\cR(f,g)$.
In the next section we will show that the domain $U$ of the
pre-renormalization can be canonically chosen.

\subsection{Essentially period tripling}
\label{sec:parabRenormTriple}

In this section we describe a construction from \cite{DD}
for finding a canonical representation of 
the parabolic renormalization
in the essentially period tripling case.
For simplicity we will state the construction for the quadratic map
$P_{-1.75}$. However, it is clear how to generalize this construction
to any map $f \in \cH(-1.75)$.

Recall from \secref{sec:special} the sequence of maximal tuned Mandelbrot sets
$c_n \star M$ with essentially period tripling combinatorics accumulate
at the root of the period three tuned copy, $c = -1.75$.
Let $f = P_{-1.75}$ and choose $f_0$ and $D_\pm$ as above. Let 
$B = \myint(K(f_0))$ and let $f_n = P_{c_n}$.
Choose $n_0$ sufficiently large and choose 
$$
U_- \in \bigcap_{n \ge n_0} Comp(\myint(K(f)) \setm B, P(f_n))
$$
such that $U_- \Subset D_-$.
Let $t$ be the landing time of $U_-$ to $B$ under $f$.

Let 
$$
\cD_f = \{g:\cC_+ \to \cC_- | g \mbox{ is a transit map and }
	g(\pi_+(0)) \in \pi_-(U_-)\}.
$$
The phase map gives a conformal isomorphism of
$\cD_f$ to a disk $D_f \Subset \C / \Z$.
Note that $D_f$ is a Jordan domain.
Choose a branch of $\pi_-^{-1}$ so that $Range(\pi_-^{-1}) \supset U_-$.
For $g \in \cD_f$ let $W_g$ be the connected component
of $(\pi_+^{-1} \circ g^{-1} \circ \pi_-)(U_-)$ containing $0$.
Since $\pi_-(U_-)$ is a topological disk, it follows
the map $R_{\bar{a}}:W_g \to B$ given by
$$
	R_{\bar{a}} = f^t \circ \pi_-^{-1} \circ g_{\bar{a}} \circ \pi_+
$$
is quadratic-like with possibly disconnected Julia set.
If $J(R_{\bar{a}})$ is connected then we have
constructed a parabolic pre-renormalization of $(f,g)$.

Fix any $\ast \in D_f$.
Define the holomorphic motion 
$$
h_{\bar{a}}: (\bd B,\bd W_{g_\ast}) \to (\bd B,\bd W_{g_{\bar{a}}})
$$ 
on $\bd B$ by the identity
and locally on $\bd W_{g_{\bar{a}}}$ by pulling back
under $R_{\bar{a}}$.
%$$
%\pi_+^{-1} \circ g_{\bar{a}}^{-1} \circ g_\ast \circ \pi_+
%$$
%for an appropriate branch of $\pi_+^{-1}$.
Let $\bV = \{(\bar{a},z): \bar{a} \in D_f, z \in B\}$ and
$\bU = \{(\bar{a},z): \bar{a} \in D_f, z \in W_{g_{\bar{a}}}\}$.
Let $\bff:\bU \to \bV$ be defined by 
$$
\bR(\bar{a},z) = (\bar{a},R_{\bar{a}}(z)).
$$

\begin{lem}
\label{parabfamily}
The family $(\bR,\bfh)$
is a proper DH quadratic-like family with winding number 1.
\end{lem}
\begin{pf}
The map $f^t$ is a conformal isomorphism of a neighborhood of $U_-$
onto a neighborhood of $B$.
There is a branch of $\pi_-^{-1}$ such that
the map $(\pi_-^{-1} \circ g_{\bar{a}} \circ \pi_+)(0)$
is a conformal isomorphism of a neighborhood of $D_f$
onto a neighborhood of $U_-$.
The lemma follows.
\end{pf}

The following lemma states that the renormalization operators
$\cR_{\sigma^{(3)}_n}$ converge to essentially period tripling
parabolic renormalization.

\begin{lem}
\label{specialconverge}
Let $f \in \cH(-1.75)$.
Suppose $f_k \in Quad$ is a sequence of renormalizable maps
with $f_k \to f$ and $\sigma(f_k) \to \tau$.
Let $g_k:\cC_{f_k,+} \to \cC_{f_k,-}$
be the induced transit maps with phases $\bar{a}_k$.
Then 
\begin{enumerate}
\item $\{ \bar{a}_k \}$ is pre-compact
\item if $\bar{a}_{k_j} \to \bar{a}$ is a convergent subsequence then
$J(R_{\bar{a}})$ is connected and
$$
[h_{k_j}] \to [R_{\bar{a}}]
$$
where $h_k$ is a pre-renormalization of $f_k$
\item $\cF(f,g_{\bar{a}})$ is contained in any
geometric limit of $\seq{f_{k_j}}$
\end{enumerate}
\end{lem}
\begin{pf}
Let $h_k$ be a pre-renormalization of $f_k$.
Since $\sigma(f_k) \to \tau$ and $f_k \to f$ we can write 
\begin{equation}
\label{eq:expr1}
h_k = f_k^{N_1} \circ \wtl{g}_k \circ f_k^{N_2}
\end{equation}
on some neighborhood of the origin
for some fixed $N_1$, $N_2$ and some choice of local lift $\wtl{g}_k$
of the induced transit maps $g_k:\cC_{f_k,+} \to \cC_{f_k,-}$. 
The first claim is that $h_k$ can be chosen in $Quad(m')$ for some $m' > 0$.
Let $V'$ be an $\eps$-neighborhood of the central basin $B$ of $f$
for some small $\eps > 0$. 
Choose $\eps$ small enough and $N_1$ and $N_2$ large enough so that 
for large $k$ the right-hand side of (\ref{eq:expr1}) can be used
to define a pre-renormalization $h_k$ with range $V'$.
Let $U'_k = h_k^{-1}(V')$.
By taking $k$ larger still we can assume $U'_k$ is
contained in an $\eps/2$ neighborhood of $B$.
It follows there is an $m' > 0$ so that $\mymod(U'_k,V') \ge m'$.
Moreover, $\diam(U'_k) \ge C > 0$ for some $C$ independent of $k$.
Hence (\ref{eq:expr1}) holds on a definite neighborhood of the origin.

From the convergence of Fatou coordiantes and the convergence of $f_k$
it follows that
$\{ \bar{a}_k \}$ is pre-compact. Let $\bar{a}_{k_j} \to \bar{a}$
be a convergent subsequence. 
Then $h_{k_j}$ converges on a definite neighborhood of the origin to
the map $f^{N_1} \circ \wtl{g}_{\bar{a}} \circ f^{N_2}$ for an
appropriate local lift $\wtl{g}_{\bar{a}}$ of $g_{\bar{a}}$.
Since the origin
is non-escaping under all $h_k$ it follows $J(R_{\bar{a}})$
is connected.
The last statement follows from the fact that $f_k \to f$ and
\lemref{converge1}.
\end{pf}

Moreover, the proof of the previous lemma can be modified to prove the following
\begin{lem}
\label{specialconverge2}
Suppose $f \in \cH(-1.75)$ and
$f_k \in \cH(-1.75)$ satisfy $f_k \to f$. 
Suppose $g_k:\cC_{f_k,+} \to \cC_{f_k,-}$
is a transit map with phase $\bar{a}$ such that 
$R_k = R_{\bar{a}_k}$ is defined.
Then 
\begin{enumerate}
\item $\{ \bar{a}_k \}$ is pre-compact
\item if $\bar{a}_{k_j} \to \bar{a}$ is a convergent subsequence then
$$
R_{k_j} \to R_{\bar{a}}.
$$
\item $\cF(f,g_{\bar{a}})$ is contained in any
geometric limit of $\cF(f_{k_j},g_{k_j})$
\end{enumerate}
\end{lem}

We finish this section with two useful properties of parabolic renormalization.
The first property is that open
sets intersecting the Julia set of the parabolic pre-renormalization
iterate under $\cF(f,g)$
to open sets intersecting $J(f)$.

\begin{lem}
\label{specialaccum}
Let $f \in \cH(-1.75)$ and $g:\cC_+ \to \cC_-$ be a
transit map with phase $\bar{a}$ such that $J(R_{\bar{a}})$ is connected.
Suppose $U$ is an open set satisfying 
$$
U \inter J(R_{\bar{a}}) \neq \es.
$$
Then there is an 
$h \in \cF(f,g)$ such that $U \inter Dom(h) \neq \es$ and
$$
h(U) \supset J(f).
$$
\end{lem}
\begin{pf}
From the construction of $R_{\bar{a}}$ it is clear that there is an
$h \in \cF(f,g)$ such that
$h$ is a quadratic-like extension of $R_{\bar{a}}$
to a small neighborhood of $B = Range(R_{\bar{a}})$.
It follows that there an $m \ge 0$ 
such that $h^m(U) \inter \bd B \neq \es$.
But $\bd B \subset J(f)$. Iterating $f$ further covers all of $J(f)$.
\end{pf}

The second property is that no quadratic-like representative of
$[R_{\bar{a}}]$ can have too large a domain.

\begin{lem}
\label{specialnest}
Let $f \in \cH(-1.75)$ and $g:\cC_+ \to \cC_-$ be a
transit map with phase $\bar{a}$ such that $R_{\bar{a}}$ is defined.
If $(\wtl{f}:U \to V) \in Quad$ satisfies $[\wtl{f}] = [R_{\bar{a}}]$
then 
$$
U \subset Range(R_{\bar{a}}).
$$
\end{lem}
\begin{pf}
Let $f_1$ be a pre-renormalization of $f$ and let $B = Range(R_{\bar{a}})$.
Suppose $U \inter \bd B \neq \es$ and
let $U'$ be the connected component of $U \inter B$ containing $0$. 
Since $f_1$-preimages of $0$ accumulate on $J(f_1) = \bd B$
there exists an $n > 0$ and $z_0 \in U'$ such that
$f_1^n(z_0) = 0$.
Since $[\wtl{f}] = [R_{\bar{a}}]$ it follows $\wtl{f}$ has a critical
point at $z_0$, which is a contradiction.
\end{pf}

\subsection{Generalized parabolic renormalization}
\label{sec:genParabolicRenorm}

In this section we modify the
construction of parabolic renormalization to 
act on generalized quadratic-like maps.

Let $f:\union_j U_j \to V$ be a generalized quadratic-like map
with $f_0 = f|_{U_0} \in \cH(\qtr)$.
Let $\xi = \beta(f_0)$.
Choose incoming and outgoing petals
$D_\pm$ around the parabolic point $\xi$ and
let ${\cal C}_\pm$ denote the respective \`Ecalle-Voronin cylinders
and $\pi_\pm$ the projections with $\pi_+$ extended to 
$B = \myint(K(f_0))$.

For a given $g:\cC_- \to \cC_+$
let $L_0$ be the first landing map under 
$\cF(f,g)$ to $\union_{j \neq 0} U_j$.
Note that if $C$ is a connected component of $Dom(L_0)$
then there is an $h \in \cF(f,g)$ 
such that $C \Subset Dom(h)$ and $h(z) = L_0(z)$ for
all $z \in C$. 
Let $T$ be the {\it first through map} 
\begin{equation}
\label{eqn:Rdef1}
T = f \circ L_0.
\end{equation}
Note that $T$ is not a generalized quadratic-like map. However,
$T$ has at most one critical value, and, with a slight abuse 
of notation we will treat $T$ as a generalized quadratic-like map.

\realfig{yetanother}{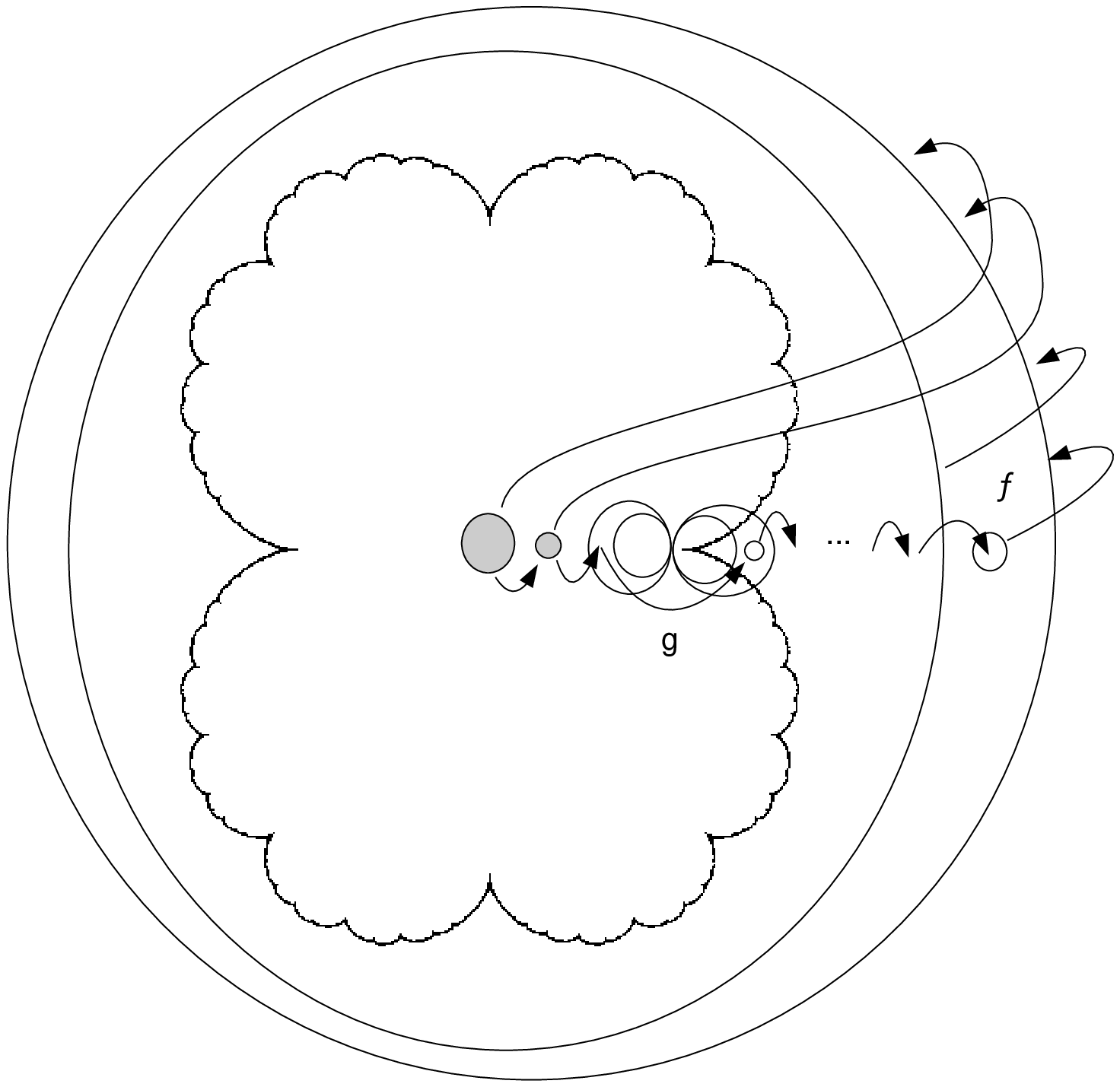}{The first through map $T$
with some components of the domain shaded.}{0.6\hsize}

Let $X \subset \C / {\Bbb Z}$ be the set of phases $\bar{a}$ such that
for $g = g_{\bar{a}}$,
$$
0 \in Dom(L_0).
$$
It is clear that $X$ is a countable pairwise disjoint collection of 
Jordan disks. Let $D$ be a connected component of $X$.
Let $\bT$ denote the family over $D$ of first through maps $T$.
The construction of the holomorphic motion $\bfh$ described before
\lemref{parabfamily} carries over unchanged to this situation.
Moreover, one can modify the proof of \lemref{parabfamily} to prove
\begin{lem}
For any connected component $D$ of $X$, 
the family $(\bT,\bfh)$ over $D$
is a proper generalized quadratic-like 
family with winding number 1.
\end{lem}

The following two lemmas are generalizations of \lemref{specialconverge}
and \lemref{specialconverge2}, respectively. We omit the proofs.
\begin{lem}
\label{genapprox1}
Suppose $f \in Gen$ satisfies $f|_{U_0} \in \cH(\qtr)$
and $f_k \in Gen$ is a sequence converging to $f$.
%such that both fixed points of $f_k|_{U_{k,0}}$ are repelling.
Let $T_k$ be the first through map for $f_k$ and suppose 
that $T_k \in Gen$.
Let $g_k$ be the induced transit maps of $f_k$ with phase $\bar{a}_k$.
Suppose $\bar{a}_{k_j} \to \bar{a}$ is a convergent subsequence
and suppose $0 \in Dom(T_{\bar{a}})$
where $T_{\bar{a}}$ is the first through map for $\cF(f,g_{\bar{a}})$,
Then
$$
T_{k_j} \to T_{\bar{a}}
$$ 
Moreover, $\cF(f,g_{\bar{a}})$ is contained in any
geometric limit of $\seq{f_{k_j}}$.
\end{lem}

\begin{lem}
\label{genapprox2}
Suppose $f \in Gen$ 
and $f_k \in Gen$ satisfy $f|_{U_0} \in \cH(\qtr)$,
$f_k|_{U_{k,0}} \in \cH(\qtr)$, and $f_k \to f$.
Let $g_k$ be a transit map of $f_k$ with phases $\bar{a}_k$.
Let $T_k$ be the first through map for $\cF(f_k,g_k)$ and suppose 
$T_k \in Gen$.
Suppose $\bar{a}_{k_j} \to \bar{a}$ is a convergent subsequence
and suppose $0 \in Dom(T_{\bar{a}})$
where $T_{\bar{a}}$ is the first through map for $\cF(f,g_{\bar{a}})$,
Then
$$
T_{k_j} \to T_{\bar{a}}
$$ 
Moreover, $\cF(f,g_{\bar{a}})$ is contained in any
geometric limit of $\cF(f_{k_j},g_{k_j})$.
\end{lem}

\section{Towers}
\label{sec:towerdef}

Let $S \subset \Z$ be a set of consecutive integers containing $\N_0$
and let $f_n$ be a sequence of maps in $Gen$ indexed by $n \in S$.
Let $U_{n,0}$ be the central component of $f_n$. Let
$$
S_\cC = \{n \in S : f_n|_{U_{n,0}} \in \cH(\qtr)\}
$$
and let 
$$
S_Q = \{n \in S : f_n \in Quad \}.
$$ 
For $n \in S_\cC$ let 
$g_n:\cC_{n,+} \to \cC_{n,-}$ be a transit map between the \`Ecalle-Voronin
cylinders $\cC_{n,\pm}$ of $f_n$.
The collection of maps 
$$
\cT = \{f_n: n \in S\} \union \{g_n: n \in S_\cC\}
$$ 
is called a {\it tower} iff
for each pair $n, n+1 \in S$ one of the following conditions hold:
\begin{enumerate}
\item[T1:] $n \in S_Q$, $f_n$ is immediately renormalizable
	and $[f_{n+1}] = [h]$ where $h$ is
	a pre-renormalization of $f_n$ of minimal period
\item[T2:] $n \in S_Q$, $f_n$ is not immediately renormalizable
	and $[f_{n+1}] = [h]$ where $h$ is a
	restriction of the first return map to the initial central puzzle piece 
	of $f_n$
\item[T3:] $n \not\in (S_Q \union S_\cC)$ and $[f_{n+1}]=[h]$ 
	where $h$ is a restriction of the first return map $R(f_n,U_{n,0})$
	or first through map $T(f_n,\union_{j \neq 0} U_{n,j})$
\item[T4:] $n \in S_\cC$ and $[f_{n+1}]=[h]$ 
	where $h$ is a restriction of the first through map
	of the pair $(f_n,g_n)$.
\end{enumerate}
We shall often identify $g_n$ with the set of local lifts of $g_n$ for some
choice of incoming and outgoing petals $D_{n,\pm}$.
If $S_\cC \neq \es$ then $\cT$ is a {\it parabolic tower}.

Let $\Tow$ be the space of towers with the following topology:
a sequence $\cT_m = \{f_{m,n},g_{m,n}\}$ converges to 
$\cT=\{f_n,g_n\}$ iff
\begin{itemize}
\item $S_m \to S$ and $S_{m,\cC} \to S'_\cC \subset S_\cC$
\item if $n \in S \setm S_\cC$ then $f_{m,n} \to f_n$ 
\item if $n \in S'_\cC$ then $f_{m,n} \to f_n$ and $g_{m,n} \to g_n$
\item if  $n \in S_\cC \setm S'_\cC$ then $f_{m,n}|_{U_{m,n,0}}$ has both
fixed points repelling, $f_{m,n} \to f_n$ and $h_{m,n} \to g_n$ 
where $h_{m,n}$ is the induced transit map on the perturbed \`Ecalle-Voronin
cylinders.
\end{itemize}

If $S = {\Bbb Z}$ then $\cT$ is a {\it bi-infinite tower} and otherwise
$\cT$ is a {\it forward tower}.
%A tower with $S = \{m' : m' \ge m\}$ for some
%$m \in {\Bbb Z}$ is a {\it forward tower}.
%For $k \in \Z \union \{-\infty\}$ let 
%$\Tow^k \subset \Tow$ be the subspace of towers
%whose index set contains $\{m \in \Z: m \ge k\}$.
The map $f_{min(S)}$ in a forward tower is called the {\it base map}.
Define $Dom(\cT)$ and $Range(\cT)$ to be the domain and range of the base map.
%A tower is {\it normalized} if $f_0$ is normalized.

Let $\cT$ be a forward tower and let $f_m$ be the base map of $\cT$.
Let
$$
\cF(\cT) = \seq{\cT \setm \{f_n : n > m\}}
$$
where recall $\seq{f_\alpha}$ denotes the set of 
restrictions of all finite compositions
of $\{f_\alpha\}$.
Define the orbit of $z \in Dom(\cT)$ by 
$$
\orb(z) = \orb(\cF(\cT),z).
$$
Note that if $Range(f_n) \subset Range(f_m)$ and $z \in Dom(f_n)$
then $f_n(z) \in \orb(z)$.
We say $\orb(z)$ {\it escapes} if 
$\orb(z) \inter (Range(\cT) \setm Dom(\cT)) \neq \es$.
Define the {\it filled Julia set}, $K(\cT)$,
the {\it Julia set}, $J(\cT)$,
and the {\it post-critical set},
$P(\cT)$, as for quadratic-like maps.
For a bi-infinite tower $\cT$ define the post-critical set
$$
P(\cT) = \cl\bigcup_{S' \subset S} P(\cT|_{S'})
$$
where $\cT|_{S'}$ ranges over forward subtowers of $\cT$ and where
the closure is taken as a subset of $\Chat$.

Two towers $\cT$ and $\cT'$ with $S(\cT) = S(\cT')$
are {\it quasi-conformally equivalent} if
there is a quasi-conformal map $\phi$ such that 
\begin{enumerate}
	\item $\phi$ is a quasi-conformal conjugacy of $f_n$ and $f'_n$ on a
		neighborhood of $K(f_n)$ to a neighborhood of $K(f'_n)$ for all
		$n \in S$,
	\item $\phi$ induces a quasi-conformal conjugacy 
		of the transit maps $g_n$ and $g'_n$ for 
		$n \in S_\cC$.
\end{enumerate}
A quasi-conformal equivalence $\phi$ between two forward towers
is a {\it hybrid} equivalence if
$\bar\partial \phi|_{K(\cT)} \equiv 0$ and is a {\it holomorphic}
equivalence if $\phi$ is holomorphic.
%Two towers $\cT$ and $\cT'$ are {\it germ-equivalent} if
%\begin{enumerate}
%	\item $[f_n] = [f'_n]$ for all $n \in S$,
%	\item $g_n = g'_n$ for $n \in S_\cC$.
%\end{enumerate}
%A measurable Beltrami differential
%$\mu = u(z) d\bar{z}/dz$ is {\it invariant} under a forward tower $\cT$
%if $h^*(\mu) = \mu$ almost everywhere for any $h \in \cT$.
The following proposition is the analogue of \propref{Straighten}
for towers.

\begin{prop}[Straightening]
\label{forStraight}
Let $\cT$ be a forward tower such that its base map
is quadratic-like. Then $\cT$ is hybird equivalent
to a tower with a quadratic base map.
%The dilatation of $\phi$ is bounded above in terms of $\mymod(f_m)$.
%Moreover, $\phi$ is a conjugacy on an 
%$\eps(\mymod(f_m))$-scaled neighborhood of $K(f_m)$.
\end{prop}
\begin{pf}
Let $f_m$ be the base map of $\cT$.
From \propref{Straighten} there is a hybrid equivalence $\phi$
between $f_m$ and a unique polynomial of the form $z^2+c$.
Let $u(z)$ be the complex dilatation of $\phi$ and
let $\mu = u(z) d\bar{z}/dz$ be the corresponding 
Beltrami differential. Since $\phi$ is quasi-conformal there is 
a $k < 1$ such that $\|u(z)\|_{\infty} \le k$.
Let $U \supset K(f_m)$ be the domain on which $\phi$ is a conjugacy.

Define the Beltrami differential $\mu'$ by
$$
\mu'|_{K(\cT)} \equiv 0
$$
and if $z \in (U \setm K(\cT))$ by
$$
\mu'|_{U'} = h^*(\mu)
$$
where $h \in \cF(\cT)$
and $U' \ni z$ satisfy
$h(U') \subset (U \setm K(f_m))$.
There are restrictions $f'_n$ of $f_n$
such that $[f'_n] = [f_n]$ and $\mu'$ is invariant under
the foward tower $\cT' = \{f'_n,g_n\}$.
%If $n \in S_\cC$
%one can choose $f'_{n+1}$ to be any restriction of $f_{n+1}$
%such that 
%$$
%Dom(f'_{n+1}) \subset Dom(\cP\R_{\psh_n}(f_n,I(f_{n+1}))).
%$$

%If $n-1 \in S_\cC$ let $U'_n = f_n^{-1}(U_n)$ 
%and $f'_n = f_n|_{U'_n}$ and otherwise let $f'_n$ be a restriction of
%$f_n$ so that 
%$Range(f'_n) \subset U'_{m(n)}$.
%From the nesting of \lemref{gennest}, 
%if $n-1 \in S_\cC$,
%$f'_n$ can be expressed as a restriction of a composition of
%$f'_{n-1}$ and $\wtl{g}_{n-1,i}$ for $i=1,\dots,|L_{n-1}|$
%and otherwise $f'_n$ is a restriction of an iterate of $f'_{m(n)}$.
%Let $\F'$ be a forward tower equivalent to $\F$
%with quadratic-like part $\seq{f'}$.
%From the nesting properties of $\seq{f'}$
%it follows $\mu'$ is invariant under $\F'$.

Write $\mu'(z)=u'(z)d\bar{z}/dz$.
Since all maps in $\cT'$ are holomorphic $\|u'(z)\|_{\infty} \le k < 1$.
Let $\phi_1$ be the solution to the Beltrami equation
$$
\bar{\partial}\phi_1 = u' \cdot \partial \phi_1
$$
and let 
$$
\cT'' = \{\phi_1 \circ h \circ \phi_1^{-1}: h \in \cT'\}.
$$

We claim $\cT''$ is again a forward tower and
that $\phi_1$ is a hybrid equivalence between $\cT$ and 
$\cT''$. Let $n \in S_\cC$ and $g_n \in \cT'$. Let 
$g_n'' = \phi_1 \circ g_n \circ \phi_1^{-1}$ and
$f_n'' = \phi_1 \circ f'_n \circ \phi_1^{-1}$. Since
$\phi_1$ conjugates forward and backward orbits of $f_n$
to orbits of $f''_n$, it follows that 
$g''_n$ is a map on the \`Ecalle-Voronin cylinders
of $f''_n$. Since $\phi_1$ is a homeomorphism, it is evident that
$g''_n$ is a homeomorphism.
Moreover, $\mu'$ is invariant under $g_n$, and so
$g''_n$ is conformal.
That is, the conjugate of a transit map in $\cT$ is a
transit map in $\cT''$.
The other properties of a tower are clear.

The base map of $\cT''$ is 
holomorphically equivalent to a polynomial. Hence
$\cT''$ is holomorphically equivalent to a
tower with a polynomial base map.
\end{pf}

\section{Limiting Towers}
\label{sec:genfortowers}

In this paper we study the parabolic towers that are
limits of certain McMullen towers. To be precise we
make the following definition. For a given $\kappa > 0$ let
$\Tow(\kappa)$ be the closure of the set $\Tow_0(\kappa)$ 
of $\cT \in \Tow$ satisfying
\begin{enumerate}
\item $S_\cC = \es$
\item $f_0 \in Quad$ is normalized
\item $f_n \in Gen(1/\kappa)$ for all $n \in S$
\item if $n \in S_Q$ then $f_n$ has real combinatorics and $p_e(f_n) \le \kappa$
\item if $n \not\in S_Q$ then $[f_n] = [h_n]$ where
$h_n$ is either the generalized renormalization of $f_{n-1}$ or 
the first through map $T$ of $f_{n-1}$ restricted to $Comp(Dom(T),P(f_{n-1}))$
\item if $f_n$ is a first through map then $f_{n-1}$ has a saddle-node 
cascade in the sense described in \secref{sec:ParabolicBackground}
\item quadratic-like levels are at most $\kappa$ apart: if $n,m \in S_Q$
are adjacent quadratic-like levels then $|m-n| \le \kappa$.
\item $V_n = Range(f_n)$ is a $\kappa$-quasidisk
\item $\diam V_n \le \kappa \diam K(f_n)$
\item Unbranched Property: $V_n \inter P(f_m) = P(f_n)$ for $n \ge m$
\end{enumerate}
We will refer to towers in $\Tow(\kappa)$ as towers with
{\it essentially bounded combinatorics and complex bounds}.
Over the next several sections we will analyze the basic properties
of towers in $\Tow(\kappa)$.

The combinatorics of a tower $\cT \in \Tow_0(\kappa)$
is the sequence $\bar{\sigma}(\cT)$
indexed by $n \in S_Q$ of shuffles $\sigma_n = \sigma(f_n)$.
Recall $\Omega^{cpt}(\kappa)$ is the compactification of
the space $\Omega(\kappa)$ of shuffles $\sigma$ with
$p_e(\sigma) \le \kappa$.
Suppose $\cT_m \in \Tow_0(\kappa)$ is a sequence of towers
converging to the parabolic tower $\cT \in \Tow(\kappa)$.
The combinatorics of $\cT$ is the 
sequence $\bar{\sigma}(\cT)$ indexed by $n \in S_Q$ of
shuffles and ends given by
$\lim_{m \to \infty} \sigma_{m,n}$.
Clearly the combinatorics of a tower is invariant under hyrbid equivalence.

Two towers $\cT=\{f_n,g_n\}$ and $\cT'=\{f'_n,g'_n\}$ 
are {\it combinatorially equivalent} if
$S(\cT) = S(\cT')$ and $\bar{\sigma}(\cT) = \bar{\sigma}(\cT')$.

\begin{prop}[Forward Combinatorial Rigidity]
\label{genCombRigid}
Let $\cT$ and $\cT'$ be forward towers hyrbid equivalent
to towers in $\Tow(\kappa)$.
Let $f_m$ and $f_m'$ be the respective base maps.
Suppose $\cT$ is combinatorially equivalent to $\cT'$ and $[f_m] = [f'_m]$.
Then $[f_n] = [f'_n]$ for all $n \in S$ and
$g_n = g'_n$ for $n \in S_\cC$.
\end{prop}
\begin{pf}
First, it is clear that
$[f_n] = [f'_n]$ for $m \le n \le \min\{S_\cC\}$ where 
$\min\{\es\} = \infty$.
Now suppose by induction that
$n \in S_\cC$ and $[f_n] = [f'_n]$. We claim $g_n = g'_n$.
Let $L_{\bar{a}}$ be the first landing map of 
$\cF(f_n,g_{\bar{a}})$ to the off-critical pieces of $f_n$
and let
$$
X = \{\bar{a}:0 \in Dom(L_{\bar{a}})\}.
$$

Let $(\bT,\bfh)$ be the holomorphic family over
the component $D \subset X$ containing $g_n$ of generalized
quadratic-like maps constructed in \secref{sec:genParabolicRenorm}.
Construct the sequence of families of first return maps
as described in \secref{sec:family}
until the next level where $f_n|_{U_{n,0}} \in Quad$. 
Similarly construct the families containing $g'_n$ using $f'_n$.
Since $\cT$ is combinatorially equivalent to $\cT'$ 
it follows from \thmref{hyper-dense},
uniqueness of root points
and \thmref{DHFamily} that $g_n = g_n'$.
\end{pf}

Combining this result with straightening we have the following
\begin{cor}
\label{CombIsHybrid}
Any two combinatorially equivalent forward towers 
$\cT, \cT' \in \Tow(\kappa)$ are hybrid equivalent.
\end{cor}
\begin{pf}
Straighten $\cT$ and $\cT'$ to the towers $\cT_1$ and $\cT_2$
with quadratic base maps. Since $\cT_1$ and $\cT_2$ are
combinatorially equivalent it follows 
from \thmref{hyper-dense} and the uniqueness of root points that 
the base maps are equal.
Hence by \propref{genCombRigid} $\cT_1$ and $\cT_2$ are
hybrid equivalent.
\end{pf}

We now prove compactness:

\begin{prop}
\label{TowCompact}
For any $\kappa > 0$ the space $\Tow(\kappa)$ is compact.
\end{prop}
\begin{pf}
Let $\cT_m = \{f_{m,n},g_{m,n}\}$ be a sequence in $\Tow(\kappa)$.
By selecting a subsequence
we may assume the index set $S(\cT_m)$ converges to some index set $S$.
If $f_{m,n}$ is a first through map then the modified landing times
$l(z)$ are bounded for all $z \in u_{f_{m,n}}$ since the essential
period is bounded.
By \lemref{returngeo}, \lemref{throughgeo} and \lemref{initialgeo},
there exists a function $C(\kappa)$ such that for all $f_{m,n}$,
$$
\geo(f_{m,n}) \ge C(\kappa) > 0.
$$
From \lemref{genCompact} we can select a subsequence $\cT_{m_k}$
so that $f_{m_k,n}$ converges on all levels $n \in S$ to some 
generalized quadratic-like maps $f_n$.
Let $S_\cC \subset S$ be the levels with $f_n|_{U_{n,0}} \in \cH(\qtr)$.
From \lemref{throughgeo} we can choose a subsequence so that the
transit maps on each level $n \in S_\cC$ converge.
From \lemref{genapprox1} and \lemref{genapprox2} the limiting collection
of maps will form a tower. The other properties of a tower are clear.
\end{pf}

\begin{lem}
\label{scaling}
Let $\cT \in \Tow(\kappa)$. Then $\diam K(f_n) \to 0$ as $n \to \infty$.
If $\cT$ is a bi-infinite tower then $\diam K(f_n) \to \infty$
as $n \to -\infty$.
\end{lem}
\begin{pf}
Let us prove the first statement.
We can assume there are an infinite number of levels $n \to \infty$ where
$f_n$ is not immediately renormalizable, for otherwise $\cT$
is eventually
a McMullen tower with period-doubling combinatorics and the result follows.
Choose a subsequence $f_{n_k}$, $n_k \to \infty$,
of generalized quadratic-like with at least one off-critical piece.

Suppose by contradiction that $\diam K(f_{n_k}) \ge \eps > 0$.
Let $\union_j U_{k,j} = Dom(f_{n_k})$ and $K_{k,j} = K(f_{n_k}) \inter U_{k,j}$.
We may assume $K_{k+1,j} \subset K_{k,0}$ by selecting levels of first return.

Then since $\geo(f_{n_k}) \ge C(\kappa) > 0$ and 
$\mymod(K_{k,j},U_{k,j}) \ge 1/\kappa$ it follows that $U_{k,j}$ contains
a definite neighborhood of $K_{k,j}$.
Hence there is eventually some $j_1, j_2 \neq 0$ and $k_2 > k_1$ with
$K_{k_2,j_2} \inter U_{k_1,j_1} \neq \es$.
But this is a contradiction since $K_{k_2,j_2} \subset K_{k_1,0}$
and $K_{k_1,0} \inter U_{k_1,j_1} = \es$.

The second statement is analogous.
\end{pf}

\begin{prop}[{\cite[Corollary 5.12]{McM2}}]
\label{Pconv}
The postcritical set
$P(\cT)$ varies continuously with $\cT \in \Tow(\kappa)$.
\end{prop}
\begin{pf}
Let $\cT_m$ be a sequence of towers in $\Tow(\kappa)$ converging 
to a tower $\cT$. Assume $\cT$ is a forward tower.
If $z \in \orb(\cT,0)$
then $d(z,P(\cT_m)) \to 0$ as $m \to \infty$
since $\cF(\cT)$ is contained in any geometric limit of $\cT_m$.
Hence $P(\cT) \subset \liminf_m P(\cT_m)$. We must show
$\limsup_m P(\cT_m) \subset P(\cT)$.

For $n \in S_Q$ let $K_n(0) = K(f_n)$ and
let $K_n(i)$ enumerate the orbit of $K(f_n)$ by $\cT$.
That is, 
$$
\union_i K_n(i) = \{h(z) : z \in K(f_n), h \in \cF(\cT)\}.
$$
Let $\delta_n = \sup_i \diam K_n(i)$. The arguments proving
$\diam K_n(0) \to 0$ can be adapted to prove $\delta_n \to 0$.
Let $\eps > 0$ and let $N$ be large enough so that $\delta_N < \eps$.
Let 
$$
\union_i K_{m,n}(i) = \{h(z) : z \in K(f_{m,n}), h \in \cF(\cT_m)\}.
$$
Since $\cT_m \to \cT$ it follows that
for $m > N$ large enough $\union_i K_{m,n}(i)$ 
is contained in an $\eps$-neighborhood of $\union_i K_n(i)$.
Hence $P(\cT_m)$ is contained in a $2\eps$-neighborhood of $P(\cT)$.

Now suppose $\cT$ is a bi-infinite tower. From the 
continuity of $P(\cT)$ for forward towers and the unbranched property
for $\Tow_0(\kappa)$ it follows that 
$P(\cT|_{S_n}) = V_n \inter P(\cT)$ where $S_n \subset S$ is
any index set of a forward tower.
Since $V_n$ contains an $\eps(\kappa)$-scaled neighborhood of
$K(f_n)$ and $\diam K(f_n) \to \infty$
as $n \to -\infty$, it follows
that $P(\cT) = \{\infty\} \union_{S' \subset S} P(\cT|_{S'})$.
\end{pf}

\subsection{Expansion of the hyperbolic metric}
\label{sec:ForwardExpands}

One of the central ideas in McMullen's arguments is that maps in a tower
expand the hyperbolic metric on the complement of the post-critical set.
In this section we prove similar propositions.

\begin{lem}
\label{LowExpand}
There are continuous increasing functions $C_1(s)$ and $C_2(s)$ such that 
if $f:X \hookrightarrow Y$ is an inclusion between two hyperbolic
Riemann surfaces and $x \in X$ then,
letting $s = d(x,Y \setm X)$,
$$0 < C_1(s) \le \| Df(x) \| \le C_2(s) < 1.$$
Moreover, $C_2(s) \to 0$ as $s \to 0$.
\end{lem}
\begin{pf}
The inequality
$\| Df(x) \| \le C_2(s) < 1$ and the properties of $C_2(s)$
are found in \cite{McM2}.
Lift $f$ to the universal cover $\pi:\D \to Y$ 
and normalize so that $x = f(x) = 0$. 
The inclusion $B_s \equiv \{z:d_{\D}(0,z) < s\} 
\hookrightarrow \D$ factors through $f$
and so 
$\| Df(0) \| \ge 1/r(s)$ where $r(s)$ is the radius of $B_s$
measured in the euclidean metric.
\end{pf}

The following Proposition states when 
maps in a forward tower $\cT \in \Tow(\kappa)$
expand the hyperbolic metric on 
$Range(\cT) \setm P(\cT)$ and
gives an estimate on the amount of expansion and the variation of expansion.

Recall if the base map of $\cT$ is $f_m:U_m \to V_m$ then
$Range(\cT) = V_m$ and $Dom(\cT) = U_m$.
We will use the notation $\rho_m$,
$\| \cdot \|_m$, $d_m(\cdot,\cdot)$ and $\ell_m(\cdot)$
to denote the hyperbolic metric, norm, distance and length on
$Range(\cT) \setm P(\cT)$.

\begin{prop}
\label{expand1}
Let $\cT \in \Tow(\kappa)$ be a forward tower with base map $f_m:U_m \to V_m$.
Suppose that $h \in \cF(\cT)$
and let $Q_h = h^{-1}(P(\cT))$. Then
$$
\| Dh(z) \|_m > 1
$$
for any $z \in (Dom(h) \setm Q_h)$.
Moreover, if $(Q_h \setm P(\cT)) \neq \es$ then
$$
C_2^{-1}(s_2) \le 
  \| Dh(z) \|_m \le C_1^{-1}(s_1)
$$
where $s_1 = d_m(z,Q_h \union \bd Dom(h))$
and $s_2 = d_m(z,Q_h)$.
Finally, if $\gamma$ is a path in $Dom(h) \setm Q_h$ 
with endpoints $z_1$ and $z_2$, then
$$
\| Dh(z_2) \|_m^{1/{\alpha}} \le \| Dh(z_1) \|_m
	\le \| Dh(z_2) \|_m^{\alpha}
$$
where $\alpha = exp(M \ell_m(h(\gamma)))$ for a universal $M > 0$.
\end{prop}
\begin{pf}
We apply McMullen's argument to the approximations of $h$.
Let $\cT_j \in \Tow_0(\kappa)$ converge to $\cT$.
We can assume $S_j = S$ and $m = 0$.
Let $\rho_0$ be the hyperbolic metric on $V_0 \setm P(\cT)$
and let $\rho_{j,0}$ be the hyperbolic metric on $V_{j,0} \setm P(f_{j,0})$.

Since
$$
f_0:(U_0 \setm f_0^{-1}(P(\cT)) \to (V_0 \setm P(\cT))
$$
is a covering map and the inclusion
$$
\imath :(U_0 \setm f_0^{-1}(P(\cT))) \hookrightarrow (V_0 \setm P(\cT))
$$
is a contraction by the Schwarz Lemma, we see $f_0$ expands $\rho_0$.
That is, $\| Df_0(z) \|_0 > 1$ for $z \in (U_0 \setm f_0^{-1}(P(\cT)))$.
Similarly $f_{j,0}$ expands $\rho_{j,0}$.

Suppose $h \in \cF(\cT)$.
Let $z \in (Dom(h) \setm h^{-1}(P(\cT))$.
Choose compact sets $K_1 \subset (Dom(h) \setm P(\cT))$
and $K_2 \subset (Range(h) \setm P(\cT))$
which contain neighborhoods of $z$ and $h(z)$, respectively.
Since the domains $\cl(V_{j,0})$ converge in the
Hausdorff topology to $\cl(V_0)$
and the post-critical sets $P(f_{j,0})$ converge to $P(\cT)$,
the hyperbolic metrics $\rho_{j,0}$ converge uniformly on $K_1$ and $K_2$ 
to $\rho_0$.
For large enough $j$ we have $Range(h) \subset V_{j,0}$
since $Range(h) \subset V_0$. Hence
there are iterates $t_j$ such that 
$$
f^{t_j}_{j,0} \to h
$$
uniformly on $K_1$ in the $C^1$ topology as $j \to \infty$.
Thus maps arbitrarily close to $h$ expand metrics arbitrarily
close to $\rho_0$.
Hence $h$ is non-contracting: $\| Dh(z) \|_0 \ge 1$. 
To prove $h$ is expanding, it suffices to assume $h$ is a local lift
of a transit map. For this we use induction on $n \in S_\cC$. 
First the base case. Let $n = \min S_\cC$
and let $h'$ be another local lift of $g_n$ such that
$$
h = f_0 \circ h'.
$$
Since $f_0$ is expanding and $\wtl{h}$ is non-contracting 
it follows $h$ is expanding and the base case holds.
Now suppose by induction that local lifts of $g_{n_1},\dots,g_{n_k}$
expand $\rho_0$ for the first $k$ levels in $S_\cC$. 
Let $h$ be a local lift of $g_{n_{k+1}}$ where $n_{k+1}$ is the next
level in $S_\cC$ after $n_k$. There is a restriction $f$ of
$f_{n_{k+1}}$ so that $f \in \cF(\cT)$ and we can assume
the attracting and repelling petals $D_\pm$ were chosen to
lie in $Dom(f)$. But then like 
before there is another local lift $h'$ so that
$h = f \circ h'$ 
and we again see $h$ must be expanding.

Now we estimate how much $h$ expands $\rho_0$.
Choose $z \in (Dom(h) \setm Q_h)$ and let $K_1$ and $K_2$ be closed
neighborhoods of $z$ and $h(z)$ as above.
Then just as above for large $j$, we can find iterates $t_j$ such that
$f_{j,0}^{t_j} \to h$ uniformly on $K_1$ as $j \to \infty$.
Let $V^{-n}_j = f_{j,0}^{-n}(V_{j,0})$
and $P^{-n}_j = f_{j,0}^{-n}(P(f_{j,0}))$.
Since
$$f_{j,0}^{t_j}:V^{-t_j}_j \setm P^{-t_j}_j \to 
	V_{j,0} \setm P(f_{j,0})$$
is a local isometry we can apply \lemref{LowExpand} to the inclusion
$$\imath:V^{-t_j}_j \setm P^{-t_j}_j \hookrightarrow 
	V_{j,0} \setm P(f_{j,0})$$
to get the inequalities
$$C_2^{-1}(s) \le \| Df_{j,0}^{t_j}(z) \|_{\rho_{j,0}} \le
	  C_1^{-1}(s)$$
where 
$$s = d_{\rho_{j,0}}(z,P^{-t_j}_j \union \bd V^{-t_j}_j).$$

Since $C_1$ and $C_2$ are increasing,
$$C_2^{-1}(s'_2) \le \| Df_{j,0}^{t_j}(z) \|_{\rho_{j,0}} \le
	  C_1^{-1}(s'_1)$$
where 
$$s'_1 = d_{\rho_{j,0}}(z,P^{-t_j}_j \union \bd K_1)
	\mbox{ and }
	s'_2 = d_{\rho_{j,0}}(z,P^{-t_j}_j \inter K_1).$$

But $f_{j,0}^{t_j} \to h$ uniformly on $K_1$,
$(P^{-t_j}_j \inter K_1) \to (Q_h \inter K_1)$ and
$\rho_{j,0} \to \rho_0$ uniformly on $K_1$ and $K_2$
as $j \to \infty$.
Thus the second statement of the Proposition follows
if we let $K_1$ range over larger and larger compact sets in
$Dom(h) \setm P(\cT)$.

To conclude let us prove the last statement about the variation of
expansion. 
From \cite[Cor 2.27]{McM1} the variation in 
$\| Df_{j,0}^{t_j}(z) \|_0$
is controlled by
the distance between $z_1$ and $z_2$ measured in the hyperbolic metric
on $V^{-t_j}_j \setm P^{-t_j}_j$. Since $f_{j,0}^{t_j}$ is a covering
map, this distance is bounded above by the length of
$f_{j,0}^{t_j}(\gamma)$ measure on $V_{j,0} \setm P(f_{j,0})$.
As $j \to \infty$ this length converges to
$\ell_0(h(\gamma))$. The statement follows.
\end{pf}

The following corollary can be used to control the expansion of the
hyperbolic metric on one level with bounds from a deeper level.

\begin{cor}
\label{expand2}
Let $\cT \in \Tow(\kappa)$ be a forward tower with base map $f_m:U_m \to V_m$.
Suppose $n \in S_Q$ is a level such that $V_n \subset V_m$ and
let $\cT'$ be the tower $\cT$ restricted to the levels $n' \ge n$.
Let $h \in \cF(\cT')$ and let $Q_h = h^{-1}(P(\cT'))$. 
Then if $(Q_h \setm P(\cT')) \neq \es$ and $z \in Dom(h) \setm Q_h$,
$$C_2^{-1}(s_2) \le \| Dh(z) \|_m$$
where $s_2 = d_n(z,Q_h)$.
\end{cor}
\begin{pf}
We may assume $m = 0$.
Since $V_n \subset V_0$ and
$P(\cT') = P(\cT) \inter V_n$ we see
$$(V_n \setm P(\cT')) \subset (V_0 \setm P(\cT))$$
and so
$$d_0(z,Q_h) \le d_n(z,Q_h).$$
Since the function $C_2$ in \propref{expand1} is increasing,
$$C_2^{-1}(d_n(z,Q_h)) \le C_2^{-1}(d_0(z,Q_h)).$$
Finally, since $Range(h) \subset V_n$ and
$V_n \inter P(\cT) = P(\cT')$,
$$h^{-1}(P(\cT')) = h^{-1}(P(\cT)).$$
Since $h \in \cF(\cT)$ it follows from \propref{expand1} that
$$C_2^{-1}(d_0(z,Q_h)) \le \| Dh(z) \|_0.$$
\end{pf}

In order to apply this corollary we need to get a bound on 
$s_2 = d_n(z,Q_h)$. This is done by compactness:

\begin{lem}
\label{BoundToQ}
Let $\cT \in \Tow(\kappa)$, $n \in S_Q$ and
$z \in f_n^{-1}(V_n \setm U_n)$. Then
$d_n(z,Q_{f_n}) \le C(\kappa)$.
\end{lem}
\begin{pf}
By shifting we may assume $n = 0$.
Since $U_0$ and $V_0$ are $\kappa$-quasidisks, the set
$V'_0 = \cl(f_0^{-1}(V_0 \setm U_0))$ varies continuously with 
$\cT \in \Tow(\kappa)$. 
Since $P(\cT)$ varies continously the hyperbolic metric $\rho_0$ 
and the set $Q_{f_0}$ vary continuously.
Therefore the function $F$ on $\Tow(\kappa)$ given by
$$
F(\cT) = \sup_{z \in V'_0} d_0(z,Q_{f_0})
$$
is continuous. Since $\Tow(\kappa)$
is compact by \lemref{TowCompact},
there is a $C(\kappa)$ such that $F(\cT) \le C$.
\end{pf}

\subsection{Equivalent definitions of the Julia set}
\label{sec:Jdef}

There are several equivalent definitions of the Julia set of a rational
map. In this section we present the analogous result for forward towers.

Fix a forward tower $\cT \in \Tow(\kappa)$.
The full orbit of a point under $\cT$,
much like the full orbit of a point very near the origin in the Feigenbaum
map, can be disected to reveal much more structure. For forward towers 
this can be done by iterating deeper maps when possible.

By shifting we may assume $S = \N_0$.
By restricting each $f_n \in \cT$ construct a tower 
$\cT' = \{f'_n,g'_n\}$ such that
\begin{enumerate}
\item $[f_n] = [f'_n]$ and $g_n = g'_n$
\item $V'_{n+1} \subset U'_{n,0}$ for each $n \in S$ except
$V'_{n+1} = V'_n$ for each $n \in S$ such that $f_{n+1}$ is a
first through map
\item $U_0 = U'_0$.
\end{enumerate}
Note that $\cT'$ may no longer be a tower in $\Tow(\kappa)$ but that
$\cF(\cT') = \cF(\cT)$.
For any non-zero $z \in U_0$ define the {\it depth} of $z$ to be 
$$
\depth(z) = \max\{n \in S: z \in U'_{n,0}\}.
$$
For a point $z \in U_0$ we say 
a (possibly finite) sequence
$(z_0,z_1,z_2,...)$ is a {\it sub-orbit} of $z$ (in $\cT'$)
if the following conditions are satisfied:
\begin{itemize}
	\item $z_0 = z$
	\item if $z_i \in V_0 \setm U_0$
		then $z_{i+1}$ is not defined
	\item if $z_i = 0$ then $z_{i+1} = 0$
	\item if $z_i \in Dom(\wtl{g}_n)$ then $z_{i+1} = \wtl{g}_n(z_i)$
		for some local lift $\wtl{g}_n \in \cT$
	\item otherwise $z_{i+1} = f'_{\depth(z_i)}(z_i)$
\end{itemize}
Note any sub-orbit of $z$ is a subset of $\orb(z)$
and $\orb(z)$ escapes iff there exists a sub-orbit $z_0,\dots,z_N$ such that
$z_N \in V_0 \setm U_0$.

A point $z \in U_0$ is called {\it periodic} (in $\cT$) if
there exists $h \in \cF(\cT)$ such that $h(z) = z$. 
Equivalently, $z \neq 0$ is periodic iff there is an $x \in \orb(z)$ such that 
$z \in \orb(x)$ and a sub-orbit $x_0,x_1 = h_1(x),...x_N = h_N(x)$
of $x$ such that $x_0 = x_N$ and $x_0 \neq x_i$ for $0 < i < N$.
The {\it multiplier}, $\lambda$, of the periodic orbit through $z$ is defined to be
$Dh_N(x)$. The multiplier does not depend on the sub-orbit.
A periodic orbit is called 
{\it superattracting, attracting, repelling, neutral} if $\lambda$  satisfies 
$\lambda=0,|\lambda| < 1, |\lambda| > 1, |\lambda| = 1$, respectively.

\begin{lem}
\label{allrepel}
Let $\cT \in \Tow(\kappa)$.
The only non-repelling periodic orbits in $\cT$
are the orbits through the 
parabolic points of $f_n$ for $n \in S_\cC$.
\end{lem}
\begin{pf}
Let $z_0,...,z_N$ be the periodic orbit. 
Since the only non-repelling periodic orbits
in $P(\cT)$ are the orbits through the parabolic
points, we can assume the orbit is disjoint from $P(\cT)$.
By Proposition \ref{expand1}, 
$$\| Dh_N(z) \|_0 > 1$$
But then 
$$|\lambda| = |Dh_N(z)| >1$$
in the euclidean metric as well.
\end{pf}

For a given level $n \in S_\cC$
let $B_n = K(f_n|_{U_{n,0}})$ be the {\it central basin} of level $n$.
A connected compact set $K \subset U_0$ is {\it iterable} if
$K \inter \bd B_n = \es$ for all central basins $B_n$.
Mimicing the definition of sub-orbits of points, we say
a (possibly finite) sequence of compact sets
$(K_0,K_1,K_2,\dots)$ is a {\it sub-orbit} of $K$ 
(in $\cT'$) if the following conditions are satisfied:
\begin{itemize}			
	\item $K_0 = K$
	\item all $K_i$ are iterable except possibly the last one, if it
		exists
	\item if $K_i \subset Dom(\wtl{g}_n)$ then $K_{i+1} = \wtl{g}_n(K_i)$
		for some local lift $\wtl{g}_n \in \cT$ 
	\item otherwise	$K_{i+1} = f'_d(K_i)$ where 
		$d = \min_{z \in K_i} \depth(z)$.
\end{itemize}

Now that we have said what it means to iterate an iterable compact set, 
we can prove the following
\begin{prop}
\label{Jdef}
Suppose $\cT \in \Tow(\kappa)$ and
let $y \in J(\cT)$. The following are two equivalent definitions of
the Julia set:
\begin{enumerate}
	\item $J(\cT) = \cl\{z \in Dom(\cT)
		 : z \mbox{ is a repelling periodic point}\}$
	\item $J(\cT) = \cl\{z \in Dom(\cT) : z
		\mbox{ is a pre-image of } y\}$
\end{enumerate}
\end{prop}
\begin{pf}
By shifting we may assume $S = \N_0$.
We may also assume $S_\cC \neq \es$.
Let $z \in \bd K(\cT)$ and let $W$ be a connected neighborhood of
$z$. We can assume $W \subset \myint(K(f_0))$.
Let $K = \cl(W)$. We can form the suborbit 
$K_i = h_i(K)$ from $K$ until the first moment when $K_i$ is not iterable.
Such a moment must exist since the orbit of $z \in K$ never escapes
but the orbit of some other point in $K$ does escape.

Case 1: Suppose $\myint(K_i) \inter \bd B_n \neq \es$
for some $n \in S$.
Then by arguing as in \lemref{specialaccum}
there is a open set $W' \subset K_i$ 
and composition $h \in \cF(\cT)$ defined on $W'$ such that
$h(W') \inter J(f_0) \neq \es$.
There is then an open set $W'' \subset h(W')$ and an $N \ge 0$ such that
$K(f_0) \subset f_0^N(W'')$.
Since $W \subset K(f_0)$ there exists a point $z_0 \in W$ such that
$$
(f^N_0 \circ h \circ h_i)(z_0) = z_0.
$$
By \lemref{allrepel}, if we chose $W$
to be small enough, $z_0$ must be repelling.

Case 2: If $K_i$ is not iterable because 
$K_i \inter (V_0 \setm U_0) \neq \es$, then
by perhaps choosing a smaller neighborhood $W$ and iterating $f_0$ 
more, we can assume that the moment when $K_i$ is not iterable is
because $\myint(K_i) \inter \bd B_n \neq \es$ for $n = \min S_\cC$
and we can argue as in case 1.

Case 3: Suppose $\myint(K_i) \inter \bd B_n = \es$ 
for some $n \in S_\cC$ but that
$\bd K_i \inter \bd B_n \neq \es$.
Then by choosing a slightly smaller neighborhood $W$ we can assume
$K_i$ is iterable and continue iterating the sub-orbit. We claim this case can
only happen a finite number of times. For otherwise every time $K_i$ is not
iterable $K_i$ falls into this case. Then by choosing the slightly smaller
neighborhoods so that they all contain some definite neighborhood $W'$ of
$z$ we see that the orbit of $W'$ is defined for all iterates. But this is
impossible since then $W'$ never escapes, contradicting the fact that
$z \in \bd K(\cT)$.
Thus after a finite number of restrictions, the non-iterable set $K_i$ must fall
into the cases considered above.
Thus 
$$
J(\cT) \subset \cl\{z \in U_0 :
	 z \mbox{ is a repelling periodic point}\}.
$$

Let $z \in K(\cT)$ and let $W$ be a connected neighborhood of $z$.
Suppose $W$ contains a repelling periodic point $z_0$. Again let 
$K = \cl(W)$ and start forming the sub-orbit $K_i = h_i(K)$
through $K$. Claim there
is a moment when $K_i$ is not iterable. For otherwise the maps $h_i$ form a
normal family on $W$ and that contradicts the fact that $W$ contains a
repelling periodic point. Thus there is a non-iterable iterate $K_i$.

Case 1: Just as case 1 above, there is a open set $W' \subset K_i$ 
and composition $h \in \cF(\cT)$ defined on $W'$ such that
$h(W') \inter J(f_0) \neq \es$.
But then there is a point in $h(W')$ that escapes and thus there is a point
in $W$ that escapes as well.

Case 2: If $K_i$ is not iterable because 
$K_i \inter (V_0 \setm U_0) \neq \es$, then we have found a point
in $W$ that escapes.

Case 3: 
Suppose $\myint(K_i) \inter \bd B_n = \es$ but that
$\bd K_i \inter \bd B_n \neq \es$.
Then by choosing a slightly smaller neighborhood $W$ that still contains the
repelling periodic point $z_0$,
we can assume
$K_i$ is iterable and continue iterating the sub-orbit. We claim this case can
only happen a finite number of times. For otherwise every time $K_i$ is not
iterable $K_i$ falls into this case. Then by choosing the slightly smaller
neighborhoods so that they all contain some definite neighborhood $W'$ 
containing $z_0$
we see that the orbit of $W'$ is defined for all iterates. 
But this is impossible since the iterates of $W'$ cannot form a normal family.
Thus after a finite number of restrictions, the non-iterable set $K_i$ must fall
into the cases considered above.
Thus 
$$
J(\cT) \supset \cl\{z \in U_0 :
	 z \mbox{ is a repelling periodic point}\}.
$$

To prove the second statement, notice that the argument proving
the first also proves that 
if $y \in J(\cT)$ then
any point in $U_0$ has a pre-image arbitrarily close to $y$. That is,
$$J(\cT) \subset \cl\{z \in U_0 : 
	\mbox{ there is an $h$ such that } h(z) = y\}.$$
The reverse inclusion follows from the fact that $J(\cT)$ is closed and
backward invariant and that
$y \in J(\cT)$.
\end{pf}

\subsection{The interior of the filled Julia set}
An infinitely renormalizable quadratic-like map $f \in RQuad$ has a filled
Julia set with empty interior. The same statement holds for forward towers:

\begin{prop}
\label{JEqualK}
For any $\cT \in \Tow(\kappa)$,
$$
\myint(K(\cT)) = \es.
$$
\end{prop}

As a corallary we have
\begin{prop}
\label{JCont}
The Julia set $J(\cT)$ varies continuously with $\cT \in \Tow(\kappa)$.
\end{prop}
\begin{pf}
\end{pf}

The proof of \propref{JEqualK}
is broken into propositions \propref{NoPeriodicComp} and 
\propref{NoWanderingComp} and will occupy the rest of this section.

By shifting we may assume $S = \N_0$.
Suppose by contradiction that 
$$
{\cal O} = Comp(\myint(K(\cT)))
$$ 
is non-empty.
Let $U \in {\cal O}$ and $z \in U$.
Let $K \subset U$ be a compact and connected neighborhood of $z$. 
Recall $B_n$ are the central basins of $\cT$.
Since $\bd B_n \subset J(\cT)$ 
for all $n \in S_\cC$ it follows that $K$ is iterable.
Since $J(\cT)$ is backward invariant we see that 
all the iterates of $K$ are iterable as well. Thus the orbit of $K$
is well defined and contains the orbit of $z$ and so,
letting $K$ range over larger and larger
compact subset of $U$, we can define the orbit of $U$,
$\orb(U)$, to be components containing the orbit of $K$.

A component $U \in {\cal O}$ is called {\it periodic} if
$U' \in \orb(U)$ implies  $U \in \orb(U')$.
A component $U \in {\cal O}$ is called {\it pre-periodic} if
$U$ is not itself periodic but
there is a peridic component in $\orb(U)$.

The classification of periodic components is based on the following

\begin{prop}\cite{L1,M1}
\label{lyub1}
Let $h:U \to U$ be an analytic transform of a hyperbolic Riemann
surface $U$. The we have one of the following possibilities:
\begin{enumerate}
	\item h has an attracting or superattracting fixed point in $U$ to
		which all orbits converge
	\item all orbits tend to infinity
	\item h is conformally conjugate to an irrational rotation of the disk,
		the punctured disk or an annulus
	\item h is a conformal homeomorphism of finite order
\end{enumerate}
\end{prop}

The following proposition expands on case 2) above

\begin{prop}\cite{L1,M1}
\label{lyub2}
Let $U$ be a hyperbolic domain on the sphere, and $h:U \to U$
an analytic transform continuous up to the boundary. Suppose that the set
of fixed points of $h$ on $\bd U$ is totally disconnected. Then in
case 2) of Proposition \ref{lyub1} there is a fixed point $\alpha \in
\bd U$ such that $h_m(z) \to \alpha$ for every
$z \in U$.
\end{prop}

We shall use these two propositions to prove
\begin{prop}
\label{NoPeriodicComp}
No $U \in {\cal O}$ is periodic or pre-periodic.
\end{prop}
\begin{pf}
Suppose $U \in {\cal O}$ is a periodic component.
Suppose $\cl(U)$ is iterable and all iterates of $\cl(U)$ are interable.
Then since $U$ is periodic
there exists a univalent map $h \in \cF(\cT)$
defined on a neighborhood of 
$\cl(U)$ such that $h(\cl(U)) = \cl(U)$.
Let us examine the possiblities from \propref{lyub1}.

Since  $\cl(U)$ is disjoint from $P(\cT)$, \lemref{allrepel} 
implies any periodic point in $\cl(U)$ must be repelling.
Thus there cannot be an superattracting or attracting orbits.
Suppose all iterates tend to $\bd U$.
Now the set of points on $\bd U$ fixed by $h$ are isolated,
since otherwise $h$ would be the identity on an open set
and that would contradict \propref{expand1}.
Applying \propref{lyub2} again contradicts \lemref{allrepel}.

The other possibilities in \propref{lyub1}
are ruled out because $h$ expands the hyperbolic
metric on $U_0 \setm P(\cT)$ and any
map conjugate to a rotation will have high iterates
arbitrarily close to the identity.

Now suppose there is a component $U'$ from $\orb(U)$ 
such that $\cl(U')$ is not iterable.
To simplify the exposition
we will assume $\cT$ is a real-symmetric tower. However, this is not essential.
Since $U$ is periodic we may assume $U = U'$.
Since $U \subset K(f_0)$ there must be an $n \in S_\cC$ 
such that $\cl(U) \inter \bd(B_n) \neq \es$.
Since $U \inter J(\cT) = \es$ it follows that $U \subset B_n$ and
if $n' \in S_\cC$ is the next parabolic level after $n$
then $\cl(U) \inter B_{n'} = \es$. 

Let $K = \cl(U)$, $f = f_n|_{U_{n,0}}$ and $\xi = \beta(f)$.
Since $B_n$ and $\bd B_n$ are invariant by $f$, it follows
$K_k = f^k(K) \subset \cl(B_n) \setm B_{n'}$ and 
$\bd K_k \inter \bd B_n \neq \es$ for all $k \ge 0$.
Let
$$
\cB = Comp(B_n \setm \left(\bigcup_{k \ge 0}^{\infty} f^{-k}(\R)\right))
$$
be the collection of components of the partition pictured in \figref{tile}.
\realfig{tile}{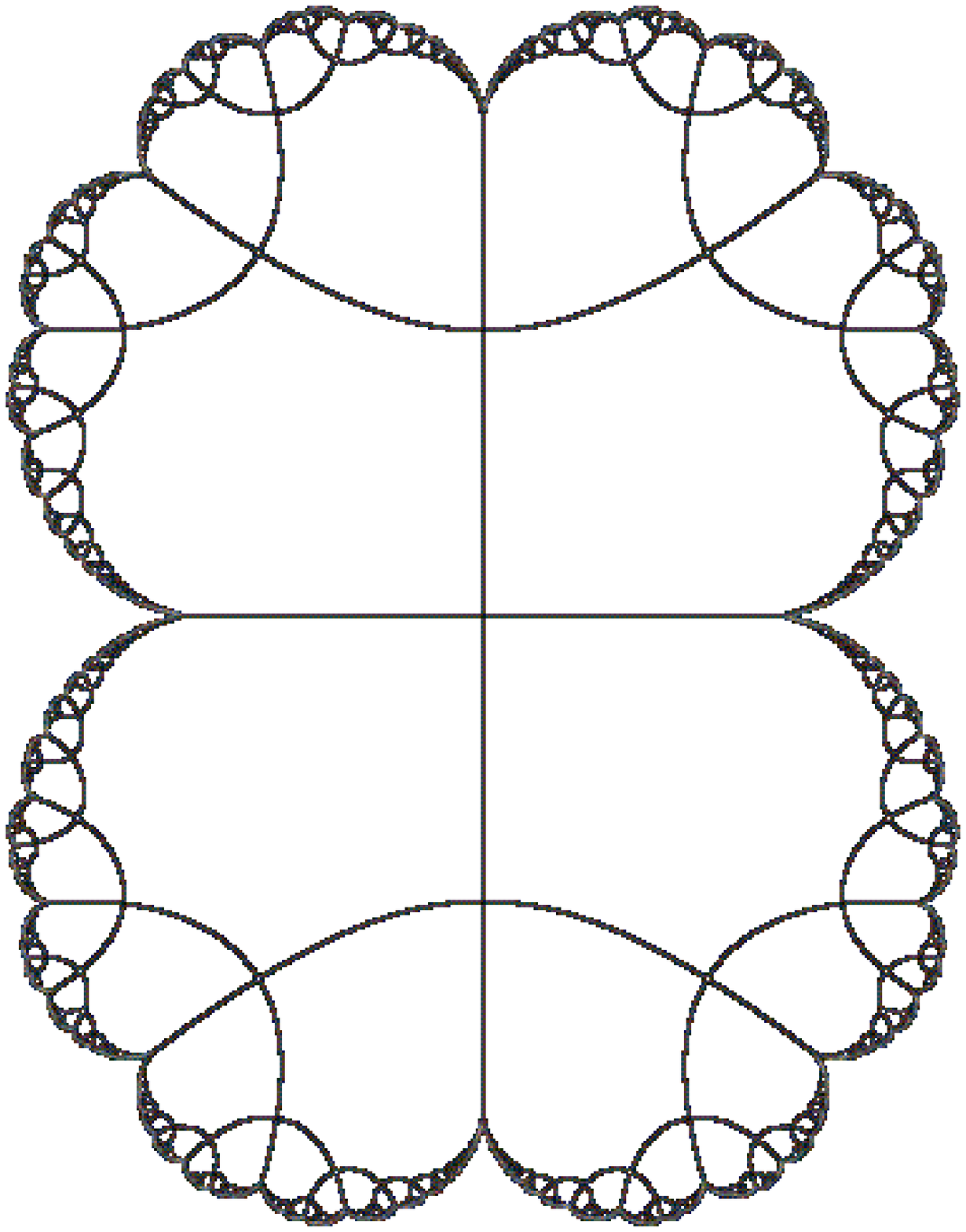}{The tiling of $B_n$.}{0.4\hsize}

First we claim that $U \inter \R = \es$.
Let $n' \in S_Q$ be the largest quadratic-like level before $n$.
Let $B_{n''}$ be the central basin of the first level $n'' \in S_\cC$
after $n'$. Then the $f_{n'}$ pre-images of $B_{n''}$ cover a dense
subset of $\R \inter K(\cT_{n'})$ where $\cT_{n'} \subset \cT$ is the 
forward tower with levels $m \ge n'$.
It follows that the pre-images by $\cF(\cT)$ cover a dense
subset of $\R \inter B_n$ and accumulate at $\xi$.
Since $\bd B_{n''} \subset J(\cT)$ the claim is established.
Since $U$ is periodic under $\cT$, 
we can assume $U \subset A$ where $A \in \cB$
satisfies $\xi \in \bd A$. Without loss of generality assume 
$A \subset \Hp$.

Let $\gamma = \bd A$. Let
$$
\gamma_1 = \bigcup_{\wtl{g}_n} \wtl{g}_n^{-1}(\gamma).
$$
Since $g_n$ is a real translation, $\gamma_1 \subset \Hp$. Let
$$
\gamma_2 = \bigcup_{k \ge 0} (f^{-1})^k(\gamma_1)
$$
where the branch of 
$f^{-1}$ is chosen so that 
$f^{-1}(\Hp \inter B_n) \subset \Hp \inter B_n$
(see \figref{tile2}).
\realfig{tile2}{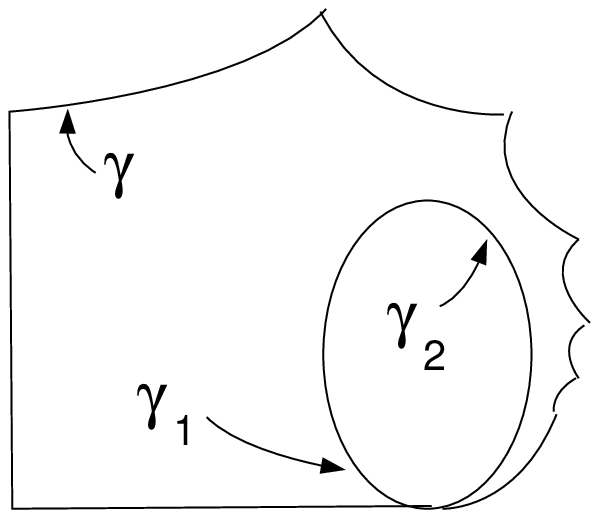}{The curves $\gamma$, $\gamma_1$ and
	 $\gamma_2$.}{0.25\hsize}
It follows from \lemref{specialaccum} that
$U$ is contained in the domain $A_1$
bounded by $\gamma_2$. Continue this process. That is, 
the pre-image of $\gamma_2$ by
$\wtl{g}_n$ is contained in $A_1$ and pulling back by
$f^{-1}$ we see that $U$ is contained in a domain
$A_2 \subset A_1$. By \lemref{EcalleExpands},
$$
\bigcap_{m \ge 1}^\infty A_m = \es
$$
and so a non-iterable periodic component $U$ cannot exist.
\end{pf}

A component $U \in {\cal O}$ that is neither periodic nor
pre-periodic is called {\it wandering}.

\begin{prop}
\label{NoWanderingComp}
No $U \in {\cal O}$ is wandering.
\end{prop}
\begin{pf}
Suppose $U \in {\cal O}$ is wandering.
Let $K \subset U$ be compact and connected. Then $K$ is iterable and
all iterates of $K$ are iterable.
Fix an $z \in \myint(K)$.
Since each map $h$ from the orbit of $K$ is defined on a neighborhood of $K$
and since $Q_{h} = h^{-1}(P(\cT)) \subset J(\cT)$, it follows from
Proposition \ref{expand1} that
\begin{equation}
\label{eqn:boundedExpand}
\sup_h \| Dh(z) \|_0 < \infty.
\end{equation}

Suppose there is an $\eps > 0$ such that there are 
an infinite number of iterates
$h_n$ satisfying
$$
d(h_n(z),P(\cT)) > \eps
$$
where the distance is just the euclidean distance. Order the $h_n$ to match
the ordering on the orbit. That is, if $n < m$ then $h_m(z) \in
\orb(h_n(z))$. Since each $h_n(z)$ lies in a compact subset of the hyperbolic
surface $V_0 \setm P(\cT)$,
$$
d_0(h_n(z),Q_{f_0}) \le C_\eps,
$$
and so from \propref{expand1},
\begin{equation}
\label{eqn:boundedExpand2}
\| Df_0(h_n(z)) \|_0 \ge C > 1.
\end{equation}
But then 
\begin{equation}
\label{eqn:boundedExpand3}
  \| Dh_{n+1}(z) \|_0
  \ge \| D(f_0 \circ h_n)(z) \|_0
  = \| Df_0(h_n(z)) \|_0 \cdot
  \| Dh_n(z) \|_0
  \ge C \| Dh_n(z) \|_0
\end{equation}
which as $n \to \infty$
contradicts equation \ref{eqn:boundedExpand}.

So we can assume
$$
\limsup_h d(h(z),P(\cT)) = 0.
$$
Let $\cT_n = \cT|_{S_n}$ be the restriction of $\cT$ to levels $m \ge n$.
Let $\cK_n$ be the collection of little filled Julia sets
$\cK_n = \orb(\cT,K(\cT_n))$.

From \lemref{EcalleExpands} we see $\orb(z)$
must accumulate on some
$z' \not\in \xi_0$ where $\xi_0$ is the parabolic orbit of $f_0$. But then
$z'$ is contained in a little filled Julia set in $\cK_1$.
By iterating forward we can assume $z' \in K(\cT_1)$.
It follows that there is a $y_1 \in \orb(z)$ such that $y_1 \in K(\cT_1)$.
Now again there is an accumulation point $\orb(y_1)$ disjoint from $\xi_1$,
the parabolic orbit of $f_1$, and,
repeating the whole argument inductively, there is
a sequence of iterates $y_n \in K(\cT_n)$.

Each $y_n$ has a moment $x_n \in \orb(z)$
when $\orb(z)$ enters the collection of little filled Julia sets 
$\cK_n$. Once $\orb(z)$ enters $\cK_n$ it never leaves.
It can happen that different $y_n$ have the same moment $x_n$. However,
since 
$$
\bigcap_{n \ge 0} K(\cT_n) = \{0\}
$$
there must be an infinite number of distinct entry moments $x_{n_i}$.

Let $z_n \in \orb(z)$ satisfy $f_0(z_n) = x_n$. Thus the relation between the
points $z$, $x_n$, $y_n$ and $z_n$ is given by:
$z_n \in \orb(z)$, $x_n = f_0(z_n)$ is the time $\orb(z)$ enters ${\cal K}_n$
and $y_n \in \orb(x_n)$ is the first time $x_n$ enters $K(\cT_n)$.
Claim
$$
d_0(z_n,Q_{f_0}) \le C'.
$$

Let $K'_n$ be the component of
$f_0^{-1}({\cal K}_n) \setm {\cal K}_n$
containing $z_n$.
The set $K'_n$ is called a {\it companion} filled Julia set of level $n$.
Since $Q_{f_0} \inter K'_n \neq \es$, it is enough to show
$$
\diam_0(K'_n) \le C'.
$$
Consider the sets $U'_n$ and $V'_n$ containing $K'_n$ which are 
pull-backs of $Dom(f_n)$ and $V_n$ by the map sending
$z_n$ to $y_n$. 
By the unbranched property this pull-back is univalent. Since 
$\mymod(Dom(f_n),V_n) \ge 1/\kappa$, we have
$\mymod(U'_n,V'_n) \ge 1/\kappa$ and so, from \cite[Theorem 2.4]{McM1},
the diameter $D_n$ of $U'_n$ in the hyperbolic metric on $V'_n$ is bounded.
But $V'_n \subset (V_0 \setm P(\cT))$. Thus
$$
\diam_0(K'_n) \le \diam_0(U'_n) \le D_n \le C(\kappa)
$$
and the claim is established.

But then equations \ref{eqn:boundedExpand2} and \ref{eqn:boundedExpand3} 
hold along the sequence $z_{n_i} = h_{n_i}(z)$,
and we again get a contradiction to \ref{eqn:boundedExpand}.
\end{pf}

\subsection{Line fields and forward towers}
\label{sec:forrigid}

A {\it line field} is a measurable Beltrami differential with
$|u(z)| = 1$ on a set of positive measure and 
$|u(z)| = 0$ otherwise. A line field is {\it invariant} under $\cT$ iff
for every $h \in \cT$, $Dh$ maps the line at $x$ to the line at $h(x)$
for almost every $x \in Dom(h)$.
Using \propref{Jdef} and \propref{JEqualK}
we can rephrase \propref{genCombRigid} in terms
of invariant line fields.

Before doing so, we need the following
\begin{lem}\cite{L1}
\label{totalDiscont}
Let $\cT \in \Tow(\kappa)$.
The group $G$ of homeomorphisms of $J(\cT)$
that commute with all maps $h \in \cT$ is totally disconnected.
\end{lem}
\begin{pf}
Let $\phi \in G$ be a map in the connected component of the identity.
Suppose $z_0$ is a repelling periodic point with $h(z_0) = z_0$
for some $h \in \cF(\cT)$.
Since the solutions to $h(z) = z$ are isolated $\phi$ must fix $z_0$.
The lemma follows from density of repelling cycles: \propref{Jdef}.
\end{pf}

We now prove the following version of forward tower rigidity:

\begin{prop}[No Line Fields for Forward Towers]
\label{forRigid}
No forward tower $\cT$ hybrid equivalent to a tower in $\Tow(\kappa)$
supports an invariant line field on its filled Julia set.
\end{prop}
\begin{pf}
By \propref{forStraight} it suffices to consider a forward tower $\cT$
having a base map of the form $z^2+c_0$ and by
shifting we may assume $S = \N_0$. Since $\cT$ is hybrid equivalent to
a tower in $\Tow(\kappa)$ it follows from
\propref{JEqualK} that $K(\cT) = J(\cT)$.
Suppose by contradiction that $\cT$ did admit an invariant line field 
$$\mu = u(z) d \bar{z}/dz$$
supported on $J(\cT)$. 
For any $w \in \D$ consider the invariant Beltrami differential
$$
\mu_w = w \cdot  u(z) d\bar{z}/dz
$$ on $\Chat$.
Let $\phi_w$ be a solution to the corresponding Beltrami equation normalized
so that the map 
$$
f_{w,0} = \phi_w \circ f_0 \circ \phi_w^{-1}
$$
is again a rational map of the form $z^2+c_w$ for some $c_w \in \C$.
Let $\cT_w$ be the tower 
$$
\cT_w = \{\phi_w \circ h \circ \phi_w^{-1}: h \in \cT\}.
$$

From \propref{hyper-dense} and the uniqueness of root points,
$c_w = c_0$ for all $w \in \D$.
\propref{genCombRigid} implies $f_{w,n} = f_n$ for all $n \in S$
and $g_{w,n} = g_n$ for all $n \in S_\cC$ and
$w \in \D$. So $\phi_w$ is a holomorphic family of quasi-conformal
maps with $\phi_0 = id$ and $\phi_w$
mapping $J(\cT)$ homeomorphically to itself commuting
with the dynamics of $\cT$. From \lemref{totalDiscont}
$\phi_w|_{J(\cT)} = id$.
But then the complex dilatation of $\phi_w$ is zero at
all points of Lebesgue density of
$J(\cT)$ and so $\mu$ is not supported on $J(\cT)$, a 
contradiction.
\end{pf}

\subsection {Line fields and bi-infinite towers}
\label{sec:biTowers}

In this section we move from studying properties of
forward towers to studying bi-infinite towers. 
The plan of attack again follows \cite{McM2}.

Let $\Tow^\infty(\kappa)$ denote the set of 
bi-infinite towers in $\Tow(\kappa)$. 
Given $\cT \in \Tow^{\infty}(\kappa)$ define 
$S_\cN \subset S_Q$ as follows. Let $S_{\cN,0} = \{0\}$.
Then inductively 
let $S_{\cN,n+1} = S_{\cN,n} \union \{m_{n+1}\}$
where $m_{n+1} = \max \{m \in S_Q | m < m', U_m \supset V_{m'}\}$
and $m' = \min S_{\cN,n}$.
Define $S_\cN = \union_{n \to \infty} S_{\cN,n}$.
That is, $S_\cN$ is the minimal set of nested levels
approaching $-\infty$.
From \lemref{scaling} we see
$S_\cN$ is unbounded below.

Define the depth of a non-zero point $z \in \C$ by
$$
\depth(z) = \max\{m \in S_\cN: z \in U_m\}.
$$
For a point $z \in \C$ we say 
a (possibly finite) sequence
$(z_0,z_1,z_2,...)$ is a {\it sub-orbit} of $z$ (in $\cT$)
if the following conditions are satisfied:
\begin{itemize}
	\item $z_0 = z$
	\item if $z_i = 0$ then $z_{i+1} = 0$
	\item if $z_i \in Dom(\wtl{g}_n)$ then $z_{i+1} = \wtl{g}_{n}(z_i)$
		for some $\wtl{g}_{n} \in \cT$
	\item otherwise $z_{i+1} = f_{\depth(z_i)}(z_i)$
\end{itemize}
%Any sub-orbit $(z_0,z_1,\dots)$ such that 
%$\union_i z_i \in U_n$ for some $n \in S_\cN$ 
%is also a sub-orbit of the forward tower $\cT|_{S_n}$.

Let $\rho_{-\infty}$ be the hyperbolic metric on
$\C \setm P(\cT)$ and as in \secref{sec:ForwardExpands} let 
$\rho_n$ be the hyperbolic metric on
$V_n \setm P(\cT|_{S_n})$.
From \lemref{scaling} and the unbranched property
the metrics $\rho_n$ converge uniformly on compact sets to 
$\rho_{-\infty}$. Using the expansion from 
\secref{sec:ForwardExpands}, we now prove

$\bold{Theorem}$ $\bold{\ref{dense-lemma}}$. {\it
For any $\cT \in \Tow^{-\infty}(\kappa)$
$$
\lim_{n \to -\infty} J(\cT|_{S_n}) = \Chat
$$
in the Hausdorff topology.}
\begin{pf}
Let $\cT_n = \cT|_{S_n}$.
Let $z \not\in \union_{s \le 0} J(\cT_s)$.
Without loss of generality we may assume $z \in U_0$.
Then $\orb(\cT_s,z)$ escapes $U_s$
for any $s \in S_\cN$. 
Let $z_s = h_s(z)$ be the orbit point just 
before the first moment of escape on level
$s$. That is, $f_s(z_s) \in V_s \setm U_s$ and if $z' \in \orb(z)$ also
satisfies $f_s(z') \in V_s \setm U_s$ then $z' \in \orb(z_s)$.
For a given $s \in S_\cN$ let $\gamma'_s$ be a hyperbolic geodesic in
$V_s \setm P(\cT_s)$ connecting $z_s$ with $J(\cT_s)$. 
From \lemref{BoundToQ}, there is a $C$ independent of $s$
such that $\ell_s(\gamma'_s) \le C.$
Fix a small $\eps > 0$ and let $A$ be an $\eps$-scaled neighborhood 
of $P(\cT_s)$.
Then $h_s$ has an extension $h \in \cF(\cT)$ that is a 
covering map onto $V_s \setm A$.
Let $\gamma_s$ be the connected component of $h^{-1}(\gamma'_s)$
containing $z$.

We now argue $\ell_s(\gamma_s)$ shrinks as
$s \to -\infty$. The proposition would follow since $\rho_s$ converges to
$\rho_{-\infty}$ near $z$ and since Julia sets are backward invariant.
Fix an $s \in S_\cN$ and let $N_s = |\{s,\dots,0\} \inter S_\cN |$ 
be the minimal number of moments when the orbit of $z$ escapes a nested
level. 
It follows from \lemref{BoundToQ} and \corref{expand2} that there is a 
$C > 1$ such that
$$
C \le \| Df_t(z_t) \|_s
$$
for any $t \in \{s,\dots,0\} \inter S_\cN$.
Hence 
\begin{equation}
\label{eqn:ExpGrow}
C^{N_s} \le \| Dh_s(z) \|_s.
\end{equation}
Hence the derivative at the endpoint $z$ grows exponentially in $N_s$. From
\propref{expand1}, there exists a $C > 1$ such that
equation \ref{eqn:ExpGrow} holds along $\gamma_s$
and hence the length of $\gamma_s$ shrinks as $s \to -\infty$.
\end{pf}

A measurable line field $\mu$ on an open set $U$ is called
{\it univalent} if there is a univalent
map $h: U \to \C$ such that $\mu = h^*(d{\bar z}/dz)$.
The main statement in this section is the following extension of 
\propref{forRigid}.

\begin{thm}[No Line Fields for Bi-infinite Towers]
\label{bi-rigidity}
Let $\cT \in Tow_{\infty}(\kappa)$.
There does not exist a measurable line field $\mu$
in the plane such that $h_*(\mu) = \mu$ for all 
$h \in \cF(\cT_n)$, $n \in S_\cN$.
\end{thm}
\begin{pf}
Suppose to the contrary that $\mu = u(z) d{\bar z}/dz$
is a measurable invariant line field which
is non-zero on a set, $B$, of positive measure.
Let $z \in B$ be a point of almost continuity of $u$ and satisfying
$|u(z)| = 1$.
That is, for each $\eps > 0$, the chance of randomly
choosing a point $y$ a distance
$r$ from $z$ that satisfies $|u(y) - u(z)| > \eps$ tends to $0$ as $r$
tends to $0$:
$$
\lim_{r \to 0} \frac{\area(\{y \in B(z,r)\ : |u(y)-u(z)| > \eps\})}
	{\area B(z,r)} = 0
$$
where $B(z,r)$ is the euclidean ball of radius $r$ centered at $z$.
By \propref{forRigid}, we can assume $z \notin K(\cT_n)$ for any $n$.
Let $z_n$ be an infinite sub-orbit from $z$ and 
for each $s \in S_\cN$ let $z_{n_s} = h_{n_s}(z)$
denote the moments in the sub-orbit when $z_{{n_s}+1}$ first satisfies
$z_{{n_s}+1} \in V_s \setm U_s$.

For a given $s \in S_\cN$ let $\cT_s$ denote $\cT$ shifted so that
level $s$ is moved to level $0$ and let
$w_s$ and $u_s$ denote $z_{n_s}$ and $u$ shifted by $s$.
That is, if $|B(f_s)| = \alpha_s$, then
$w_s = \alpha^{-1}_s z_{n_s}$ and $u_s(z) = u(\alpha_s z)$.
Then since 
$\Tow^{\infty}(\kappa)$ 
is compact the sequence $\cT_s$ has a subsequence which as
$s \to -\infty$ converges to some $\cT' \in \Tow^\infty(\kappa)$.
By choosing a further subsequence we may assume 
$w_s$ converges to a $w \in \cl((f'_0)^{-1}(V'_0 \setm U'_0))$ and,
from \cite{McM2},
$\mu_s$ converges $\mbox{weak}^*$, and hence
pointwise almost everywhere, to a measurable line field 
$\mu'$ invariant by $\cT'$
in the sense that $h_*(\mu) = \mu$ for all
$h \in \cF(\cT'_n)$, $n \in S_\cN(\cT')$.

Let $D$ be a small disk around $w$ in $V'_0 \setm P(\cT'_0)$. The
hyperbolic diameter of $D$ in $V'_0 \setm P(\cT'_0)$ is close to that of
$D_s = \alpha^{-1}_s(D)$ in the metric on $V_s \setm P(\cT_s)$ for $s$
near $-\infty$. Since $D_s$ is disjoint from $P(\cT_s)$, there is, by the
argument given in \propref{dense-lemma}, a
univalent pullback $D'_s$ of $D_s$ by the map $h_{n_s}$. By equation
\ref{eqn:ExpGrow} and the variation of expansion in \propref{expand1},
we see $D'_s$ is a sequence of open sets containing $z$ such that
in the euclidean metric
$\diam(D'_s) \to 0$ and $B(z,C \diam(D'_s)) \subset D'_s$ 
as $s \to -\infty$ for some constant $C$. 
Therefore from \cite[Theorem 5.16]{McM1}
we can choose $\mu'$ to be univalent on $D$.

By \propref{dense-lemma}, there is an $s \in S_\cN(\cT')$ such that
$J(\cT'_s) \inter D \neq \es$.
By invariance, if $Dh(z) \neq 0$ and $\mu'$ is locally univalent around $z$
then $\mu'$ agrees almost everywhere with a locally univalent 
line field around $h(z)$ for any composition $h \in \cF(\cT'_s)$.
From \propref{Jdef}, the orbit of $D$ by $\cT'_s$ covers all of $V'_s$.
So $\mu'$ agrees almost everywhere with a line field that is locally
univalent on the set $V'_s \setm P(\cT'_s)$. Since $f'_s$ is
injective on $P(\cT'_s)$ every point in $P(\cT'_s)$ except 
$f'_s(0)$ has an $f'_s$ pre-image around which $\mu'$ agrees (a.e.) with a
locally univalent line field. Hence $\mu'$ agrees (a.e.) with a 
locally univalent line field around $(f'_s)^2(0)$ and $0$, which is a
contradiction, since then we obtain contradictory behavior of $\mu'$
around $f'_s(0)$.
\end{pf}

As a corollary we obtain

%\begin{cor}[Bi-infinite Tower Rigidity]
$\bold{Theorem}$ $\bold{\ref{affine}}$. {\it
If $\cT,\cT' \in \Tow^{\infty}(\kappa)$ are
normalized combinatorially equivalent towers then
$[f_n] = [f'_n]$ for all $n \in S$ and $g_n = g'_n$ for $n \in S_\cC$.
}
\begin{pf}
Let $S_\cN$ be the set of nested levels of $\cT$ as constructed above.
For each $n \in S_\cN$,
let $\phi_n$ be a hybrid equivalence between
$\cT_n$ and $\cT'_n$ coming from straightening 
(see \corref{CombIsHybrid}). The dilatation of $\phi_n$ is bounded 
above by a constant depending only on
$\kappa$ and $\phi_n$ fixes $0$ and $\infty$ and maps
$\beta(f_0)$ to $\beta(f'_0)$. Thus we can pass to a convergent
subsequence $\phi_{n_k} \to \phi$ as $n \to -\infty$.
Since $\phi_n$ restricts to a quasi-conformal equivalence
of $f_s$ and $f'_s$ for $s > n$, $s \in S_\cN$,
on a definite neighborhood of 
$K(f_s)$, it follows that $\phi$ is a quasi-conformal equivalence.
Let $\mu$ be the line field defined by $\phi$ and $\mu_n$ the line field
defined by $\phi_n$.
Since $h_*(\mu_n) = \mu_n$ for all $h \in \cF(\cT_n)$ 
it follows that 
$h_*(\mu) = \mu$ for all $h \in \cF(\cT_n)$, $n \in S_\cN$.

%There are restrictions $\wtl{f}'_n$
%of $\wtl{f}_n$ for $n \in {\Bbb Z}$ such that
%$\seq{\wtl{f}',\sigma/L} \in \T^{-\infty}_{\infty}(\bar{p})$
%and such that
%$\mu$ is invariant under any tower $\wtl{\F'}$ built from
%$\seq{\wtl{f}',\sigma/L}$.
From \thmref{bi-rigidity}, $\mu = 0$ and so $\phi$ is conformal.
Since $\lim_{n \to -\infty} U_n = \Chat$, $\phi$ is linear and
since $\cT$ and $\cT'$ are normalized $\phi$ is the identity.
\end{pf}

\section{Proof of \thmref{per3} and \thmref{generalThm}}

Let $p > 1$.
Let $f$ be an $\infty$-renormalizable real quadratic polynomial 
with $\bar{p}_e(f) \le p$.
The first step in the proof of \thmref{generalThm} is to construct
a tower $\cT \in \Tow(\kappa)$ from $f$.

It follows from \lemref{unbranched}
that there is a $\kappa > 0$ depending only on $p$, and
a forward tower $\cT = \{f_n\} \in \Tow_0(\kappa)$ with the 
following property. For $n \in S_Q$ let $[f'_n]$ be $[f_n]$ normalized and
let $k(n) = | S_Q \inter \{1,\dots,n\} |$.
Then
$$
[f'_n] = \cR^{k(n)}(f).
$$
Hence renormalization acts on towers by shifting. Let $\cT_n$ denote the
tower $\cT$ shifted by $n$ so that $f_n$ is normalized and has index $0$.
By compactness there exists a limiting tower $\cT'$ 
and by \thmref{affine} the germ $[f_0]$ is
uniquely specified by the combinatorics of $\cT'$: a
bi-infinite sequence of $\sigma \in \Omega^{cpt}_p$.
Hence if $f$ has essentially period tripling combinatorics the
germs $\cR^k(f)$ converge to a unique germ $F$, which proves
\thmref{per3}.

To prove \thmref{generalThm} suppose $\bar{\sigma} \in \Sigma_p$ is a 
bi-infinite sequence of shuffles and ends in $\Omega^{cpt}_p$.
Let $\sigma_n = \pi_n(\bar{\sigma})$ denote the $n$-th element of $\bar{\sigma}$.
For each $\sigma_n$ let $\sigma_{m,n}$ be a sequence in $\Omega_p$
converging to $\sigma_n$.
Define the sequence $\bar(\tau) \in \Pi_0^{\infty} \Omega_p$ by
$$
\bar{\tau} = (\sigma_{0,0},\sigma_{1,-1},\sigma_{1,0},\sigma_{1,1},
\sigma_{2,-2},\sigma_{2,-1},\sigma_{2,0},\sigma_{2,1},\sigma_{2,2},
\dots,\sigma_{n,-n},\dots,\sigma_{n,n},\dots)
$$
and let $\bar{\tau}_n  = \omega^{j(n)}(\bar{tau})$ 
where $\omega$ is the left-shift 
operator and $j(n) = 1+3+5+\cdots+(2n-1)+n$. 
Then by construction $\bar{\tau}_n \to \bar{\sigma}$.
Let $f$ be a real quadratic polynomial with shuffle sequence $\bar{\tau}$
and let $\cT$ be a tower in $\Tow(\kappa)$ constructed from $f$.
By compactness of towers 
let $\cT' = \{f'_n,g'_n\}$ be a limiting tower of $\cT_{j(n)}$.
Define the function $h:\Sigma_p \to GQuad(m)$ by
$$
h(\bar{\sigma}) = [f'_0].
$$
From \thmref{affine} $h$ is well-defined and is continuous.
The other properties of $h$ are clear.

\end{document}